\documentclass[11pt,reqno]{amsart}
\usepackage[T1]{fontenc}
\usepackage{lmodern}
\usepackage{microtype}
\usepackage[a4paper,textwidth=150mm,textheight=229mm,centering]{geometry}
\usepackage{amsmath,amssymb,amsthm,mathtools}
\usepackage{booktabs,array,tabularx}
\usepackage{enumitem}
\usepackage[hidelinks]{hyperref}
\hypersetup{pdftitle={Nonlocal Fisher--KPP propagation with heavy tails},pdfauthor={Yuanyang Hu}}
\newtheorem{theorem}{Theorem}[section]
\newtheorem{lemma}[theorem]{Lemma}
\newtheorem{proposition}[theorem]{Proposition}
\newtheorem{corollary}[theorem]{Corollary}
\theoremstyle{definition}

\theoremstyle{remark}
\newtheorem{remark}[theorem]{Remark}
\numberwithin{equation}{section}
\newcommand{\R}{\mathbb R}
\newcommand{\N}{\mathbb N}
\newcommand{\ee}{\mathrm e}
\newcommand{\dd}{\,\mathrm d}
\newcommand{\ind}{\mathbf 1}
\newcommand{\Prob}{\mathbb P}
\newcommand{\E}{\mathbb E}
\newcommand{\Lop}{\mathcal L}
\newcommand{\Yg}{\mathcal Y_\gamma}
\newcommand{\Xg}{\mathcal X_\gamma}
\DeclareMathOperator{\supp}{supp}
\DeclareMathOperator{\Lip}{Lip}
\setlist[enumerate]{label=\textup{(\roman*)},leftmargin=*,itemsep=3pt}
\allowdisplaybreaks[2]
\title[Nonlocal Fisher--KPP propagation]{Nonlocal Fisher--KPP propagation with heavy tails}
\author{Yuanyang Hu}

\makeatletter
\let\old@setauthors\@setauthors
\renewcommand{\@setauthors}{%
  \old@setauthors
  \vspace{-4pt}%
  \begin{center}
    \small School of Mathematics and Statistics, Henan University\\
    475004 Kaifeng, China\\
    \texttt{yuanyhu@henu.edu.cn}
  \end{center}%
}
\makeatother

\date{}
\begin{document}
\begin{abstract}
We study a nonlocal Fisher--KPP equation on the real line with compactly supported initial data. The reaction satisfies $f(s)\le f^{\prime}(0)s$ for $0\le s\le1$, with a power bound on $f^{\prime}(0)s-f(s)$ near zero. For continuous even probability kernels $J$ comparable to $|x|^{-\gamma}$, $\gamma>1$, we prove a global two-sided estimate and locate every fixed level set with a bounded error in $\log|x|$. The propagation scale is $\exp(f'(0)t/\gamma)$ for every $\gamma>1$. This differs sharply from the free-boundary propagation studied by Du and Ni \cite{DuNiJEMS}. For polynomial kernel tails, spreading in their model has linear, algebraic or $t\log t$ growth. We also obtain arrival-time and weighted-integral estimates and identify a spatial tail coefficient when $|x|^\gamma J(x)$ has a positive limit.
\end{abstract}
\maketitle
\section{Introduction}\label{sec:intro}

We study the location of fixed level sets of a Fisher--KPP solution with a heavy-tailed dispersal kernel. The solution satisfies the Cauchy problem
\begin{equation}\label{eq:cauchy}
\begin{cases}
 u_t(t,x)=d\left[\displaystyle\int_{\R}J(x-y)u(t,y)\dd y-u(t,x)\right]
              +f(u(t,x)),& t>0,\ x\in\R,\\
 u(0,x)=u_0(x),&x\in\R,
\end{cases}
\end{equation}
where $d>0$ is the dispersal rate. For each fixed level, we seek bounds on the logarithmic positions of all points in its level set, with an error independent of time. We also ask how the finiteness of $\int_{\R}|z|J(z)\dd z$ and $\int_{\R}z^2J(z)\dd z$ enters the level-set estimates.

Throughout the paper the basic kernel assumption is
\begin{equation}\label{eq:J}
 J\in C(\R)\cap L^\infty(\R),\qquad J\ge0,\qquad J(0)>0,
 \qquad \int_{\R}J(z)\dd z=1,\qquad J(-z)=J(z).
\end{equation}
For $p>0$, define $m_p(J)=\int_{\R}|z|^pJ(z)\dd z\in[0,\infty]$. We consider heavy tails with two-sided power-law bounds:
\begin{equation}\label{eq:tail}
 J(x)\sim |x|^{-\gamma}\quad\text{for }|x|\gg1,\qquad \gamma>1.
\end{equation}
For positive quantities $a$ and $b$, the notation $a\sim b$ means that $cb\le a\le Cb$ for constants $c,C>0$ independent of the variables in the stated range. Thus \eqref{eq:tail} means that there are $j_1,j_2,R_J>0$ such that
\begin{equation}\label{eq:tail-explicit}
 j_1|x|^{-\gamma}\le J(x)\le j_2|x|^{-\gamma}
 \text{ for } |x|\ge R_J.
\end{equation}
The bounds in \eqref{eq:tail-explicit} allow $J(x)|x|^\gamma$ to oscillate. Additional kernel hypotheses, including an actual tail limit, are stated with the results that use them.

For a bounded Borel function $v:\R\to\R$, convolution by $J$ means $(J*v)(x)=\int_{\R}J(z)v(x-z)\dd z$. Here $C_c(\R)$ denotes the continuous functions with compact support. The initial function satisfies
\begin{equation}\label{eq:data}
 u_0\in C_c(\R),\qquad 0\le u_0\le1,\qquad u_0\not\equiv0.
\end{equation}
For the reaction, we assume
\begin{equation}\label{eq:F}
\begin{gathered}
 f\in C^1([0,1]),\qquad f(0)=f(1)=0,\qquad f(s)>0\text{ for } 0<s<1,\\
 r:=f'(0)>0,\qquad
 0\le rs-f(s)\le K s^{1+\sigma}\text{ for } 0\le s\le1
\end{gathered}
\end{equation}
for some $K>0$ and $0<\sigma\le1$. The bound $0\le rs-f(s)\le Ks^{1+\sigma}$ holds, for example, when $f(s)\le rs$ on $[0,1]$ and $f$ is $C^{1,\sigma}$ near zero.

We denote the solution of \eqref{eq:cauchy} by $u(t,x)$. For $0<\theta<1$, define
\begin{equation}\label{eq:level-def}
 E_\theta(t)=\{x\in\R:u(t,x)=\theta\},\qquad
 \Omega_\theta(t)=\{x\in\R:u(t,x)\ge\theta\}.
\end{equation}
Whenever $E_\theta(t)$ is nonempty and compact, set $x_\theta^-(t)=\min E_\theta(t)$ and $x_\theta^+(t)=\max E_\theta(t)$. Theorem~\ref{thm:main}\textup{(ii)} establishes the nonemptiness and compactness of $E_\theta(t)$ for all sufficiently large times. In the notation $O_\theta(1)$, the bound may depend on $\theta$ and is uniform in the other variables in the stated range.

\subsection{Previous results}

The Fisher--KPP equation was introduced in \cite{Fisher1937,KPP1937}, and its spreading theory was developed by Aronson and Weinberger \cite{AronsonWeinberger1975,AronsonWeinberger1978}. Schumacher \cite{Schumacher1980} studied travelling fronts for integral dispersal, while Weinberger \cite{Weinberger1982} developed spreading theory for order-preserving recursions. The role of dispersal tails was emphasized by Kot, Lewis and van den Driessche \cite{KotLewisDriessche1996} and Medlock and Kot \cite{MedlockKot2003}. Subsequent results include asymmetric dispersal \cite{CovilleDavilaMartinez2008,XuLiRuan2021}, probability-measure kernels \cite{Yagisita2009}, and spreading formulas for recursion models \cite{WeinbergerZhao2010}.

For the continuous-time convolution equation, Garnier \cite[Theorems~1--3]{Garnier2011} proved acceleration and quantitative level-set bounds under hypotheses that include $m_1(J)<\infty$ and $J^{\prime}(x)/J(x)\to0$ as $x\to+\infty$. His polynomial examples have exponent greater than two. Hamel and Roques \cite{HamelRoques2010} studied acceleration caused by slowly decaying initial data for local diffusion. In both settings, reaction amplifies a small spatial tail.

Finkelshtein and Tkachov \cite{FinkelshteinTkachov2019} treated nonsymmetric heavy tails and general monostable reactions. Their Theorem~1 places the transition between tail scales evaluated at $(1-\varepsilon)t$ and $(1+\varepsilon)t$ for every fixed $\varepsilon>0$. The hypotheses permit oscillation between suitable envelopes with equivalent logarithms. Finkelshtein, Kondratiev and Tkachov \cite{FinkelshteinEtAl2019} obtained multidimensional propagation results.

For fractional diffusion with symbol $|\xi|^{2s}$ in $\R^n$, $0<s<1$, Cabr\'e and Roquejoffre \cite{CabreRoquejoffre2009,CabreRoquejoffre2013} established the compact-data exponential rate $f'(0)/(n+2s)$. Their Theorem~1.6 in \cite{CabreRoquejoffre2013} gives constant-factor localization for the one-dimensional half-Laplacian with logistic reaction. Cabr\'e, Coulon and Roquejoffre \cite{CabreCoulonRoquejoffre2012} proved constant-factor bounds in periodic logistic media. Roquejoffre and Tarfulea \cite[Theorem~1.3 and Lemma~5.1]{RoquejoffreTarfulea2017} obtained a limiting rescaled radius and fixed-time tail asymptotics under their hypotheses.

Long-time and large-space rescaling has also been used to describe propagation. M\'el\'eard and Mirrahimi \cite{MeleardMirrahimi2015} studied singular limits for fractional equations, and Souganidis and Tarfulea \cite{SouganidisTarfulea2019} treated periodic KPP equations with singular integral diffusion. Tao, Li, Sun and Xu \cite{TaoEtAl2025} obtained interface limits for nonlocal dispersal equations with a small scaling parameter.

For bounded convolution diffusion, Bouin, Garnier, Henderson and Patout \cite[Proposition~2.1]{BouinEtAl2018} constructed comparison functions from the logistic profile $(1+e^{-t}/J(x))^{-1}$. The relative convolution estimate in \cite[Proposition~2.1]{BouinEtAl2018} produces factors involving $\exp(\pm\int_0^t\vartheta(s)\dd s)$, with $\vartheta(t)\to0$. Remark~2.2 of \cite{BouinEtAl2018} identifies $\int_0^\infty\vartheta(s)\dd s<\infty$ as a sufficient condition for constant-factor bounds. We follow this comparison construction and prove explicit integrable error estimates for the reference tail $(1+|x|)^{-\gamma}$.

Li \cite[Theorems~1.1--1.2 and Remark~1.1]{LiCooperative2026} studied cooperative systems with strictly subhomogeneous reactions and compactly supported initial data. The kernels are symmetric and nonincreasing and satisfy the long-tail condition (1.11) and doubling bound (1.13) in \cite{LiCooperative2026}. In that setting, let $\rho^*>0$ denote the principal eigenvalue of the reaction Jacobian at zero. When the kernels are comparable to $|x|^{-\gamma}$ with a common $\gamma>1$, the cited theorems give $\log|x|=(\rho^*/\gamma)t+o(t)$ for points in any component level set with a fixed level strictly between zero and the corresponding positive equilibrium.

Theorem~\ref{thm:main}\textup{(ii)} gives constant-factor localization under \eqref{eq:J}--\eqref{eq:F} for every $\gamma>1$. The proof combines \eqref{eq:global-main} with Lemma~\ref{lem:amplification}, which raises a positive lower bound to any fixed level below one. Estimate \eqref{eq:global-main} also yields the arrival-time bounds in \eqref{eq:arrival} and the integral bounds in \eqref{eq:moment-growth}. Section~\ref{sec:ranges} relates $m_1(J)$ and $m_2(J)$ to \eqref{eq:truncated-error} and \eqref{eq:smooth-error}; Section~\ref{sec:tail} proves the formula for $\mathcal A(t)$ in \eqref{eq:tail-formula}.

For degenerate monostable reactions, Alfaro and Coville \cite{AlfaroCoville2017} studied the interaction between the reaction near zero and the kernel tail, including finite-speed regimes. Li, Pan and Su \cite{LiPanSu2025} studied sharp acceleration for recursive systems, and Li \cite{LiAge2026} treated age-structured models. Cui, Li and Xu \cite{CuiLiXu2025} examined the effect of reducible coupling in cooperative recursions.

For free-boundary models, Cao, Du, Li and Li \cite{CaoDuLiLi2019} established a spreading--vanishing dichotomy, and Du, Li and Zhou \cite{DuLiZhou2021} characterized spreading speed through a semi-wave and the condition $m_1(J)<\infty$. Du and Ni \cite{DuNiJEMS} obtained rates for polynomial kernels and identified acceleration constants under an actual tail limit in \cite{DuNiMathAnn}. Section~\ref{sec:examples} compares the flux law \eqref{eq:fb-flux-comparison} with the Cauchy level-set estimate \eqref{eq:log-main}.

For local Fisher--KPP diffusion, Bramson \cite{Bramson1978,Bramson1983} established a logarithmic delay. In the nonlocal thin-tail setting, let $\widetilde J$ be a smooth even compactly supported probability density and set
\begin{equation}\label{eq:thin-tail-speed}
 c(a)=\frac{d\left[\int_{\R}e^{az}\widetilde J(z)\dd z-1\right]+r}{a},
 \qquad c_*=\min_{a>0}c(a)=c(\lambda_*).
\end{equation}
Graham \cite{Graham2022} proved the nonlinear delay $-(3/(2\lambda_*))\log t$ for compactly supported jump laws and step-like data. For compactly supported data, Roquejoffre \cite[Theorem~4.1.1]{Roquejoffre2024} gives a rightmost level of the form $c_*t-(3/(2\lambda_*))\log t+C_\theta+o(1)$ under the hypotheses of Chapter~4. Boutillon \cite{Boutillon2024} obtained the coefficient $1/(2\lambda_*)$ for the linearized step-data problem. Roquejoffre and Zhang \cite[Theorem~1.1]{RoquejoffreZhang2026} proved a multidimensional expansion with directional speed, logarithmic delay and a direction-dependent constant for smooth compactly supported kernels. Jin and Zhou \cite{JinZhou2026} studied time-dependent asymmetric kernels under exponential-boundedness and uniform-mean assumptions.

\subsection{The main estimate}

Set
\begin{equation}\label{eq:scale-envelope}
 S_\gamma(t)=\ee^{rt/\gamma},\qquad
 H_\gamma(t,x)=\min\left\{1,\frac{\ee^{rt}}{(1+|x|)^\gamma}\right\}.
\end{equation}
The next theorem answers the opening question by locating every fixed level between two time-independent multiples of $S_\gamma(t)$, throughout the full range of integrable polynomial tails.

\begin{theorem}\label{thm:main}
Assume \eqref{eq:J}, \eqref{eq:tail}, \eqref{eq:data}, and \eqref{eq:F}. Then the following hold.
\begin{enumerate}
\item There are constants $c,C>0$ such that
\begin{equation}\label{eq:global-main}
 cH_\gamma(t,x)\le u(t,x)\le CH_\gamma(t,x)
 \text{ for } t\ge1,\ x\in\R.
\end{equation}
\item For every $\theta\in(0,1)$ there exist $c_\theta,C_\theta,T_\theta>0$ such that, for $t\ge T_\theta$,
\begin{equation}\label{eq:main-inclusions}
 [-c_\theta S_\gamma(t),c_\theta S_\gamma(t)]
 \subset\{x:u(t,x)>\theta\}
 \subset\Omega_\theta(t)
 \subset[-C_\theta S_\gamma(t),C_\theta S_\gamma(t)].
\end{equation}
For every $t\ge T_\theta$, the set $E_\theta(t)$ is nonempty and compact and intersects both $(-\infty,0)$ and $(0,\infty)$. There is also a constant $K_\theta>0$ such that
\begin{equation}\label{eq:log-main}
 \sup_{x\in E_\theta(t)}
 \left|\log|x|-\frac r\gamma t\right|\le K_\theta
 \text{ for } t\ge T_\theta.
\end{equation}
In particular, $|x_\theta^\pm(t)|\sim S_\gamma(t)$ and $|x_\theta^\pm(t)|/t\to\infty$ as $t\to\infty$.
The constants and starting times in \eqref{eq:main-inclusions} and \eqref{eq:log-main} can be chosen uniformly when $\theta$ ranges over a compact subinterval of $(0,1)$.

\item If positive functions $a(t)$ and $b(t)$ satisfy $a(t)/S_\gamma(t)\to0$ and $b(t)/S_\gamma(t)\to\infty$ as $t\to\infty$, then
\begin{equation}\label{eq:inner-outer}
 \lim_{t\to\infty}\sup_{|x|\le a(t)}|u(t,x)-1|=0,
 \qquad
 \lim_{t\to\infty}\sup_{|x|\ge b(t)}u(t,x)=0.
\end{equation}
\end{enumerate}
\end{theorem}

We also record an integral consequence of \eqref{eq:global-main}. For $p\ge0$ and $q>0$, define the possibly infinite quantity
\begin{equation}\label{eq:moment-def}
 M_{p,q}(t)=\int_{\R}|x|^p u(t,x)^q\dd x.
\end{equation}
When $p=0$, the factor $|x|^p$ is understood to equal one everywhere. We write $\mathcal M(t)=M_{0,1}(t)$ for the mass of $u(t,\cdot)$ and $M_0=\mathcal M(0)=\int_{\R}u_0(x)\dd x$ for its initial mass.

\begin{corollary}\label{cor:moments}
Under the assumptions of Theorem~\ref{thm:main}, for each fixed $t>0$,
\begin{equation}\label{eq:moment-criterion}
 M_{p,q}(t)<\infty\quad\Longleftrightarrow\quad p+1<\gamma q.
\end{equation}
If $p+1<\gamma q$, there are positive constants $c_{p,q},C_{p,q}$ such that
\begin{equation}\label{eq:moment-growth}
 c_{p,q}\ee^{(p+1)rt/\gamma}
 \le M_{p,q}(t)\le C_{p,q}\ee^{(p+1)rt/\gamma}
 \text{ for } t\ge1.
\end{equation}
In particular, $\mathcal M(t)\sim S_\gamma(t)$ for $t\ge1$. For each $t>0$, the first absolute moment $M_{1,1}(t)=\int_{\R}|x|u(t,x)\dd x$ is finite exactly when $\gamma>2$, and the second moment $M_{2,1}(t)=\int_{\R}x^2u(t,x)\dd x$ is finite exactly when $\gamma>3$.
\end{corollary}

By Corollary~\ref{cor:moments}, the integrals $M_{p,q}(t)$ with $p+1\ge\gamma q$ are infinite at every positive time, although $M_{p,q}(0)<\infty$ for every $p\ge0$ and $q>0$ under \eqref{eq:data}.

\subsection{Further estimates}

Define
\begin{equation}\label{eq:W-def}
 W_R(x)=\left(1+\frac{x^2}{R^2}\right)^{-\gamma/2},\qquad R\ge2,
\end{equation}
and
\begin{equation}\label{eq:omega-def}
 \omega_\gamma(R)=
 \begin{cases}
 R^{1-\gamma},&1<\gamma<3,\\
 R^{-2}\log R,&\gamma=3,\\
 R^{-2},&\gamma>3.
 \end{cases}
\end{equation}

Proposition~\ref{thm:smooth-error} estimates $\sup_{x\in\R}|(J*W_R)(x)-W_R(x)|/W_R(x)$ by $\omega_\gamma(R)$. The order of $\omega_\gamma(R)$ changes at $\gamma=3$, the threshold for $m_2(J)<\infty$.

\begin{proposition}\label{thm:smooth-error}
Assume \eqref{eq:J} and \eqref{eq:tail}. There are constants $c,C,R_*>0$ such that
\begin{equation}\label{eq:smooth-error}
 c\omega_\gamma(R)
 \le \sup_{x\in\R}\frac{|(J*W_R)(x)-W_R(x)|}{W_R(x)}
 \le C\omega_\gamma(R)\text{ for } R\ge R_*.
\end{equation}
The lower estimate in \eqref{eq:smooth-error} also holds with the supremum replaced by $|(J*W_R)(0)-W_R(0)|/W_R(0)$. Consequently, the order in \eqref{eq:omega-def} cannot be improved for this family under the stated kernel hypotheses.
\end{proposition}

Section~\ref{sec:ranges} relates \eqref{eq:smooth-error} to the Fourier transform of $J$.

A tail coefficient can be identified if comparison in \eqref{eq:tail} is strengthened to an actual limit:
\begin{equation}\label{eq:exact-tail}
 \lim_{|x|\to\infty}|x|^\gamma J(x)=\lambda\in(0,\infty).
\end{equation}
Condition \eqref{eq:exact-tail} is imposed only in results that explicitly invoke it.

Set
\begin{equation}\label{eq:A-def}
	\mathcal A(t)=d\int_0^t\ee^{-rs}\mathcal M(s)\dd s,
	\qquad \beta_\gamma=r\left(1-\frac1\gamma\right)>0.
\end{equation}
\begin{theorem}\label{thm:tail-coefficient}
Assume the hypotheses of Theorem~\ref{thm:main} and \eqref{eq:exact-tail}. 
For every fixed $t>0$,
\begin{equation}\label{eq:tail-formula}
 \lim_{x\to+\infty}\frac{u(t,x)}{\ee^{rt}J(x)}
 =\lim_{x\to-\infty}\frac{u(t,x)}{\ee^{rt}J(x)}
 =\mathcal A(t).
\end{equation}
More precisely, for every $T>0$,
\begin{equation}\label{eq:tail-uniform-main}
 \lim_{X\to\infty}\sup_{\substack{0\le t\le T\\ |x|\ge X}}
 \left|\frac{u(t,x)}{\ee^{rt}J(x)}-\mathcal A(t)\right|=0.
\end{equation} Moreover,
\begin{equation}\label{eq:A-infty}
 \mathcal A_\infty=d\int_0^\infty\ee^{-rs}\mathcal M(s)\dd s\in(0,\infty),
 \qquad
 \mathcal A_\infty-\mathcal A(t)\sim\ee^{-\beta_\gamma t}
 \text{ for } t\ge1.
\end{equation}
These assertions hold for every $\gamma>1$.
\end{theorem}

Sections~\ref{sec:cauchy} and~\ref{sec:convolution} establish the comparison principles, positive-time tail bounds and relative convolution estimates. Section~\ref{sec:barriers} uses them to prove Theorem~\ref{thm:main} and its consequences. Section~\ref{sec:ranges} treats smooth comparison functions, and Section~\ref{sec:tail} proves Theorem~\ref{thm:tail-coefficient} and gives examples. Constants may depend on $d,f,u_0,J,\gamma$ and $j_1,j_2,R_J$ in \eqref{eq:tail-explicit}; other dependencies are specified where needed.

\section{The Cauchy problem}\label{sec:cauchy}

Assumptions \eqref{eq:J}, \eqref{eq:tail}, \eqref{eq:data}, and \eqref{eq:F} remain in force throughout this section, except where explicitly weakened.

To allow interval indicators as initial data in comparison arguments, let $\mathcal B$ be the Banach space of bounded Borel functions on $\R$, with norm $\|v\|_\infty=\sup_{x\in\R}|v(x)|$. Let $C_0(\R)$ be its closed subspace of continuous functions that tend to zero as $|x|\to\infty$. We write $\|v\|_1=\int_{\R}|v(x)|\dd x$ for $v\in L^1(\R)$, and use Lebesgue measure in all spatial integrals unless another measure is displayed. The symbol $\ind_A$ denotes the indicator of a Borel set $A$, and $\supp v$ denotes the closure of $\{v\ne0\}$ for a continuous function $v$.
For a function $F:\R\to\R$, write
\[
 \Lip(F)=\sup_{\substack{a,b\in\R\\a\ne b}}
 \frac{|F(a)-F(b)|}{|a-b|}.
\]
Thus $F$ is $L$-Lipschitz on $\R$ if $|F(a)-F(b)|\le L|a-b|$ for all $a,b\in\R$. Define
\begin{equation}\label{eq:operator}
 \Lop v=d(J*v-v).
\end{equation}
For a finite Borel measure $\nu$, set $(\nu*v)(x)=\int_{\R}v(x-z)\dd\nu(z)$. When applied to a time-dependent function, $\Lop$ acts on its spatial variable with time fixed. The measure $\delta_0$ has unit mass at zero, $J^{*1}=J$, and $J^{*(n+1)}=J*J^{*n}$ for $n\ge1$. We set $\N=\{1,2,\ldots\}$.

\subsection{The semigroup and comparison}

\begin{lemma}\label{lem:semigroup}
Under \eqref{eq:J}, the operator $\Lop$ is bounded on $\mathcal B$, on $C_0(\R)$, and on $L^1(\R)$. Its semigroup is
\begin{equation}\label{eq:P}
 P_t v=\ee^{-dt}\sum_{n=0}^{\infty}\frac{(dt)^n}{n!}J^{*n}*v,
 \qquad J^{*0}=\delta_0.
\end{equation}
The series converges in operator norm on each bounded time interval. The operators $P_t$ are positive, satisfy $P_t1=1$, and are contractions on $\mathcal B$ and $L^1(\R)$. They preserve the integral of an $L^1$ function.
\end{lemma}

\begin{proof}
Formula \eqref{eq:P} is the Poissonized transition formula in \cite[Preface, pp.~xii--xiii]{Andreu2010}, with time replaced by $dt$. We verify its realization on the three spaces in the statement.

Set $Tv=J*v$. For $v\in\mathcal B$, parameter integration gives Borel measurability and $\|Tv\|_\infty\le\|v\|_\infty$. For $v\in L^1(\R)$, Tonelli's theorem gives
\[
 \int_{\R}\int_{\R}J(z)|v(x-z)|\dd z\dd x=\|v\|_1.
\]
Hence $\|Tv\|_1\le\|v\|_1$, and Fubini's theorem yields $\int_{\R}Tv=\int_{\R}v$. If $v\in C_0(\R)$, dominated convergence, with majorant $\|v\|_\infty J(z)$, proves both continuity of $Tv(x)$ and $Tv(x)\to0$ as $|x|\to\infty$. Thus $T$ is a positive contraction on each space and $\|\Lop\|\le2d$.

By Tonelli's theorem, each $J^{*n}$ is nonnegative and has integral one. Fubini's theorem then identifies $T^nv=J^{*n}*v$: the absolute multiple integral is bounded by $\|v\|_\infty$ for bounded $v$, and its integral also over $x$ is bounded by $\|v\|_1$ for $v\in L^1(\R)$. Consequently, on every bounded time interval $[0,T_0]$, the exponential series and its derivative converge in operator norm, since
\[
 \sum_{n=0}^\infty\frac{(dT_0)^n}{n!}\|T^n\|\le\ee^{dT_0},
 \qquad
 \sum_{n=1}^\infty\frac{d(dT_0)^{n-1}}{(n-1)!}\|T^n\|\le d\ee^{dT_0}.
\]
The bounded-operator exponential therefore gives $P_t=\ee^{td(T-I)}=\ee^{-dt}\ee^{dtT}$ and \eqref{eq:P}. The exponential identity gives $P_tP_s=P_{t+s}$; termwise differentiation gives $\partial_tP_t=\Lop P_t$.

Finally, \eqref{eq:P} is a convex combination of the positive contractions $T^n$, and $T^n1=1$. It follows that $P_t$ is positive, contractive and satisfies $P_t1=1$. For $v\in L^1(\R)$, integration is a bounded linear functional and the series converges in $L^1(\R)$, so
\[
 \int_{\R}P_tv\dd x
 =\ee^{-dt}\sum_{n=0}^\infty\frac{(dt)^n}{n!}\int_{\R}T^nv\dd x
 =\int_{\R}v\dd x.
\]
\end{proof}

To express $P_t$ as an average over spatial displacements, write it as convolution by
\begin{equation}\label{eq:mu}
 \mu_t=\ee^{-dt}\delta_0+
 \ee^{-dt}\sum_{n=1}^{\infty}\frac{(dt)^n}{n!}J^{*n}(z)\dd z.
\end{equation}
The series in \eqref{eq:mu} consists of nonnegative finite measures. By Tonelli's theorem, its sum is countably additive, and Lemma~\ref{lem:semigroup} gives
\[
 \mu_t(\R)=\ee^{-dt}\sum_{n=0}^{\infty}\frac{(dt)^n}{n!}=1.
\]
The total variation norm of a tail of this series equals the corresponding tail of the scalar exponential series and therefore tends to zero.

The measure $\mu_t$ in \eqref{eq:mu} is a probability measure with an atom of mass $\ee^{-dt}$ at zero. Moreover, $\mu_t$ is the time-$t$ law of a compound Poisson process; see Applebaum \cite[Chapter~1]{Applebaum2009} and the Poissonized transition kernel in \cite[Preface, pp.~xii--xiii]{Andreu2010}. On a probability space with probability $\Prob$ and expectation $\E$, let $N_t$ be a Poisson counting process of rate $d$, starting at zero: its independent increments satisfy
\begin{equation}\label{eq:Poisson-count}
 \Prob(N_t-N_s=n)=\ee^{-d(t-s)}\frac{[d(t-s)]^n}{n!}
 \text{ for } 0\le s<t,\ n=0,1,\ldots.
\end{equation}
Let $(Z_k)_{k\ge1}$ be independent random variables with common density $J$, independent of the Poisson process $(N_t)_{t\ge0}$, and put $X_t=\sum_{k=1}^{N_t}Z_k$, with the empty sum equal to zero. Conditioning on $N_t=n$ gives, for every bounded Borel function $v$,
\begin{equation}\label{eq:Poisson-representation}
 \Prob(X_t\in A)=\mu_t(A),\qquad
 \E[v(x-X_t)]=(P_tv)(x),
 \quad A\subset\R\ \text{Borel}.
\end{equation}
For $n\ge1$, the sum $Z_1+\cdots+Z_n$ has density $J^{*n}$ by independence and the convolution identity verified in Lemma~\ref{lem:semigroup}; for $n=0$, its law is $\delta_0$. Thus conditioning on $N_t$ gives \eqref{eq:mu} and \eqref{eq:Poisson-representation}. For bounded Borel $v$, each conditional expectation is bounded in absolute value by $\|v\|_\infty$, and the weights $\Prob(N_t=n)$ sum to one; hence the weighted series converges absolutely. The same conditioning is used for the linearized step-data equation in \cite[Proposition~2.1]{Boutillon2024}.

We use the following convention for lower and upper solutions. Outside $[0,1]$, extend $f$ by the tangent lines $rs$ for $s<0$ and $f'(1)(s-1)$ for $s>1$. A lower solution is a bounded function $\underline v(t,\cdot)\in\mathcal B$ with the time regularity stated in the next lemma and satisfying $\partial_t\underline v\le\Lop\underline v+f(\underline v)$ for almost every $t$ at each $x$; an upper solution satisfies the reverse inequality. The initial values must be ordered for comparison. All comparison functions used below take values in $[0,1]$. For a Banach space $B$, the notation $C^1([0,\infty);B)$ means that difference quotients converge in $\|\cdot\|_B$ to a continuous $B$-valued derivative, with a right derivative at zero. Banach-space-valued time integrals below are Bochner integrals; a norm-continuous integrand on a compact time interval is Bochner integrable.

\begin{lemma}\label{lem:wellposed}
Suppose \eqref{eq:J} and \eqref{eq:F} hold. For every $v_0\in\mathcal B$ with $0\le v_0\le1$, equation \eqref{eq:cauchy} with datum $v_0$ has a unique solution in $C^1([0,\infty);\mathcal B)$. The solution lies in $[0,1]$. If $v_0\in C_0(\R)$, it belongs to $C^1([0,\infty);C_0(\R))$.

Comparison holds for bounded lower and upper solutions whose values stay in a common bounded interval, provided they are continuous in time in the supremum norm, locally absolutely continuous in time at each spatial point, and, for every $T>0$, there is $g_T\in L^1(0,T)$ such that
\[
 |\partial_t\underline v(t,x)|+|\partial_t\overline v(t,x)|\le g_T(t)
\]
for almost every $t\in(0,T)$ at each $x\in\R$. If $\underline v(0,x)\le\overline v(0,x)$ for all $x\in\R$, then $\underline v(t,x)\le\overline v(t,x)$ for all $t\ge0$ and $x\in\R$.
\end{lemma}

\begin{proof}
The Banach-space Cauchy--Lipschitz (Picard--Lindel\"of) theorem gives a unique global norm-$C^1$ solution for a globally Lipschitz vector field; see \cite[Section~2.2 and Remark~2.2.1]{Roquejoffre2024}. Here the tangent-line extension of $f$ is $C^1(\R)$ and has Lipschitz constant $L=\max_{[0,1]}|f'|$. Since $f(0)=0$, composition by $f$ preserves both $\mathcal B$ and $C_0(\R)$. By Lemma~\ref{lem:semigroup}, on either space $B$ we have
\[
 \|\Lop v+f(v)-\Lop w-f(w)\|_B
 \le(2d+L)\|v-w\|_B.
\]
The Cauchy--Lipschitz theorem therefore gives a unique local solution in $C^1([0,T];B)$. Its integral equation and Gronwall's inequality yield $\|v(t)\|_B\le\ee^{(2d+L)t}\|v_0\|_B$, so the solution extends globally. Uniqueness in $\mathcal B$ identifies the two constructions when $v_0\in C_0(\R)$.

For comparison, we adapt the positive-part argument of \cite[Proposition~2.2.3]{Roquejoffre2024} to the time regularity in the statement; comparison in $C^1([0,T];L^\infty(\R))$ also follows from \cite[Lemma~15]{Yagisita2009}. Fix $T>0$ and put
\[
 w(t,x)=\underline v(t,x)-\overline v(t,x),\qquad
 q(t)=\sup_{x\in\R}\max\{w(t,x),0\},\qquad H(s)=f(s)+Ls.
\]
The function $q$ is finite and continuous by the assumed supremum-norm continuity. Moreover, $H$ is nondecreasing and $2L$-Lipschitz, so $H(a)-H(b)\le2L\max\{a-b,0\}$ for all $a,b\in\R$. Since $J\ge0$ and $\int_{\R}J=1$, subtraction of the differential inequalities gives
\[
 w_t+(d+L)w
 \le dJ*w+H(\underline v)-H(\overline v)
 \le(d+2L)q(t)
\]
for almost every time at each fixed $x$. Absolute continuity and the integrable derivative bound in the statement justify multiplication by $\ee^{(d+L)t}$ and time integration. Using $w(0,x)\le0$, we obtain
\[
 w(t,x)\le(d+2L)\int_0^t\ee^{-(d+L)(t-s)}q(s)\dd s.
\]
Taking positive parts and the spatial supremum yields $q(t)\le(d+2L)\int_0^tq(s)\dd s$. Gronwall's inequality gives $q\equiv0$ on $[0,T]$, proving comparison. In the estimate $q(t)\le(d+2L)\int_0^tq(s)\dd s$, the properties of $J$ used are its nonnegativity and unit mass, both supplied by \eqref{eq:J}.

Finally, $0$ and $1$ are constant solutions. The norm-$C^1$ solution constructed above satisfies the stated time regularity on every bounded interval; comparison with these constants gives $0\le v\le1$. Hence it solves the equation with the original reaction on $[0,1]$.
\end{proof}

Lemma~\ref{lem:wellposed} applies to $\overline u$ in \eqref{eq:upper-def}. At each fixed $x$, the formula for $\partial_t\overline u(t,x)$ may change once, at $1+|x|=R_+(t)$, while $t\mapsto\overline u(t,\cdot)$ is Lipschitz in $\|\cdot\|_\infty$. The differential inequality is verified on both sides of this time in Lemma~\ref{lem:upper}.

\begin{lemma}\label{lem:stability}
Let $v,w\in C^1([0,\infty);\mathcal B)$ solve the evolution equation in \eqref{eq:cauchy} and take values in $[0,1]$. Set $L=\max_{[0,1]}|f'|$. Then
\begin{equation}\label{eq:stability}
 |v(t,x)-w(t,x)|
 \le \ee^{Lt}P_t(|v(0,\cdot)-w(0,\cdot)|)(x)
 \text{ for } t\ge0,\ x\in\R.
\end{equation}
The solution $u$ of \eqref{eq:cauchy} also satisfies
\begin{equation}\label{eq:linear-sandwich}
 P_tu_0\le u(t,\cdot)\le\ee^{rt}P_tu_0.
\end{equation}
\end{lemma}

\begin{proof}

We first derive the variation-of-constants formula. If $z\in C^1([0,T];\mathcal B)$ satisfies $z_t=\Lop z+g(t)$ with continuous $g$, fix $t\le T$ and differentiate $P_{t-s}z(s)$ in the norm of $\mathcal B$. Lemma~\ref{lem:semigroup} and the product rule give
\[
\frac{\dd }{\dd s}[P_{t-s}z(s)]
=-\Lop P_{t-s}z(s)+P_{t-s}(\Lop z(s)+g(s))=P_{t-s}g(s).
\]
The cancellation uses commutation of $\Lop $ with its operator exponential. The map $s\mapsto P_{t-s}g(s)$ is norm continuous on $[0,t]$, hence Bochner integrable. The fundamental theorem of calculus for norm-$C^1$ paths therefore gives
\[
z(t)=P_tz(0)+\int_0^tP_{t-s}g(s)\dd s.
\]

Apply the variation-of-constants formula just obtained to $z=v-w$, with $g=f(v)-f(w)$. Because $v,w\in[0,1]$, the mean value theorem gives $|g(s,x)|\le L|v(s,x)-w(s,x)|$. Positivity of $P_t$ implies $|P_t\varphi|\le P_t|\varphi|$; this also follows directly from integration against the probability measure $\mu_t$. Thus, putting $a(t,x)=|v(t,x)-w(t,x)|$ and $a_0(x)=a(0,x)$,
\[
a(t)\le P_ta_0+L\int_0^tP_{t-s}a(s)\dd s.
\]
Evaluation at $x$ is a bounded linear functional on $\mathcal B$, so the norm integral and its pointwise evaluation agree.

We next iterate the integral inequality for $a$. Define $\mathcal K$ on $C([0,T];\mathcal B)$ by
\[
(\mathcal Ka)(t)=L\int_0^tP_{t-s}a(s)\dd s,
\]
and put $b(t)=P_ta_0$. Positivity of $\mathcal K$ gives, for every integer $n\ge0$,
\[
a(t)\le\sum_{k=0}^n(\mathcal K^kb)(t)+(\mathcal K^{n+1}a)(t).
\]
The semigroup identity and integration over an ordered time simplex yield
\[
(\mathcal K^kb)(t)=\frac{L^kt^k}{k!}P_ta_0.
\]
For example, the induction step uses
$L\int_0^tP_{t-s}[L^ks^kP_sa_0/k!]\dd s
=L^{k+1}P_ta_0\int_0^ts^k\dd s/k!$.
All integrands are norm continuous and bounded on the compact simplex. The norms of their multiple integrals are bounded by $L^{k+1}\|a_0\|_\infty$ times the simplex volume $t^{k+1}/(k+1)!$; hence the required changes of integration order are also justified by absolute integrability. The volume formula follows by induction, integrating $s^k/k!$ from $0$ to $t$.

Let $M_t=\sup_{0\le s\le t}\|a(s)\|_\infty$, which is finite (indeed at most one). Contraction of $P_t$ gives
\[
\|\mathcal K^{n+1}a(t)\|_\infty
\le M_t\frac{(Lt)^{n+1}}{(n+1)!}\longrightarrow0.
\]
For fixed $t$, the factorial dominates the power, as is seen by the ratio $Lt/(n+2)\to0$. Letting $n\to\infty$ proves
\[
|v(t,x)-w(t,x)|\le \ee^{Lt}(P_ta_0)(x),
\]
which is \eqref{eq:stability}, including $t=0$.

It remains to prove \eqref{eq:linear-sandwich}. Since $u\in[0,1]$, \eqref{eq:F} gives $0\le f(u)\le ru$. Variation of constants with generator $\Lop $ first yields
\[
u(t)=P_tu_0+\int_0^tP_{t-s}f(u(s))\dd s\ge P_tu_0.
\]
For the other inequality use the bounded generator $\Lop +rI$, whose exponential is $\ee^{rt}P_t$. The equation can be written
$u_t=(\Lop +rI)u-[ru-f(u)]$, so
\[
u(t)=\ee^{rt}P_tu_0-\int_0^te^{r(t-s)}P_{t-s}[ru(s)-f(u(s))]\dd s
\le \ee^{rt}P_tu_0.
\]
The integrands are norm continuous and nonnegative before the subtraction; \eqref{eq:F} and positivity of $P_t$ justify the order inequalities. This proves \eqref{eq:linear-sandwich}.

\end{proof}

\subsection{Bounds at positive times}

Introduce the reference function
\begin{equation}\label{eq:h}
 h_\gamma(x)=(1+|x|)^{-\gamma}.
\end{equation}
By the boundedness of $J$ in \eqref{eq:J} and the tail bounds in \eqref{eq:tail-explicit}, there are constants $c_J,C_J>0$ such that
\begin{equation}\label{eq:J-h}
 J(x)\le C_Jh_\gamma(x)\text{ for } x\in\R,\qquad
 J(x)\ge c_Jh_\gamma(x)\text{ for } |x|\ge R_J.
\end{equation}

\begin{lemma}\label{lem:weight}
For $\gamma>1$, there is $C_h>0$ such that
\begin{equation}\label{eq:h-convolution}
 (h_\gamma*h_\gamma)(x)\le C_hh_\gamma(x)\text{ for } x\in\R.
\end{equation}
The space
\begin{equation}\label{eq:Y}
 \Yg=\left\{v\in C_0(\R):
 \|v\|_\gamma:=\sup_{x\in\R}\frac{|v(x)|}{h_\gamma(x)}<\infty\right\}
\end{equation}
is a Banach space. Convolution by $J$ is bounded on $\Yg$, and the solution of \eqref{eq:cauchy} lies in $C^1([0,\infty);\Yg)$.
\end{lemma}

\begin{proof}

We first prove \eqref{eq:h-convolution}. Since $\gamma>1$,
\[
\int_{\R}h_\gamma(y)\dd y=2\int_0^\infty(1+y)^{-\gamma}\dd y
=\frac{2}{\gamma-1}<\infty.
\]
For fixed $x$, split $\R$ into $A_x=\{y:|y|\ge|x|/2\}$ and its complement. On $A_x$,
$1+|y|\ge(1+|x|)/2$, hence $h_\gamma(y)\le2^\gamma h_\gamma(x)$.
On the complement, the triangle inequality gives $|x-y|>|x|/2$ and $h_\gamma(x-y)\le2^\gamma h_\gamma(x)$. Therefore
\[
\begin{aligned}
(h_\gamma*h_\gamma)(x)
&\le2^\gamma h_\gamma(x)\int_{A_x}h_\gamma(x-y)\dd y
+2^\gamma h_\gamma(x)\int_{\R\setminus A_x}h_\gamma(y)\dd y\\
&\le2^{\gamma+1}h_\gamma(x)\int_{\R}h_\gamma(y)\dd y
=\frac{2^{\gamma+2}}{\gamma-1}h_\gamma(x).
\end{aligned}
\]
The integrands are nonnegative, so enlarging the two integration sets is valid. This proves \eqref{eq:h-convolution} with $C_h=2^{\gamma+2}/(\gamma-1)$.

To prove that $\Yg$ is complete, let $C_b(\R)$ be the bounded continuous functions with the supremum norm. The map
\[
v\longmapsto v/h_\gamma
\]
is linear and isometric from $\Yg$ into $C_b(\R)$. It is onto: for $g\in C_b(\R)$, the function $v=h_\gamma g$ is continuous, tends to zero at infinity and has $\|v\|_\gamma=\|g\|_\infty$. If $g_n$ is Cauchy in $C_b$, its pointwise limit exists, the convergence is uniform by the Cauchy criterion, and the limit is bounded and continuous. Thus $C_b$ is complete, and the isometry proves completeness of $\Yg$ in $\|\cdot\|_\gamma$.

We next estimate convolution by $J$ in $\|\cdot\|_\gamma$. For $v\in\Yg$, $|v(x-z)|\le\|v\|_\gamma h_\gamma(x-z)$. By \eqref{eq:J-h} and \eqref{eq:h-convolution},
\[
\begin{aligned}
|J*v(x)|
&\le C_J\|v\|_\gamma\int_{\R}h_\gamma(z)h_\gamma(x-z)\dd z\\
&\le C_JC_h\|v\|_\gamma h_\gamma(x).
\end{aligned}
\]
Lemma~\ref{lem:semigroup} also gives $J*v\in C_0(\R)$. Hence $J*v\in\Yg$, and
\begin{equation}\label{eq:J-weight-bound}
 \|J*v\|_\gamma\le C_JC_h\|v\|_\gamma.
\end{equation}

For the tangent-line extension of $f$, put $L=\max_{[0,1]}|f^{\prime}|$, which is a global Lipschitz constant by Lemma~\ref{lem:wellposed}. Then $f(v)$ is continuous, tends to zero and satisfies
\[
\|f(v)\|_\gamma\le L\|v\|_\gamma,\qquad
\|f(v)-f(w)\|_\gamma\le L\|v-w\|_\gamma.
\]
Both inequalities follow by dividing the pointwise scalar inequalities by the positive $h_\gamma(x)$ before taking the supremum.

Thus $v\mapsto\Lop v+f(v)$ is globally Lipschitz on $\Yg$ with constant at most $d(C_JC_h+1)+L$.

Choose $R_u$ with $\operatorname{supp}u_0\subset[-R_u,R_u]$. Then
$\|u_0\|_\gamma\le(1+R_u)^\gamma\|u_0\|_\infty<\infty$.
The Cauchy--Lipschitz argument used in Lemma~\ref{lem:wellposed}, now with the norm $\|\cdot\|_\gamma$ and time step satisfying $\Delta[d(C_JC_h+1)+L]<1$, gives a global $C^1$ solution in $\Yg$.
Since $h_\gamma\le1$, $\|v\|_\infty\le\|v\|_\gamma$. The continuous embedding $\Yg\hookrightarrow C_0(\R)$ shows that the $\Yg$-valued solution also satisfies \eqref{eq:cauchy} as an equation in $\mathcal B$. Uniqueness in Lemma~\ref{lem:wellposed} identifies it with $u$.

The integral equation and Gronwall further give the explicit finite-time bound
\[
\|u(t)\|_\gamma\le\|u_0\|_\gamma
\exp\{[d(C_JC_h+1)+L]t\}.
\]
Also $\|v\|_1\le[2/(\gamma-1)]\|v\|_\gamma$ for every $v\in\Yg$, by integration of $|v|\le\|v\|_\gamma h_\gamma$. Thus $\sup_{0\le t\le T}\|u(t,\cdot)\|_\gamma<\infty$ for every $T>0$, and integration is controlled by the displayed $L^1$ bound.

\end{proof}

By Lemma~\ref{lem:weight}, $\sup_{0\le t\le T}\|u(t,\cdot)\|_\gamma<\infty$ for every $T>0$. Since $\|v\|_1\le\|h_\gamma\|_1\|v\|_\gamma$, the path $t\mapsto u(t,\cdot)$ also belongs to $C^1([0,\infty);L^1(\R))$, which justifies differentiating $\mathcal M(t)=\int_{\R}u(t,x)\dd x$ below.

\begin{lemma}\label{lem:seed}
For every $t>0$ and $x\in\R$, one has $0<u(t,x)<1$. There is $b_0>0$ such that
\begin{equation}\label{eq:seed}
 u(1,x)\ge b_0h_\gamma(x)\text{ for } x\in\R.
\end{equation}
More generally, for $0<\tau<T<\infty$ there are $b_{\tau,T},B_{\tau,T}>0$ such that
\begin{equation}\label{eq:finite-tail}
 b_{\tau,T}h_\gamma(x)\le u(t,x)\le B_{\tau,T}h_\gamma(x)
 \text{ for } \tau\le t\le T,\ x\in\R.
\end{equation}
\end{lemma}

\begin{proof}

Fix $R_u>0$ such that $\operatorname{supp}u_0\subset[-R_u,R_u]$ for this proof.

We first prove $u(t,x)>0$ for $t>0$ and $x\in\R$. Continuity of $J$ and $J(0)>0$ supply $a_0,\delta>0$ such that $J(z)\ge a_0$ for $|z|\le\delta$. For $n\ge2$ and $|s|<n\delta$, set $b=s/n$ and choose
\[
\eta=\frac{\delta-|b|}{2n}>0.
\]
In the integral formula
\[
J^{*n}(s)=\int_{\R^{n-1}}J(z_1)\cdots J(z_{n-1})
J\left(s-\sum_{i=1}^{n-1}z_i\right)\dd z_1\cdots\dd z_{n-1},
\]
restrict each $z_i$ to $(b-\eta,b+\eta)$. Then $|z_i|<\delta$, and
\[
\left|s-\sum_{i=1}^{n-1}z_i\right|
\le |b|+(n-1)\eta<\delta.
\]
The product is at least $a_0^n$ on a box of volume $(2\eta)^{n-1}>0$. Therefore $J^{*n}(s)>0$ for $|s|<n\delta$. The case $n=1$ follows from the initial lower bound. Tonelli applies to this nonnegative iterated integral; boundedness of $J$ and unit masses of the other factors also bound the integral by $\|J\|_\infty$.

By continuity and nontriviality of $u_0$, some interval $I=[y_0-\eta_0,y_0+\eta_0]$ and $b>0$ satisfy $u_0(y)\ge b$ on $I$. For a fixed $x$, choose an integer $n$ so that $|x-y_0|+\eta_0<n\delta$. Then $J^{*n}(x-y)>0$ for every $y\in I$, and
\[
(J^{*n}*u_0)(x)\ge b\int_IJ^{*n}(x-y)\dd y>0.
\]
The last integral is positive because its measurable integrand is positive throughout an interval of positive length. For $t>0$, the $n$-th coefficient of \eqref{eq:P} is strictly positive. Thus $P_tu_0(x)>0$, and \eqref{eq:linear-sandwich} proves $u(t,x)>0$.

To prove $u(t,x)<1$, define $w(t,x)=1-u(t,x)$. Lemma~\ref{lem:wellposed} gives $w\ge0$, and $w$ is $C^1$ in $\mathcal B$. Since $f(1)=0$ and $L=\max_{[0,1]}|f'|$,
\[
f(1-w)\le Lw\text{ for }0\le w\le1,\qquad
w_t=\Lop w-f(1-w)\ge(\Lop -LI)w.
\]
Variation of constants with generator $\Lop -LI$ gives
\[
w(t)=\ee^{-Lt}P_t(1-u_0)
+\int_0^te^{-L(t-s)}P_{t-s}[Lw(s)-f(1-w(s))]\dd s
\ge \ee^{-Lt}P_t(1-u_0).
\]
The integrand is norm continuous and nonnegative, which justifies both the integral and the order inequality.

Since $u_0$ vanishes outside $[-R_u,R_u]$, the function $1-u_0$ equals one there. Fix $x$ and choose $A>\max\{R_J,|x|+R_u+1\}$. For $z\in[A,A+1]$, $x-z<-R_u$ and \eqref{eq:J-h} gives $J(z)\ge c_Jh_\gamma(z)>0$. Hence
\[
J*(1-u_0)(x)\ge c_J\int_A^{A+1}h_\gamma(z)\dd z>0.
\]
The one-jump term of \eqref{eq:P} then gives
$w(t,x)\ge \ee^{-Lt}\ee^{-dt}dt\,[J*(1-u_0)](x)>0$ for $t>0$. Thus $u(t,x)<1$.

We next establish the lower estimate in \eqref{eq:finite-tail}. Fix $0<\tau<T$. For $|x|\ge R_J+R_u$ and $|y|\le R_u$, $|x-y|\ge R_J$. Since
$1+|x-y|\le1+|x|+R_u\le(1+|x|)(1+R_u)$, \eqref{eq:J-h} implies
\[
J(x-y)\ge c_J(1+R_u)^{-\gamma}h_\gamma(x).
\]
Using \eqref{eq:linear-sandwich} and the one-jump term of \eqref{eq:P}, we obtain
\begin{equation}\label{eq:one-jump}
 u(t,x)\ge dt\ee^{-dt}\int_{-R_u}^{R_u}J(x-y)u_0(y)\dd y
 \ge c_Jdt\ee^{-dt}(1+R_u)^{-\gamma}
       h_\gamma(x)\int_{\R}u_0(y)\dd y.
\end{equation}
Here $M_0=\int_{\R}u_0(y)\dd y>0$ by \eqref{eq:data}.
For $\tau\le t\le T$, $dt \ee^{-dt}\ge d\tau \ee^{-dT}>0$. Thus \eqref{eq:one-jump} gives a uniform positive lower bound for $u(t,x)/h_\gamma(x)$ when $\tau\le t\le T$ and $|x|\ge R_J+R_u$.

On the remaining compact region, joint continuity follows from continuity of $t\mapsto u(t,\cdot)$ in the supremum norm and spatial continuity: as $(t,x)\to(s,y)$,
\[
|u(t,x)-u(s,y)|\le\|u(t)-u(s)\|_\infty+|u(s,x)-u(s,y)|\longrightarrow0.
\]
The positive function $u(t,x)/h_\gamma(x)$ is therefore continuous on
$[\tau,T]\times[-R_J-R_u,R_J+R_u]$. The extreme value theorem gives a positive minimum there, since a zero minimum would be attained and contradict the already proved strict positivity. For the lower constant in \eqref{eq:finite-tail}, take the minimum of
\[
 \min_{\substack{\tau\le t\le T\\ |x|\le R_J+R_u}}
       \frac{u(t,x)}{h_\gamma(x)}
 \quad\text{and}\quad c_Jd\tau\ee^{-dT}(1+R_u)^{-\gamma}M_0.
\]
The first quantity is positive by compactness, and the second controls the remaining $x$ by \eqref{eq:one-jump}.

By Lemma~\ref{lem:weight}, $u(t,x)\le\|u(t)\|_\gamma h_\gamma(x)$, and continuity of $t\mapsto u(t,\cdot)$ in $\Yg$ gives $\sup_{\tau\le t\le T}\|u(t)\|_\gamma<\infty$. This proves the upper estimate in \eqref{eq:finite-tail}. Taking, for example, $\tau=1$ and $T=2$ and evaluating at $t=1$ proves \eqref{eq:seed}.

\end{proof}

By \eqref{eq:seed}, the lower comparison can start at time one with a positive multiple of $h_\gamma$. For related positive-time lower bounds, see Garnier \cite[proof of Theorem~2]{Garnier2011} and \cite[Remark~1.4]{BouinEtAl2018}.

\subsection{A lower bound at any fixed level}

\begin{lemma}\label{lem:amplification}
For each $\delta,\theta\in(0,1)$, there exist $T,R_*>0$ such that, for every $R\ge R_*$, the solution $v_R$ of \eqref{eq:cauchy} with $u_0=\delta\ind_{[-R,R]}$ satisfies
\begin{equation}\label{eq:amplification}
 v_R(T,x)>\theta\text{ for } R\ge R_*,\ |x|\le R/2.
\end{equation}
For every $x_0\in\R$, the solution of \eqref{eq:cauchy} with $u_0=\delta\ind_{[x_0-R,x_0+R]}$ is $v_R(t,x-x_0)$. Hence $v_R(T,x-x_0)>\theta$ for $R\ge R_*$ and $|x-x_0|\le R/2$, with the same $T$ and $R_*$.
\end{lemma}

\begin{proof}

Set $L=\max_{[0,1]}|f^{\prime}|$. Let $y$ be the solution of $y'=f(y)$ with $y(0)=\delta$. The Lipschitz extension used in Lemma~\ref{lem:wellposed} supplies a unique global scalar solution, and comparison with the equilibria gives $0\le y\le1$. For $z=1-y$, the Lipschitz bound at one gives
$z'=-f(y)\ge-Lz$, and multiplication by $\ee^{Lt}$ yields
$1-y(t)\ge(1-\delta)\ee^{-Lt}>0$ at every finite time. Because $f(s)>0$ on $(0,1)$, $y$ is strictly increasing and $y(t)\ge\delta>0$.

The bounded monotone function $y$ has a limit $\ell\in[\delta,1]$. If $\ell<1$, continuity and $f(\ell)>0$ give an interval about $\ell$ on which $f\ge c>0$. For all sufficiently large times, $y$ lies in that interval, so integrating $y'\ge c$ forces $y$ above one, a contradiction. Therefore $y(t)\to1$. Choose $T>0$ with $y(T)>\theta$.

As $\int J=1$, a spatially constant function has $\Lop  y=0$. Thus $(t,x)\mapsto y(t)$ is the solution of \eqref{eq:cauchy} with constant initial datum $\delta$, in the class supplied by Lemma~\ref{lem:wellposed}. The initial difference from $v_R$ is exactly
\[
|v_R(0,x)-\delta|=\delta\ind_{\{|x|>R\}}.
\]
Apply \eqref{eq:stability} at time $T$ and use $P_T\varphi(x)=\int_{\R}\varphi(x-z)\dd \mu_T(z)$. For $|x|\le R/2$, the event $|x-z|>R$ implies $|z|\ge|x-z|-|x|>R/2$. Consequently,
\begin{align}\label{eq:amp-error}
 |v_R(T,x)-y(T)|
 &\le\delta\ee^{LT}\int_{\R}\ind_{\{|x-z|>R\}}\dd\mu_T(z)\notag\\
 &\le\delta\ee^{LT}\mu_T(\{|z|>R/2\}).
\end{align}

At the fixed time $T$, $\mu_T$ is a probability measure. The indicators $\ind_{\{|z|>R/2\}}$ tend pointwise to zero as $R\to\infty$, and are dominated by the integrable function $1$ on $(\R,\mu_T)$. The dominated convergence theorem gives
$\mu_T(\{|z|>R/2\})\to0$. Thus one can choose $R_*>0$ so that
\[
\delta \ee^{LT}\mu_T(\{|z|>R_*/2\})<\frac{y(T)-\theta}{2}.
\]
For $R\ge R_*$, inclusion $\{|z|>R/2\}\subset\{|z|>R_*/2\}$ gives $\delta\ee^{LT}\mu_T(\{|z|>R/2\})<(y(T)-\theta)/2$. Consequently
\[
v_R(T,x)\ge y(T)-\delta \ee^{LT}\mu_T(\{|z|>R/2\})
>\frac{y(T)+\theta}{2}>\theta
\text{ for }R\ge R_*,\ |x|\le R/2.
\]
This proves \eqref{eq:amplification} with constants independent of $R$ and $x$ in the stated ranges.

For any $x_0\in\R$, set $v_R^{x_0}(t,x)=v_R(t,x-x_0)$. The identity
\[
(J*v_R^{x_0}(t,\cdot))(x)
=\int_{\R}J(z)v_R(t,x-x_0-z)\dd z
=(J*v_R(t,\cdot))(x-x_0)
\]
shows that $v_R^{x_0}$ satisfies the evolution equation in \eqref{eq:cauchy}. Its initial datum is
$\delta\ind_{[x_0-R,x_0+R]}$. Uniqueness in Lemma~\ref{lem:wellposed} identifies it with the solution in the last sentence of Lemma~\ref{lem:amplification}. Substituting $x-x_0$ in \eqref{eq:amplification} proves the translated estimate with the same $T,R_*$.

\end{proof}

The convergence $\mu_T(\{|z|>R\})\to0$ used above keeps $T$ fixed. To treat growing time intervals in Proposition~\ref{prop:interior-rate}, we next make the dependence on time explicit. Estimate \eqref{eq:jump-moment} follows from subadditivity for powers at most one and Markov's inequality; see \cite[Chapter~1]{Durrett2019}.

\begin{lemma}\label{lem:jump-moment}
Let $0<p\le1$ and assume $m_p(J)<\infty$. Then, for every $t\ge0$ and $R>0$,
\begin{equation}\label{eq:jump-moment}
 \int_{\R}|z|^p\dd\mu_t(z)\le dt\,m_p(J),\qquad
 \mu_t(\{|z|>R\})\le dt\,m_p(J)R^{-p}.
\end{equation}
Under \eqref{eq:tail}, these bounds apply to any $0<p<\min\{1,\gamma-1\}$.
\end{lemma}

\begin{proof}

Recall the defined kernel integral $m_p(J)=\int_{\R}|z|^pJ(z)\dd z$ and the probability measure $\mu_t$ in \eqref{eq:mu}.

First, for $0<p<1$ and $s>0$,
\[
\frac{\dd }{\dd s}[(1+s)^p-1-s^p]
=p[(1+s)^{p-1}-s^{p-1}]\le0.
\]
The bracketed function has limit zero at $s=0$, so $(a+b)^p\le a^p+b^p$ for $a,b\ge0$, by scaling when $a>0$ and by direct evaluation when $a=0$. For $p=1$ equality holds. The triangle inequality and induction give
\[
|z_1+\cdots+z_n|^p\le\sum_{k=1}^n|z_k|^p.
\]
Using the product density for $n$ independent jumps and the convolution identity verified in Lemma~\ref{lem:semigroup}, Tonelli's theorem yields
\[
\begin{aligned}
\int_{\R}|z|^pJ^{*n}(z)\dd z
&=\int_{\R^n}|z_1+\cdots+z_n|^p\prod_{k=1}^nJ(z_k)\dd z_1\cdots\dd z_n\\
&\le\sum_{k=1}^n\left[\int_{\R}|z_k|^pJ(z_k)\dd z_k\right]
\prod_{j\ne k}\int_{\R}J(z_j)\dd z_j
=nm_p(J).
\end{aligned}
\]
All multiple integrands used before the inequality are nonnegative, so Tonelli applies even before finiteness is known. The finite right-hand side then proves integrability under the hypothesis $m_p(J)<\infty$.

We now integrate $|z|^p$ against $\mu_t$ using \eqref{eq:mu}. The atom at zero contributes zero because $p>0$. Tonelli for the nonnegative measure series \eqref{eq:mu} gives
\[
\begin{aligned}
\int_{\R}|z|^p\dd \mu_t(z)
&=\ee^{-dt}\sum_{n=1}^{\infty}\frac{(dt)^n}{n!}
\int_{\R}|z|^pJ^{*n}(z)\dd z\\
&\le m_p(J)\ee^{-dt}\sum_{n=1}^{\infty}\frac{(dt)^nn}{n!}\\
&=m_p(J)\ee^{-dt}dt\sum_{k=0}^{\infty}\frac{(dt)^k}{k!}
=dt\,m_p(J).
\end{aligned}
\]
The identity also holds at $t=0$, where $\mu_0=\delta_0$ and both sides vanish.

For $R>0$, $|z|^p\ge R^p$ on $\{|z|>R\}$. Integration gives
\[
R^p\mu_t(\{|z|>R\})
\le\int_{|z|>R}|z|^p\dd \mu_t(z)
\le\int_{\R}|z|^p\dd \mu_t(z)
\le dt\,m_p(J).
\]
Dividing by $R^p>0$ proves the second inequality of \eqref{eq:jump-moment}. This is Markov's inequality with its nonnegative random variable and finite expectation explicitly verified.

Equivalently, conditioning on $N_t=n$ and using subadditivity gives
$\mathbb E(|X_t|^p\mid N_t=n)\le nm_p(J)$. Independence of the jumps from $N_t$ justifies the unchanged product density. Summing over $n$ gives
\[
\mathbb E|X_t|^p\le m_p(J)\mathbb EN_t=dt\,m_p(J),\qquad
\mathbb P(|X_t|>R)\le R^{-p}\mathbb E|X_t|^p.
\]
Here $\mathbb EN_t=dt$ follows from the scalar series computed above. The law identity $\mathbb P(X_t\in\cdot)=\mu_t$ identifies these with the analytic inequalities.

Finally, set $A=\max\{1,R_J\}$. The compact contribution to $m_p(J)$ is bounded by $A^p\int_{|z|\le A}J(z)\dd z\le A^p$. On the complement, \eqref{eq:tail-explicit} gives
\[
\int_{|z|>A}|z|^pJ(z)\dd z
\le2j_2\int_A^\infty z^{p-\gamma}\dd z
=\frac{2j_2A^{p+1-\gamma}}{\gamma-p-1}<\infty
\quad\text{if }p<\gamma-1.
\]
Thus every $0<p<\min\{1,\gamma-1\}$ satisfies the hypotheses of \eqref{eq:jump-moment}.

\end{proof}

\section{Comparison functions}\label{sec:convolution}

For a positive bounded Borel function $v$, we estimate $|(J*v)(x)-v(x)|/v(x)$ under \eqref{eq:J} and \eqref{eq:J-h}. The comparison functions are defined in \eqref{eq:profiles}. Constants are independent of $x$ and $R$ in the stated ranges.

\subsection{Tail integrals}

For $R\ge1$, write
\begin{equation}\label{eq:tail-moments}
 T_J(R)=\int_{|z|>R}J(z)\dd z,
 \qquad V_k(R)=\int_{|z|\le R}|z|^kJ(z)\dd z\text{ for } k>0.
\end{equation}
We first estimate $T_J(R)$, $V_k(R)$ and $m_k(J)$.

\begin{lemma}\label{lem:kernel-moments}
Assume \eqref{eq:J} and \eqref{eq:tail}. For all sufficiently large $R$,
\begin{equation}\label{eq:tail-mass}
 T_J(R)\sim R^{1-\gamma}.
\end{equation}
For each fixed $k>0$ and all sufficiently large $R$,
\begin{equation}\label{eq:truncated-moment}
 V_k(R)\sim
 \begin{cases}
 R^{k+1-\gamma},&1<\gamma<k+1,\\
 \log R,&\gamma=k+1,\\
 1,&\gamma>k+1.
 \end{cases}
\end{equation}
Moreover,
\begin{equation}\label{eq:moment-exact-criterion}
 \int_{\R}|z|^kJ(z)\dd z<\infty
 \quad\Longleftrightarrow\quad \gamma>k+1.
\end{equation}
For every $a>0$, $\int_{\R}\ee^{a|z|}J(z)\dd z=\infty$, while
\begin{equation}\label{eq:log-moment}
 \int_{\R}J(z)\log(2+|z|)\dd z<\infty.
\end{equation}
\end{lemma}

\begin{proof}
By the evenness of $J$ in \eqref{eq:J} and the bounds in \eqref{eq:tail-explicit},
\[
 \frac{2j_1}{\gamma-1}R^{1-\gamma}
 \le T_J(R)=2\int_R^\infty J(z)\dd z
 \le\frac{2j_2}{\gamma-1}R^{1-\gamma}
 \quad\text{for }R\ge\max\{1,R_J\}.
\]
Here $\int_R^\infty z^{-\gamma}\dd z=R^{1-\gamma}/(\gamma-1)$ is finite because $\gamma>1$. This proves \eqref{eq:tail-mass}.

For fixed $k>0$, $\int_{|z|\le R_J}|z|^kJ(z)\dd z$ is finite, since
$|z|^k\le R_J^k$ on $[-R_J,R_J]$ and $\int_{\R}J=1$. By \eqref{eq:tail-explicit},
\[
 2j_1\int_{R_J}^R z^{k-\gamma}\dd z
 \le V_k(R)-\int_{|z|\le R_J}|z|^kJ(z)\dd z
 \le 2j_2\int_{R_J}^R z^{k-\gamma}\dd z
 \quad\text{for }R\ge\max\{1,R_J\}.
\]
If $k+1-\gamma\ne0$, the last power integral equals
\begin{equation}\label{eq:power-integral}
 \int_{R_J}^R z^{k-\gamma}\dd z
 =\frac{R^{k+1-\gamma}-R_J^{k+1-\gamma}}{k+1-\gamma}.
\end{equation}
At $\gamma=k+1$, $\int_{R_J}^Rz^{-1}\dd z=\log(R/R_J)$. If $\gamma<k+1$, \eqref{eq:power-integral} lies between positive constant multiples of $R^{k+1-\gamma}$ for all sufficiently large $R$. If $\gamma=k+1$, $\log(R/R_J)$ is comparable to $\log R$ for large $R$. In both cases the preceding two-sided bounds imply \eqref{eq:truncated-moment} and $V_k(R)\to\infty$.

If $\gamma>k+1$, then $\int_{R_J}^{\infty}z^{k-\gamma}\dd z<\infty$. Thus $V_k(R)$ increases to $m_k(J)<\infty$ by the monotone convergence theorem, applied to $|z|^kJ(z)\ind_{\{|z|\le R\}}$. Continuity and $J(0)>0$ give $a>0$ with $J(z)\ge J(0)/2$ for $|z|\le a$; consequently $m_k(J)\ge J(0)a^{k+1}/(k+1)>0$. Hence $V_k(R)$ is bounded above and below by positive constants for large $R$. This completes \eqref{eq:truncated-moment}; monotone convergence and divergence in the other two cases also prove \eqref{eq:moment-exact-criterion}.

For $a>0$, $\ee^{az}z^{-\gamma}\ge1$ for sufficiently large $z$, because $az-\gamma\log z\to\infty$. The lower bound in \eqref{eq:tail-explicit} then gives $\int_{\R}\ee^{a|z|}J(z)\dd z=\infty$. To prove \eqref{eq:log-moment}, put $\varepsilon=(\gamma-1)/2>0$. Since $\log(2+z)\le C_\varepsilon z^\varepsilon$ for $z\ge1$,
\[
 \int_{|z|>A}J(z)\log(2+|z|)\dd z
 \le 2j_2C_\varepsilon\int_A^\infty z^{-(\gamma+1)/2}\dd z<\infty,
 \qquad A=\max\{1,R_J\}.
\]
The integral on $[-A,A]$ is at most $\log(2+A)\int_{\R}J=\log(2+A)$. In particular, \eqref{eq:moment-exact-criterion} with $k=1$ and $k=2$ gives the respective thresholds $\gamma=2$ for $m_1(J)$ and $\gamma=3$ for $m_2(J)$.
\end{proof}

\subsection{Two comparison functions}

For $R\ge4$, define
\begin{equation}\label{eq:profiles}
 \Phi_R(x)=\left[1+\left(\frac{1+|x|}{R}\right)^{\gamma\sigma}\right]^{-1/\sigma},
 \qquad
 \Psi_R(x)=\min\left\{1,\left(\frac R{1+|x|}\right)^\gamma\right\}.
\end{equation}
Here $\sigma$ is fixed by \eqref{eq:F}. Both functions are positive and continuous on $\R$.

\begin{lemma}\label{lem:profile-properties}
For $R\ge4$, the profiles in \eqref{eq:profiles} satisfy
\begin{equation}\label{eq:profile-comparison}
 2^{-1/\sigma}\Psi_R\le\Phi_R\le\Psi_R,
 \qquad \int_{\R}\Psi_R(x)\dd x
 =2(R-1)+\frac{2R}{\gamma-1}.
\end{equation}
Moreover,
\begin{equation}\label{eq:profile-lip}
 \Lip(\log\Phi_R)+\Lip(\log\Psi_R)\le C R^{-\min\{1,\gamma\sigma\}}.
\end{equation}
For each choice $V_R\in\{\Phi_R,\Psi_R\}$,
$V_R(x)\sim R^\gamma(1+|x|)^{-\gamma}$ for $|x|\ge R$, uniformly in $R\ge4$ and $x$.
\end{lemma}

\begin{proof}
If $0<a\le1$, then $1\le1+a^\sigma\le2$, so
$2^{-1/\sigma}\le(1+a^\sigma)^{-1/\sigma}\le1$. If $a\ge1$, write
$(1+a^\sigma)^{-1/\sigma}=a^{-1}(1+a^{-\sigma})^{-1/\sigma}$ and apply the same bounds. Substituting $a=((1+|x|)/R)^\gamma$ proves the pointwise inequalities in \eqref{eq:profile-comparison}. By evenness and the two branches defining $\Psi_R$ in \eqref{eq:profiles},
\[
 \int_{\R}\Psi_R(x)\dd x
 =2\int_0^{R-1}1\dd x+2R^\gamma\int_{R-1}^{\infty}(1+x)^{-\gamma}\dd x.
\]
Both terms are finite since $R\ge4$ and $\gamma>1$. The tail integral on one side is
$R^\gamma\int_{R-1}^{\infty}(1+x)^{-\gamma}\dd x=R/(\gamma-1)$.

For $s>0$, logarithmic differentiation gives
\begin{equation}\label{eq:log-derivative}
 \left|\frac{\dd}{\dd s}\log\Phi_R(s)\right|
 =\frac{\gamma}{1+s}
   \frac{((1+s)/R)^{\gamma\sigma}}
        {1+((1+s)/R)^{\gamma\sigma}}.
\end{equation}
When $1+s\le R$, the right-hand side is bounded by
$\gamma R^{-\gamma\sigma}(1+s)^{\gamma\sigma-1}$. Its maximum for $1\le1+s\le R$ is at most $CR^{-\min\{1,\gamma\sigma\}}$. When $1+s\ge R$, it is at most $\gamma/R$. The derivative bound for $\log\Phi_R$ on the negative half-line is the same by evenness. For any $x<y$, split $[x,y]$ at zero when necessary. On each resulting interval the fundamental theorem of calculus bounds the increment of $\log\Phi_R$ by $CR^{-\min\{1,\gamma\sigma\}}$ times the interval length; continuity at zero permits addition of those increments. Hence
$|\log\Phi_R(y)-\log\Phi_R(x)|\le CR^{-\min\{1,\gamma\sigma\}}|y-x|$.

For $\log\Psi_R$, the derivative is zero on $(-R+1,R-1)$ and has absolute value $\gamma/(1+|x|)\le\gamma/R$ on either exterior interval. Splitting $[x,y]$ at $\pm(R-1)$ and using continuity at both endpoints gives
$|\log\Psi_R(y)-\log\Psi_R(x)|\le(\gamma/R)|y-x|$.
Since $R^{-1}\le R^{-\min\{1,\gamma\sigma\}}$ for $R\ge1$, these two inequalities prove \eqref{eq:profile-lip}.

For $|x|\ge R$, the minimum in \eqref{eq:profiles} gives $\Psi_R(x)=R^\gamma(1+|x|)^{-\gamma}$. By \eqref{eq:profile-comparison}, the tail comparison also holds for $\Phi_R$.
\end{proof}

Put
\begin{equation}\label{eq:basic-error}
 \rho=\frac12\min\{1,\gamma\sigma\},\qquad
 e_\gamma(R)=R^{-\rho}+R^{-(\gamma-1)\rho}\text{ for } R\ge4.
\end{equation}
Since $\rho>0$ and $\gamma>1$, $e_\gamma(R)\to0$ and $\int_4^\infty e_\gamma(R)\dd R/R<\infty$.
Indeed, by \eqref{eq:basic-error},
\[
 \int_4^\infty\frac{e_\gamma(R)}R\dd R
 =\frac{4^{-\rho}}\rho+
   \frac{4^{-(\gamma-1)\rho}}{(\gamma-1)\rho}<\infty.
\]

\begin{lemma}\label{lem:relative}
Under \eqref{eq:J} and \eqref{eq:J-h}, there is $C_0>0$ such that
\begin{equation}\label{eq:relative}
 |(J*V_R)(x)-V_R(x)|\le C_0e_\gamma(R)V_R(x)
 \text{ for } R\ge4,\ x\in\R
\end{equation}
for $V_R=\Phi_R$ and for $V_R=\Psi_R$.
\end{lemma}

\begin{proof}
Set $\ell=R^\rho$. Since $0<\rho\le1/2$ and $R\ge4$, we have $\ell\le R/2$. By \eqref{eq:profile-lip} and $2\rho=\min\{1,\gamma\sigma\}$,
\[
 |\log V_R(x-z)-\log V_R(x)|
 \le CR^{-2\rho}|z|\le CR^{-\rho}
 \text{ for } |z|\le\ell.
\]
The elementary inequality $|\ee^a-1|\le\ee^{|a|}|a|$ therefore yields
\begin{equation}\label{eq:short-ratio}
 \left|\frac{V_R(x-z)}{V_R(x)}-1\right|\le CR^{-\rho}
 \text{ for } |z|\le\ell,\ x\in\R.
\end{equation}
For fixed $R$ and $x$, the integrand below is absolutely integrable: since $0<V_R\le1$,
\[
 \int_{\R}J(z)\left|\frac{V_R(x-z)}{V_R(x)}-1\right|\dd z
 \le\frac1{V_R(x)}+1<\infty.
\]
Thus linearity of the integral and $\int_{\R}J=1$ give
\begin{equation}\label{eq:relative-formula}
 \frac{(J*V_R)(x)-V_R(x)}{V_R(x)}
 =\int_{\R}J(z)\left[\frac{V_R(x-z)}{V_R(x)}-1\right]\dd z.
\end{equation}
By \eqref{eq:short-ratio} and $\int_{|z|\le\ell}J(z)\dd z\le1$, the absolute value of the integral in \eqref{eq:relative-formula} restricted to $|z|\le\ell$ is at most $CR^{-\rho}$. For the other regions, \eqref{eq:J-h} gives the uniform bound
\[
 T_J(a)\le2C_J\int_a^\infty(1+z)^{-\gamma}\dd z
 =\frac{2C_J}{\gamma-1}(1+a)^{1-\gamma}
 \le Ca^{1-\gamma}\quad\text{for }a\ge1.
\]

If $|x|\le2R$, then $(1+|x|)/R\le9/4$ because $R\ge4$. Formula \eqref{eq:profiles} therefore gives
\[
 V_R(x)\ge c_*:=\min\bigl\{[1+(9/4)^{\gamma\sigma}]^{-1/\sigma},(4/9)^\gamma\bigr\}>0.
\]
Since $0<V_R\le1$, $|V_R(x-z)/V_R(x)-1|\le c_*^{-1}+1$, independently of $z$. Hence
\[
 \int_{|z|>\ell}J(z)
 \left|\frac{V_R(x-z)}{V_R(x)}-1\right|\dd z
 \le C\int_{|z|>\ell}h_\gamma(z)\dd z
 \le C\ell^{1-\gamma}.
\]
This proves \eqref{eq:relative} for $|x|\le2R$.

Now suppose $|x|>2R$. Split the remaining integral into
$\ell<|z|\le|x|/2$ and $|z|>|x|/2$. On the first region, $|x-z|\ge|x|/2>R$. By \eqref{eq:profile-comparison}, both $V_R(x-z)$ and $V_R(x)$ lie between fixed multiples of the respective functions $R^\gamma(1+|x-z|)^{-\gamma}$ and $R^\gamma(1+|x|)^{-\gamma}$. Thus
\[
 \frac{V_R(x-z)}{V_R(x)}
 \le C\left(\frac{1+|x|}{1+|x-z|}\right)^\gamma\le C2^\gamma.
\]
It follows that the absolute value of the integral in \eqref{eq:relative-formula} over $\ell<|z|\le|x|/2$ is at most $CT_J(\ell)\le C\ell^{1-\gamma}$.

On $|z|>|x|/2$, the contribution from the term $-1$ in \eqref{eq:relative-formula} is bounded by $T_J(|x|/2)\le C\ell^{1-\gamma}$, by \eqref{eq:J-h}. For the positive contribution, use
$J(z)\le C(1+|x|)^{-\gamma}$ and
$V_R(x)\ge cR^\gamma(1+|x|)^{-\gamma}$. The first inequality follows from \eqref{eq:J-h} and $1+|z|\ge(1+|x|)/2$; the second follows from \eqref{eq:profile-comparison} and $|x|>R$. With $y=x-z$, Lebesgue measure is unchanged. Extending the resulting nonnegative integral to all $y\in\R$, and using $\int_{\R}V_R\le\int_{\R}\Psi_R\le CR$ from \eqref{eq:profile-comparison}, we obtain
\begin{equation}\label{eq:remote-mass}
 \frac1{V_R(x)}\int_{|z|>|x|/2}J(z)V_R(x-z)\dd z
 \le CR^{-\gamma}\int_{\R}V_R(y)\dd y
 \le CR^{1-\gamma}.
\end{equation}
Since $R\ge\ell$ and $1-\gamma<0$, \eqref{eq:remote-mass} is at most $C\ell^{1-\gamma}$. Combining the three regions gives
$CR^{-\rho}+C\ell^{1-\gamma}=Ce_\gamma(R)$ and proves \eqref{eq:relative}. Estimate \eqref{eq:relative} uses the unit mass of $J$ and the upper bound in \eqref{eq:J-h}.
\end{proof}

\subsection{Estimate for the truncated function}

For $\Psi_R$ in \eqref{eq:profiles}, the proof below bounds $|(J*\Psi_R)(x)-\Psi_R(x)|/\Psi_R(x)$ using $V_1(R)/R+T_J(R)$. Define
\begin{equation}\label{eq:first-error}
 \eta_\gamma(R)=
 \begin{cases}
 R^{1-\gamma},&1<\gamma<2,\\
 R^{-1}\log R,&\gamma=2,\\
 R^{-1},&\gamma>2.
 \end{cases}
\end{equation}

\begin{proposition}\label{prop:truncated-error}
Under \eqref{eq:J} and \eqref{eq:tail}, there are $c,C,R_*>0$ such that
\begin{equation}\label{eq:truncated-error}
 c\eta_\gamma(R)\le
 \sup_{x\in\R}\frac{|(J*\Psi_R)(x)-\Psi_R(x)|}{\Psi_R(x)}
 \le C\eta_\gamma(R)\text{ for } R\ge R_*.
\end{equation}
The lower estimate holds at the junction $x=R-1$.
\end{proposition}

\begin{proof}
From \eqref{eq:profiles}, the derivative of $\Psi_R$ is zero in the central interval and has absolute value at most $\gamma/R$ on both exterior intervals. Continuity at $\pm(R-1)$, followed by splitting any segment at these two points, gives
$|\Psi_R(x-z)-\Psi_R(x)|\le\gamma|z|/R$ for all $x,z\in\R$.
Also $0\le\Psi_R\le1$, so $|\Psi_R(x-z)-\Psi_R(x)|\le1$. If $|x|\le2R$ and $R\ge4$, then $\Psi_R(x)\ge(4/9)^\gamma$. The identity \eqref{eq:relative-formula}, applied to $V_R=\Psi_R$, therefore yields
\[
 \frac{|(J*\Psi_R)(x)-\Psi_R(x)|}{\Psi_R(x)}
 \le(9/4)^\gamma\left[\frac\gamma R\int_{|z|\le R}|z|J(z)\dd z
             +\int_{|z|>R}J(z)\dd z\right].
\]
By \eqref{eq:tail-moments}, the right-hand side is at most a fixed multiple of $V_1(R)/R+T_J(R)$, giving
\begin{equation}\label{eq:first-core}
 \frac{|(J*\Psi_R)(x)-\Psi_R(x)|}{\Psi_R(x)}
 \le C\left[\frac{V_1(R)}R+T_J(R)\right]
 \le C\eta_\gamma(R).
\end{equation}
For the last inequality in \eqref{eq:first-core}, \eqref{eq:truncated-moment} with $k=1$ and \eqref{eq:tail-mass} give respectively
$V_1(R)/R\le CR^{1-\gamma}$ for $1<\gamma<2$,
$V_1(R)/R\le C(\log R)/R$ for $\gamma=2$, and
$V_1(R)/R\le C/R$ for $\gamma>2$.
The term $T_J(R)\le CR^{1-\gamma}$ is bounded by the same right-hand sides when $R$ is sufficiently large. These are precisely the three branches of \eqref{eq:first-error}.

Let $X=|x|>2R$. If $|z|\le X/2$ and $y$ lies between $x$ and $x-z$, then $X/2\le|y|\le3X/2$. Hence \eqref{eq:profiles} gives
\[
 \frac{|\Psi_R'(y)|}{\Psi_R(x)}
 =\frac{\gamma}{1+|y|}
   \left(\frac{1+X}{1+|y|}\right)^\gamma
 \le\frac{2^{\gamma+1}\gamma}{X}.
\]
The fundamental theorem of calculus on this interval then yields
\[
 \int_{|z|\le X/2}J(z)
 \frac{|\Psi_R(x-z)-\Psi_R(x)|}{\Psi_R(x)}\dd z
 \le C\frac{V_1(X/2)}X.
\]
By \eqref{eq:truncated-moment}, the last expression is bounded by
$CX^{1-\gamma}$, $C(\log X)/X$ or $C/X$, according as $1<\gamma<2$, $\gamma=2$ or $\gamma>2$.
For $R$ sufficiently large and $X>2R$, these quantities are at most $C\eta_\gamma(R)$: the power functions decrease, and $(\log X)/X$ decreases on $[\ee,\infty)$. On $|z|>X/2$, the term containing $-\Psi_R(x)$ contributes at most $T_J(X/2)$ to $|(J*\Psi_R)(x)-\Psi_R(x)|/\Psi_R(x)$. By \eqref{eq:remote-mass}, the term containing $\Psi_R(x-z)$ contributes at most $CR^{1-\gamma}$. Thus this region contributes at most
$C[X^{1-\gamma}+R^{1-\gamma}]\le C\eta_\gamma(R)$. This proves the upper inequality in \eqref{eq:truncated-error}.

For the lower estimate, take $x_R=R-1$. Then $\Psi_R(x_R)=1$ and $\Psi_R\le1$, so every integrand in
\[
 1-(J*\Psi_R)(x_R)
 =\int_{\R}J(z)[1-\Psi_R(x_R-z)]\dd z
\]
is nonnegative. Restrict the integral to $-R\le z\le0$ and change variables from $z$ to $-z$, using the evenness of $J$. For $0\le z\le R$,
\[
 \Psi_R(x_R+z)=(1+z/R)^{-\gamma},\qquad
 1-(1+z/R)^{-\gamma}\ge c_\gamma z/R.
\]
Indeed, for $0\le s\le1$,
\[
 1-(1+s)^{-\gamma}
 =\gamma\int_0^s(1+a)^{-\gamma-1}\dd a
 \ge\gamma2^{-\gamma-1}s.
\]
Thus $c_\gamma=\gamma2^{-\gamma-1}$ is admissible. By evenness of $J$ and \eqref{eq:tail-moments}, $\int_0^RzJ(z)\dd z=V_1(R)/2$, and hence
\begin{equation}\label{eq:first-lower}
 1-(J*\Psi_R)(x_R)
 \ge\frac{c_\gamma}{R}\int_0^R zJ(z)\dd z
 =\frac{c_\gamma V_1(R)}{2R}
 \ge c\eta_\gamma(R).
\end{equation}
The last inequality follows from \eqref{eq:truncated-moment} with $k=1$. Thus \eqref{eq:first-lower} proves the lower inequality in \eqref{eq:truncated-error} at $x=R-1$.
\end{proof}

By \eqref{eq:truncated-error}, $\sup_{x\in\R}|(J*\Psi_R)(x)-\Psi_R(x)|/\Psi_R(x)$ is comparable to $V_1(R)/R+T_J(R)$ as $R\to\infty$, even when $J$ is even. For $W_R$ in \eqref{eq:W-def}, cancellation in \eqref{eq:central-difference} yields \eqref{eq:smooth-error}.

\section{Proof of the level-set estimates}\label{sec:barriers}\label{sec:consequences}

Assumptions \eqref{eq:J}, \eqref{eq:tail}, \eqref{eq:data}, and \eqref{eq:F} hold throughout this section. For $V_R=\Phi_R$ and $V_R=\Psi_R$, estimate \eqref{eq:relative} bounds $|(J*V_R)(x)-V_R(x)|/V_R(x)$ by $C_0e_\gamma(R)$. We use \eqref{eq:relative} to choose $R_-(t)$ and $R_+(t)$ in \eqref{eq:Rpm}; the finiteness of $\int_4^\infty e_\gamma(R)\dd R/R$ keeps both radii comparable to $S_\gamma(t)$.

\subsection{Choice of the radii}

Fix $C_0$ from Lemma~\ref{lem:relative}. Choose $R_0\ge4$ so large that
\begin{equation}\label{eq:R0}
 \supp u_0\subset[-R_0+1,R_0-1],\qquad
 dC_0e_\gamma(R_0)\le r/2.
\end{equation}
Set
\begin{equation}\label{eq:epsilon-B}
 \nu=\frac r{2\gamma},\qquad
 \varepsilon(t)=dC_0e_\gamma(R_0\ee^{\nu t}),\qquad
 B(t)=\int_0^t\varepsilon(s)\dd s,
\end{equation}
and define
\begin{equation}\label{eq:Rpm}
 R_\pm(t)=R_0\exp\left(\frac{rt\pm B(t)}\gamma\right)
 \text{ for } t\ge0.
\end{equation}
By \eqref{eq:basic-error} and \eqref{eq:epsilon-B},
\begin{equation}\label{eq:B-explicit}
\begin{split}
 B(t)={}&\frac{dC_0R_0^{-\rho}}{\nu\rho}
                  (1-\ee^{-\nu\rho t})\\
       &+\frac{dC_0R_0^{-(\gamma-1)\rho}}{\nu(\gamma-1)\rho}
                  (1-\ee^{-\nu(\gamma-1)\rho t}).
\end{split}
\end{equation}
In particular, $B(t)$ increases to the finite constant
\begin{equation}\label{eq:B-infty}
 B_\infty=\frac{dC_0R_0^{-\rho}}{\nu\rho}
          +\frac{dC_0R_0^{-(\gamma-1)\rho}}{\nu(\gamma-1)\rho}.
\end{equation}

\begin{lemma}\label{lem:scales}
For every $t\ge0$,
\begin{equation}\label{eq:scale-control}
 0<\varepsilon(t)\le r/2,\qquad
 R_\pm(t)\ge R_0\ee^{\nu t},\qquad
 dC_0e_\gamma(R_\pm(t))\le\varepsilon(t).
\end{equation}
Both $R_-$ and $R_+$ are increasing. If $t\ge\tau\ge0$, then
\begin{equation}\label{eq:scale-shift}
 R_-(t-\tau)\ge\ee^{-r\tau/\gamma}R_-(t).
\end{equation}
Moreover,
\begin{equation}\label{eq:scale-comparison}
 R_0\ee^{-B_\infty/\gamma}S_\gamma(t)
 \le R_-(t)\le R_+(t)
 \le R_0\ee^{B_\infty/\gamma}S_\gamma(t).
\end{equation}
\end{lemma}

\begin{proof}
By \eqref{eq:R0} and the monotonicity of $e_\gamma$, we have $0<\varepsilon(t)\le r/2$ and $0\le B(t)\le rt/2$. Thus \eqref{eq:Rpm} gives $R_\pm(t)\ge R_0\ee^{\nu t}$. Since $e_\gamma$ is decreasing and $R_\pm(t)\ge R_0\ee^{\nu t}$, one has $dC_0e_\gamma(R_\pm(t))\le dC_0e_\gamma(R_0\ee^{\nu t})=\varepsilon(t)$, proving the last inequality in \eqref{eq:scale-control}.

By \eqref{eq:basic-error}, $e_\gamma(R)=R^{-\rho}+R^{-(\gamma-1)\rho}$ is positive and decreasing for $R\ge4$. The change of variables $R=R_0\ee^{\nu s}$ in \eqref{eq:epsilon-B} gives
\[
 \int_0^\infty\varepsilon(s)\dd s
 =\frac{dC_0}{\nu}\int_{R_0}^\infty
   \bigl(R^{-\rho}+R^{-(\gamma-1)\rho}\bigr)\frac{\dd R}{R}
 =B_\infty<\infty.
\]
Since $\rho>0$ and $(\gamma-1)\rho>0$, integrating $R^{-\rho-1}$ and $R^{-(\gamma-1)\rho-1}$ over $[R_0,\infty)$ gives \eqref{eq:B-infty}. In particular, the constants in \eqref{eq:scale-comparison} are independent of time.

Differentiation yields
\begin{equation}\label{eq:scale-derivative}
 \frac{\gamma R_\pm'(t)}{R_\pm(t)}=r\pm\varepsilon(t).
\end{equation}
In view of \eqref{eq:scale-control}, both radii are increasing and $\gamma R_-'/R_-\le r$. Integration on $[t-\tau,t]$ proves \eqref{eq:scale-shift}, while $0\le B(t)\le B_\infty$ in \eqref{eq:Rpm} proves \eqref{eq:scale-comparison}.

More explicitly, \eqref{eq:scale-derivative} gives
\[
 \log\frac{R_-(t)}{R_-(t-\tau)}
 =\frac1\gamma\int_{t-\tau}^t(r-\varepsilon(s))\dd s
 \le\frac{r\tau}{\gamma}.
\]
Exponentiation yields \eqref{eq:scale-shift}, including $\tau=0$ and $\tau=t$.
\end{proof}

\subsection{Upper and lower solutions}

\begin{lemma}\label{lem:upper}
The function
\begin{equation}\label{eq:upper-def}
 \overline u(t,x)=\Psi_{R_+(t)}(x)
 =\min\left\{1,\frac{R_+(t)^\gamma}{(1+|x|)^\gamma}\right\}
\end{equation}
is an upper solution of \eqref{eq:cauchy}, and
\begin{equation}\label{eq:upper-bound}
 u(t,x)\le\overline u(t,x)\text{ for } t\ge0,\ x\in\R.
\end{equation}
\end{lemma}

\begin{proof}
If $1+|x|>R_+(t)$, equations \eqref{eq:upper-def} and \eqref{eq:scale-derivative} give
$\overline u_t=(r+\varepsilon(t))\overline u$. Since $f(s)\le rs$ on $[0,1]$, Lemma~\ref{lem:relative} and \eqref{eq:scale-control} imply
\begin{align}\label{eq:upper-tail-inequality}
 \Lop\overline u(t,x)+f(\overline u(t,x))
 &\le[dC_0e_\gamma(R_+(t))+r]\overline u(t,x)\notag\\
 &\le[r+\varepsilon(t)]\overline u(t,x)
 =\overline u_t(t,x).
\end{align}
If $1+|x|<R_+(t)$, then $\overline u=1$ and $\overline u_t=0$. Since $\overline u\le1$ on the entire line and $J$ has unit mass, $J*\overline u\le1$, so
$\Lop\overline u+f(1)\le0=\overline u_t$ there as well.

Since $R_+$ is strictly increasing, for each fixed $x$ the equality $1+|x|=R_+(t)$ holds at most once. Consequently, $\partial_t\overline u(t,x)\ge\Lop\overline u(t,x)+f(\overline u(t,x))$ holds for almost every $t$ at each fixed $x$. Moreover,
\[
 |\min\{1,\ee^{a_1}b\}-\min\{1,\ee^{a_2}b\}|
 \le |a_1-a_2|\text{ for }a_1,a_2\in\R,\ b>0.
\]
For each $b>0$, the function $a\mapsto\min\{1,\ee^a b\}$ is continuous and piecewise differentiable, with derivative $\ee^ab\le1$ below $a=-\log b$ and derivative zero above that value. Integrating its derivative separately on the two intervals proves the displayed inequality, including pairs $a_1,a_2$ on opposite sides of $-\log b$.

Taking $a_1=\gamma\log R_+(t)$, $a_2=\gamma\log R_+(s)$ and $b=(1+|x|)^{-\gamma}$, and using \eqref{eq:scale-derivative} and $0<r+\varepsilon(t)\le3r/2$ from \eqref{eq:scale-control}, gives
\[
 \|\overline u(t,\cdot)-\overline u(s,\cdot)\|_\infty
 \le \frac{3r}{2}|t-s|\text{ for }s,t\ge0.
\]

At each fixed $x$, \eqref{eq:upper-def} is locally absolutely continuous in time, and $|\partial_t\overline u(t,x)|\le3r/2$ for almost every $t$ by \eqref{eq:scale-control}. By $u\in C^1([0,T];\mathcal B)$ from Lemma~\ref{lem:wellposed}, the constant function $g_T(t)=3r/2+\max_{0\le s\le T}\|u_t(s,\cdot)\|_\infty$ is an integrable bound on $|\partial_t\overline u|+|u_t|$. We have $\|\overline u(t,\cdot)-\overline u(s,\cdot)\|_\infty\le(3r/2)|t-s|$ for $s,t\ge0$, and $\partial_t\overline u(t,x)\ge\Lop\overline u(t,x)+f(\overline u(t,x))$ for almost every $t\ge0$ at each fixed $x\in\R$. Thus $\overline u$ satisfies the time-regularity and differential-inequality requirements of Lemma~\ref{lem:wellposed}. On $|x|\le R_0-1$, \eqref{eq:upper-def} equals one at time zero; outside this interval, \eqref{eq:R0} gives $u_0(x)=0$.

By \eqref{eq:R0} and $0\le u_0\le1$, we have $\overline u(0,\cdot)\ge u_0$. Lemma~\ref{lem:wellposed} now gives \eqref{eq:upper-bound}.
\end{proof}

\begin{lemma}\label{lem:lower}
There is $a_0\in(0,1)$ such that
\begin{equation}\label{eq:lower-def}
 \underline u(t,x)=a_0\Phi_{R_-(t)}(x)
\end{equation}
is a lower solution of the evolution equation in \eqref{eq:cauchy}, and
\begin{equation}\label{eq:lower-bound}
 \underline u(t,x)\le u(t+1,x)
 \text{ for } t\ge0,\ x\in\R.
\end{equation}
\end{lemma}

\begin{proof}
Let $b_0>0$ be given by Lemma~\ref{lem:seed}, so that $u(1,x)\ge b_0h_\gamma(x)$ for every $x\in\R$. Choose $a_0\in(0,1)$ so small that
\begin{equation}\label{eq:a0-choice}
 a_0R_0^\gamma\le b_0,\qquad
 (r/2)a_0^{-\sigma}\ge K.
\end{equation}

For example,
\[
 a_0=\frac12\min\left\{1,\frac{b_0}{R_0^\gamma},
                          \left(\frac r{2K}\right)^{1/\sigma}\right\}
\]
satisfies \eqref{eq:a0-choice}; every quantity in this minimum is positive by \eqref{eq:F} and \eqref{eq:seed}.
By \eqref{eq:profiles},
\begin{equation}\label{eq:Phi-R-identity}
 R\partial_R\Phi_R=\gamma\Phi_R(1-\Phi_R^\sigma).
\end{equation}
Writing $q=((1+|x|)/R)^{\gamma\sigma}$ in \eqref{eq:profiles}, direct differentiation gives $R\partial_R q=-\gamma\sigma q$ and $\Phi_R^\sigma=(1+q)^{-1}$. Consequently,
\[
 R\partial_R\Phi_R=\gamma q(1+q)^{-1/\sigma-1}
                  =\gamma\Phi_R(1-\Phi_R^\sigma),
\]
which verifies \eqref{eq:Phi-R-identity} at every $x$, including $x=0$.

Differentiating \eqref{eq:lower-def} and using \eqref{eq:Phi-R-identity} and \eqref{eq:scale-derivative} gives
\begin{equation}\label{eq:lower-time-identity}
 \underline u_t
 =(r-\varepsilon(t))\underline u
 -(r-\varepsilon(t))a_0^{-\sigma}\underline u^{1+\sigma}.
\end{equation}
Since $r-\varepsilon(t)\ge r/2$, equations \eqref{eq:a0-choice} and \eqref{eq:F} give
\begin{equation}\label{eq:lower-reaction}
 \underline u_t\le (r-\varepsilon(t))\underline u
                  -K\underline u^{1+\sigma}
 \le f(\underline u)-\varepsilon(t)\underline u.
\end{equation}
By linearity, Lemma~\ref{lem:relative}, and \eqref{eq:scale-control},
\[
 \Lop\underline u\ge-dC_0e_\gamma(R_-(t))\underline u
                   \ge-\varepsilon(t)\underline u.
\]
This proves $\underline u_t\le\Lop\underline u+f(\underline u)$.

To check time differentiability in the supremum norm, set $a=\log R$ and $Q(a,x)=\Phi_{\ee^a}(x)$. The identity \eqref{eq:Phi-R-identity} gives
\[
 \partial_aQ=\gamma Q(1-Q^\sigma),\qquad
 |\partial_a^2Q|
 =\gamma^2\bigl|Q(1-Q^\sigma)\bigl(1-(1+\sigma)Q^\sigma\bigr)\bigr|
 \le\gamma^2(2+\sigma).
\]
Taylor's formula in $a$, with this bound uniform in $x$, shows that $a\mapsto Q(a,\cdot)$ is continuously differentiable in the supremum norm. Since $t\mapsto\log R_-(t)$ is $C^1$ by \eqref{eq:scale-derivative}, the chain rule proves the asserted norm differentiability of \eqref{eq:lower-def}. Moreover, \eqref{eq:lower-time-identity} can be written as $(r-\varepsilon)\underline u(1-(\underline u/a_0)^\sigma)$; since $0<\underline u\le a_0$ and $0<r-\varepsilon\le r$, its values lie between zero and $ra_0$.

At time zero,
\[
 \underline u(0,x)\le a_0R_0^\gamma h_\gamma(x)
                    \le b_0h_\gamma(x)\le u(1,x).
\]
Applying Lemma~\ref{lem:wellposed} to $u(t+1,x)$ proves \eqref{eq:lower-bound}.
\end{proof}

\subsection{Proof of the level-set theorem}

\begin{lemma}\label{lem:global}
Under the assumptions of Theorem~\ref{thm:main}, estimate \eqref{eq:global-main} holds.
\end{lemma}

\begin{proof}
By \eqref{eq:upper-bound} and \eqref{eq:scale-comparison},
\[
 u(t,x)\le\min\left\{1,
  \frac{R_0^\gamma\ee^{B_\infty}\ee^{rt}}{(1+|x|)^\gamma}\right\}
 \le\max\{1,R_0^\gamma\ee^{B_\infty}\}H_\gamma(t,x).
\]
The bound $u(t,x)\le C H_\gamma(t,x)$ holds for all $t\ge0$ and $x\in\R$. For $t\ge1$, equations \eqref{eq:lower-bound}, \eqref{eq:profile-comparison}, and \eqref{eq:scale-comparison}, applied at time $t-1$, give
\[
 u(t,x)\ge a_0 2^{-1/\sigma}
  \min\left\{1,
   \frac{R_0^\gamma\ee^{-r-B_\infty}\ee^{rt}}{(1+|x|)^\gamma}\right\}.
\]
For any $k,z>0$, $\min\{1,kz\}\ge\min\{1,k\}\min\{1,z\}$. In view of \eqref{eq:scale-envelope}, \eqref{eq:global-main} holds with
\[
 C=\max\{1,R_0^\gamma\ee^{B_\infty}\},\qquad
 c=a_0 2^{-1/\sigma}\min\{1,R_0^\gamma\ee^{-r-B_\infty}\}.
\]
\end{proof}

\begin{lemma}\label{lem:all-levels}
For every $\theta\in(0,1)$, there are $c_\theta>0$ and $T_\theta>0$ such that
\begin{equation}\label{eq:interior-level}
 u(t,x)>\theta\text{ for } t\ge T_\theta,\ |x|\le c_\theta S_\gamma(t).
\end{equation}
\end{lemma}

\begin{proof}
Set $\delta=a_0 2^{-1/\sigma}$. By \eqref{eq:profiles}, $\Phi_{R_-(s)}(x)\ge2^{-1/\sigma}$ when $|x|\le R_-(s)-1$. Hence \eqref{eq:lower-bound} gives
\begin{equation}\label{eq:small-level}
 u(s+1,x)\ge\delta\text{ for } s\ge0,\ |x|\le R_-(s)-1.
\end{equation}
Choose $T,R_*>0$ from Lemma~\ref{lem:amplification} for $\delta$ and $\theta$. Once $R_-(s)-1\ge R_*$, \eqref{eq:small-level} and $u\ge0$ give $u(s+1,\cdot)\ge\delta\ind_{[-R_-(s)+1,R_-(s)-1]}$. Let $R=R_-(s)-1$. By Lemma~\ref{lem:wellposed}, $u(s+1+\tau,x)\ge v_R(\tau,x)$ for $\tau\ge0$, where $v_R$ is the solution with initial datum $\delta\ind_{[-R,R]}$ in Lemma~\ref{lem:amplification}. Evaluating at $\tau=T$ gives
\[
 u(s+1+T,x)>\theta
 \quad\text{for }|x|\le\tfrac12(R_-(s)-1)
\]
for all sufficiently large $s$. Because $R_-(s)\ge4$, this interval contains $[-R_-(s)/4,R_-(s)/4]$.

By \eqref{eq:scale-control}, $R_-(s)\ge R_0\ee^{\nu s}$, so one may choose
\[
 s_* =\max\left\{0,\frac1\nu\log\frac{R_*+1}{R_0}\right\}.
\]
For every $s\ge s_*$ the radius $R=R_-(s)-1$ meets the hypothesis $R\ge R_*$ of Lemma~\ref{lem:amplification}. The comparison with $v_R$ uses the autonomous equation \eqref{eq:cauchy}, the ordered initial functions at time $s+1$, and the bounded solution class in Lemma~\ref{lem:wellposed}. Thus the same $T$ works for all $s\ge s_*$.

Write $t=s+1+T$. The lower bound in \eqref{eq:scale-comparison} implies
\[
 \frac{R_-(s)}4\ge
 \frac{R_0}4\exp\left(-\frac{B_\infty+r(1+T)}\gamma\right)S_\gamma(t).
\]
Taking $c_\theta=(R_0/4)\exp[-(B_\infty+r(1+T))/\gamma]>0$ proves \eqref{eq:interior-level}.

Taking $T_\theta=s_*+1+T$ proves \eqref{eq:interior-level} for every $t\ge T_\theta$, rather than only along selected times.
\end{proof}

\begin{proof}[Proof of Theorem~\ref{thm:main}]
Lemma~\ref{lem:global} proves \textup{(i)}. For \textup{(ii)}, Lemma~\ref{lem:all-levels} gives the first inclusion in \eqref{eq:main-inclusions}. For the last inclusion, \eqref{eq:upper-bound} implies that $u(t,x)\ge\theta$ only if
\[
 (1+|x|)^\gamma\le\theta^{-1}R_+(t)^\gamma
 \le\theta^{-1}R_0^\gamma\ee^{B_\infty}\ee^{rt}.
\]
Thus we may take $C_\theta=\theta^{-1/\gamma}R_0\ee^{B_\infty/\gamma}$, increasing it if necessary.

For large $t$, \eqref{eq:interior-level} gives $u>\theta$ on an interval containing zero. By Lemma~\ref{lem:wellposed}, $u(t,\cdot)$ is continuous and tends to zero at both ends. The intermediate value theorem therefore gives a point of $E_\theta(t)$ on each side of that interval. More precisely, if $a=c_\theta S_\gamma(t)$, then $u(t,\pm a)>\theta$ by \eqref{eq:interior-level}. Since $u(t,x)\to0$ as $x\to\pm\infty$, choose $b_+>a$ and $b_-<-a$ with $u(t,b_\pm)<\theta$. Continuity on $[a,b_+]$ and $[b_-,-a]$ gives a point of $E_\theta(t)$ in each open interval.

The set $E_\theta(t)$ is closed by continuity of $u(t,\cdot)$ and bounded by \eqref{eq:main-inclusions}. Every $x\in E_\theta(t)$ satisfies
\begin{equation}\label{eq:all-points-annulus}
 c_\theta S_\gamma(t)<|x|\le C_\theta S_\gamma(t).
\end{equation}
By \eqref{eq:all-points-annulus}, estimate \eqref{eq:log-main} holds with $K_\theta=\max\{1,|\log c_\theta|,|\log C_\theta|\}$. Since $S_\gamma(t)/t\to\infty$, \eqref{eq:all-points-annulus} gives $|x_\theta^\pm(t)|/t\to\infty$ as $t\to\infty$.
Since $x_\theta^-(t)<0<x_\theta^+(t)$, \eqref{eq:all-points-annulus} also gives

\begin{equation}\label{eq:level-span}
 2c_\theta S_\gamma(t)
 <x_\theta^+(t)-x_\theta^-(t)
 \le2C_\theta S_\gamma(t)
 \text{ for } t\ge T_\theta.
\end{equation}

For $\theta\in[\theta_0,\theta_1]\subset(0,1)$, use $c=c_{\theta_1}$, $C=C_{\theta_0}$ and the maximum of the two endpoint starting times. For $\theta_0\le\theta\le\theta_1$, one has $u>\theta_1\Rightarrow u>\theta$ and $u\ge\theta\Rightarrow u\ge\theta_0$. Hence every $x\in E_\theta(t)$ satisfies $cS_\gamma(t)<|x|\le CS_\gamma(t)$, and $K=\max\{1,|\log c|,|\log C|\}$ is independent of $\theta$ in this interval.

To prove \textup{(iii)}, suppose that $a(t)/S_\gamma(t)\to0$. For each fixed $\theta\in(0,1)$, the interval $[-a(t),a(t)]$ is eventually contained in the first interval of \eqref{eq:main-inclusions}. Hence $u>\theta$ there.

For any $\eta\in(0,1)$, take $\theta=1-\eta$. Once $t\ge T_{1-\eta}$ and $a(t)\le c_{1-\eta}S_\gamma(t)$, the bound $u\le1$ gives $\sup_{|x|\le a(t)}|u(t,x)-1|\le\eta$. Letting $\eta\downarrow0$ proves the first limit in \eqref{eq:inner-outer}. For the second limit, the upper bound in \eqref{eq:global-main} gives
\[
 \sup_{|x|\ge b(t)}u(t,x)
 \le C\frac{S_\gamma(t)^\gamma}{(1+b(t))^\gamma}\to0.
\]
\end{proof}

\begin{remark}\label{rem:auxiliary-scales}
The comparison scales themselves have limits
\[
 \lim_{t\to\infty}\frac{R_\pm(t)}{S_\gamma(t)}
 =R_0\ee^{\pm B_\infty/\gamma}.
\]
The limits follow from \eqref{eq:Rpm} and $B(t)\to B_\infty$.
\end{remark}

We next derive arrival-time and integral estimates from \eqref{eq:global-main} and \eqref{eq:main-inclusions}.

\subsection{Arrival times}

For $0<\theta<1$ and $x\in\R$, define the first arrival time and the time after which $u(\cdot,x)$ stays at or above $\theta$ by
\begin{equation}\label{eq:arrival-def}
\begin{split}
 \tau_\theta(x)&=\inf\{t\ge0:u(t,x)\ge\theta\},\\
 \mathcal T_\theta(x)&=\inf\{t\ge0:u(s,x)\ge\theta\text{ for every }s\ge t\}.
\end{split}
\end{equation}
By Theorem~\ref{thm:main}\textup{(iii)}, $u(t,x)\to1$ for each fixed $x$. Thus $\tau_\theta(x)$ and $\mathcal T_\theta(x)$ in \eqref{eq:arrival-def} are finite, with $\tau_\theta(x)\le\mathcal T_\theta(x)$.

\begin{proposition}\label{prop:arrival}
Under the assumptions of Theorem~\ref{thm:main}, for every $\theta\in(0,1)$ there are $C_\theta,R_\theta>0$ such that
\begin{equation}\label{eq:arrival}
 \frac\gamma r\log|x|-C_\theta
 \le\tau_\theta(x)\le \mathcal T_\theta(x)
 \le\frac\gamma r\log|x|+C_\theta
 \text{ for } |x|\ge R_\theta.
\end{equation}
The constants can be chosen uniformly for $\theta$ in a compact subinterval of $(0,1)$.
\end{proposition}

\begin{proof}
The proof of Lemma~\ref{lem:global} gives an upper bound $u(t,x)\le C H_\gamma(t,x)$ for all $t\ge0$. If $u(t,x)\ge\theta$, then
\[
 \ee^{rt}\ge(\theta/C)(1+|x|)^\gamma,
 \qquad
 t\ge\frac\gamma r\log(1+|x|)-\frac1r\log(C/\theta).
\]
Taking the infimum over such $t$ gives the lower bound for $\tau_\theta$ in \eqref{eq:arrival}, since $\log(1+|x|)-\log|x|$ is bounded for $|x|\ge1$.

To prove the upper inequality in \eqref{eq:arrival}, take $c_\theta$ and $T_\theta$ from \eqref{eq:main-inclusions}. If
\[
 t\ge \max\left\{T_\theta,\frac\gamma r\log\frac{|x|}{c_\theta}\right\},
\]
then $|x|\le c_\theta S_\gamma(s)$ for every $s\ge t$, because $S_\gamma$ is increasing. The first inclusion in \eqref{eq:main-inclusions} therefore gives $u(s,x)>\theta$ for every $s\ge t$. This proves the upper inequality for $\mathcal T_\theta(x)$ in \eqref{eq:arrival} when $|x|$ is sufficiently large. Uniformity in the level follows by using the uniform constants in Theorem~\ref{thm:main}\textup{(ii)}.

For example, once $|x|\ge\max\{1,c_\theta\ee^{rT_\theta/\gamma}\}$, the maximum in the preceding display equals $(\gamma/r)\log(|x|/c_\theta)$. Thus one may enlarge $C_\theta$ to dominate both $r^{-1}|\log(C/\theta)|$ and $(\gamma/r)|\log c_\theta|$. For $\theta\in[\theta_0,\theta_1]\subset(0,1)$, use $\theta_0$ in the lower-time estimate and the common lower-inclusion constants for $\theta_1$ in the upper-time estimate. Both the resulting $C_\theta$ and $R_\theta$ are then independent of $\theta$ in this interval.
\end{proof}

By \eqref{eq:arrival}, $\mathcal T_\theta(x)-\tau_\theta(x)$ is bounded as $|x|\to\infty$, and $\tau_\theta(x)=(\gamma/r)\log|x|+O_\theta(1)$ and $\mathcal T_\theta(x)=(\gamma/r)\log|x|+O_\theta(1)$ as $|x|\to\infty$.

\subsection{Integral estimates}

We prove Corollary~\ref{cor:moments} and then describe the growth of truncated integrals when $M_{p,q}(t)$ is infinite.

\begin{lemma}\label{lem:envelope-integrals}
Let $S\ge1$ and $H_S(x)=\min\{1,(S/(1+|x|))^\gamma\}$. If $p\ge0$, $q>0$, and $p+1<\gamma q$, then
\begin{equation}\label{eq:envelope-integrals}
 c_{p,q}S^{p+1}
 \le\int_{\R}|x|^pH_S(x)^q\dd x
 \le C_{p,q}S^{p+1}
\end{equation}
with constants independent of $S$.
\end{lemma}

\begin{proof}
For the upper bound, use $H_S\le1$ on $|x|\le S$ and $H_S(x)\le(S/|x|)^\gamma$ on $|x|>S$. This yields
\begin{equation}\label{eq:envelope-upper-integral}
\begin{split}
 \int_{\R}|x|^pH_S(x)^q\dd x
 &\le2\int_0^Sx^p\dd x
     +2S^{\gamma q}\int_S^\infty x^{p-\gamma q}\dd x\\
 &=2S^{p+1}\left(\frac1{p+1}+\frac1{\gamma q-p-1}\right).
\end{split}
\end{equation}
The integral $\int_S^\infty x^{p-\gamma q}\dd x$ is finite because $p+1<\gamma q$.

For the lower bound, integrate only on $S/4\le|x|\le S/2$. Since $S\ge1$, $1+|x|\le3S/2$ on this region. Hence $H_S(x)^q\ge(2/3)^{\gamma q}$. The integral of $|x|^p$ over this region is a positive constant times $S^{p+1}$, proving the result.

In fact the lower integral is at least
\[
 2(2/3)^{\gamma q}\int_{S/4}^{S/2}x^p\dd x
 =\frac{2(2/3)^{\gamma q}}{p+1}
   \bigl(2^{-(p+1)}-4^{-(p+1)}\bigr)S^{p+1}>0.
\]
The upper integral in \eqref{eq:envelope-upper-integral} follows by splitting the nonnegative integrand into the two indicated regions; $p\ge0$ ensures integrability at zero and $p-\gamma q<-1$ ensures integrability at infinity.
\end{proof}

\begin{proof}[Proof of Corollary~\ref{cor:moments}]
For a fixed $t>0$, Lemma~\ref{lem:seed} bounds $u(t,x)$ above and below by positive constant multiples of $(1+|x|)^{-\gamma}$. The integrand in \eqref{eq:moment-def} is therefore comparable, for $|x|\ge1$, to $|x|^{p-\gamma q}$. The power-integral test proves \eqref{eq:moment-criterion}, including divergence at equality. The integrand is locally integrable at the origin because $p\ge0$ and $u$ is bounded.

For completeness, fix $0<\tau<t<T$ and take $b=b_{\tau,T}$, $B=B_{\tau,T}$ from \eqref{eq:finite-tail}. For $|x|\ge1$,
\[
 b^q2^{-\gamma q}|x|^{p-\gamma q}
 \le |x|^pu(t,x)^q
 \le B^q|x|^{p-\gamma q}.
\]
The integral of $x^{p-\gamma q}$ on $[1,\infty)$ is finite exactly when $p-\gamma q<-1$; at equality its integral up to $L$ is $\log L$. Also $\int_{-1}^1|x|^pu(t,x)^q\dd x\le2/(p+1)$. These bounds establish both directions of \eqref{eq:moment-criterion}.

For $t\ge1$, raise both sides of \eqref{eq:global-main} to the positive power $q$ and integrate against $|x|^p$. Lemma~\ref{lem:envelope-integrals} with $S=S_\gamma(t)$ gives \eqref{eq:moment-growth}.

The constants in \eqref{eq:global-main} do not depend on $t\ge1$. Consequently the integrated inequalities are $c^q\int_{\R}|x|^pH_\gamma(t,x)^q\dd x\le M_{p,q}(t)\le C^q\int_{\R}|x|^pH_\gamma(t,x)^q\dd x$; Lemma~\ref{lem:envelope-integrals} applies because $S_\gamma(t)\ge1$, and $S_\gamma(t)^{p+1}=\ee^{(p+1)rt/\gamma}$. The assertions about $\mathcal M(t)$, $M_{1,1}(t)$ and $M_{2,1}(t)$ follow by taking $(p,q)=(0,1),(1,1),(2,1)$, respectively.
\end{proof}

\begin{proposition}\label{prop:truncated-solution-moments}
Under the assumptions of Theorem~\ref{thm:main}, let $p\ge0$ and $q>0$. For a fixed $t>0$, define
$M_{p,q}(t;L)=\int_{-L}^L|x|^pu(t,x)^q\dd x$. If $p+1\ge\gamma q$, then, as $L\to\infty$,
\begin{equation}\label{eq:truncated-solution-moments}
 M_{p,q}(t;L)\sim
 \begin{cases}
 \log L,&p+1=\gamma q,\\
 L^{p+1-\gamma q},&p+1>\gamma q.
 \end{cases}
\end{equation}
The constants are uniform for $t$ in any interval $[\tau,T]$ with $0<\tau<T<\infty$.
\end{proposition}

\begin{proof}
Use \eqref{eq:finite-tail} on $|x|\ge1$ and integrate $x^{p-\gamma q}$ from $1$ to $L$. At equality the integral is $\log L$; above equality it is $(L^{p+1-\gamma q}-1)/(p+1-\gamma q)$. The contribution from $[-1,1]$ is bounded, uniformly in $t$ because $u\le1$. These power integrals give \eqref{eq:truncated-solution-moments}, uniformly for $\tau\le t\le T$.

Indeed, for $\alpha=p-\gamma q\ge-1$ and $L\ge1$, \eqref{eq:finite-tail} gives
\[
 2^{1-\gamma q}b_{\tau,T}^q\int_1^Lx^\alpha\dd x
 \le M_{p,q}(t;L)
 \le\frac2{p+1}+2B_{\tau,T}^q\int_1^Lx^\alpha\dd x
 \quad\text{for }\tau\le t\le T.
\]
For $L\ge2$, $\log L\ge\log2$ when $\alpha=-1$. If $\alpha>-1$, then $L^{\alpha+1}-1\ge(1-2^{-(\alpha+1)})L^{\alpha+1}$. The additive $2/(p+1)$ is absorbed into a fixed multiple of the corresponding positive function of $L$. This proves precisely the two-sided comparison denoted by $\sim$ in \eqref{eq:truncated-solution-moments}, uniformly on $[\tau,T]$.
\end{proof}

Thus, at each positive time, $M_{1,1}(t;L)$ grows like $\log L$ when $\gamma=2$, and $M_{2,1}(t;L)$ grows like $\log L$ when $\gamma=3$.

For $1\le q<\infty$, Corollary~\ref{cor:moments} gives
\begin{equation}\label{eq:Lq}
 \|u(t,\cdot)\|_{L^q(\R)}\sim\ee^{rt/(q\gamma)}
 \text{ for } t\ge1.
\end{equation}
For $q=\infty$, the bound $u\le1$ and the first limit in \eqref{eq:inner-outer} give $\|u(t,\cdot)\|_{L^\infty(\R)}\to1$ as $t\to\infty$. For $0<q<1$ and each fixed $t>0$, $\int_{\R}u(t,x)^q\dd x<\infty$ precisely when $q>1/\gamma$, by \eqref{eq:moment-criterion} with $p=0$.

\begin{lemma}\label{lem:mass}
The mass $\mathcal M$ defined after \eqref{eq:moment-def} belongs to $C^1([0,\infty))$ and satisfies
\begin{equation}\label{eq:mass-equation}
 \mathcal M'(t)=\int_{\R}f(u(t,x))\dd x,
 \qquad 0<\mathcal M'(t)\le r\mathcal M(t)\text{ for } t>0.
\end{equation}
In particular, $\mathcal M(t)$ is strictly increasing for $t>0$ and
\begin{equation}\label{eq:mass-order}
 c\ee^{rt/\gamma}\le\mathcal M(t)\le C\ee^{rt/\gamma}
 \text{ for } t\ge1.
\end{equation}
\end{lemma}

\begin{proof}
By Lemma~\ref{lem:weight}, $u\in C^1([0,\infty);\Yg)$, and integration is a bounded linear functional on $\Yg$ because
\[
 \left|\int_{\R}v(x)\dd x\right|
 \le\|h_\gamma\|_1\|v\|_\gamma.
\]
A bounded linear functional commutes with a norm derivative because it is continuous and linear. Applying this fact to $v\mapsto\int_{\R}v(x)\dd x$ gives $\mathcal M\in C^1([0,\infty))$ and $\mathcal M'(t)=\int_{\R}u_t(t,x)\dd x$. Moreover, $u(t,\cdot)\in L^1(\R)$, so the integral-preservation identity for convolution in Lemma~\ref{lem:semigroup} gives $\int_{\R}(\Lop u(t,\cdot))(x)\dd x=0$. Integrating \eqref{eq:cauchy} therefore proves the identity in \eqref{eq:mass-equation}.

By \eqref{eq:F}, $0\le f(u)\le ru$, proving finiteness and the upper bound for $\mathcal M'(t)$. For $t>0$, Lemma~\ref{lem:seed} gives $0<u(t,x)<1$ at every $x$. Thus $f(u(t,x))>0$, and spatial continuity gives an interval on which its integral is strictly positive. This proves $\mathcal M'(t)>0$. Finally, \eqref{eq:mass-order} is Corollary~\ref{cor:moments} with $(p,q)=(0,1)$.
\end{proof}

\subsection{Convergence and travelling solutions}

Theorem~\ref{thm:main}\textup{(iii)} gives convergence on $|x|\le a(t)$ whenever $a(t)/S_\gamma(t)\to0$. We now estimate its rate on $|x|\le e^{ct}$, $0\le c<r/\gamma$, under the additional condition $f'(1)<0$.

\begin{proposition}\label{prop:interior-rate}
Assume the hypotheses of Theorem~\ref{thm:main}. For every $c>r/\gamma$,
\begin{equation}\label{eq:outer-rate}
 \sup_{|x|\ge\ee^{ct}}u(t,x)
 \le C\ee^{-(\gamma c-r)t}\text{ for } t\ge1.
\end{equation}
If in addition $f'(1)<0$, then for every $0\le c<r/\gamma$ there exist $C_c,\kappa_c>0$ such that
\begin{equation}\label{eq:inner-rate}
 \sup_{|x|\le\ee^{ct}}(1-u(t,x))
 \le C_c\ee^{-\kappa_ct}\text{ for } t\ge1.
\end{equation}
\end{proposition}

\begin{proof}
For $|x|\ge\ee^{ct}$, \eqref{eq:global-main} gives $u(t,x)\le C\ee^{rt}(1+|x|)^{-\gamma}\le C\ee^{-(\gamma c-r)t}$, proving \eqref{eq:outer-rate}. For the inner estimate, assume $f'(1)<0$ and choose
$0<p<\min\{1,\gamma-1\}$, so that Lemma~\ref{lem:jump-moment} applies. Let $L=\max_{[0,1]}|f'|$. Choose a number $a\in(0,1)$ sufficiently small that
\begin{equation}\label{eq:inner-parameters}
 c<\frac{r(1-a)}\gamma,
 \qquad
 b:=\frac{pr(1-a)}\gamma-La>0.
\end{equation}
A number $a$ satisfying \eqref{eq:inner-parameters} exists because $c<r/\gamma$ and $p>0$.

For instance, putting $A=pr/\gamma>0$, it suffices to choose
\[
 0<a<\frac12\min\left\{1-\frac{\gamma c}{r},\frac{A}{A+L}\right\}.
\]
Then $r(1-a)/\gamma>c$ and $b=A-(A+L)a>0$, as required in \eqref{eq:inner-parameters}. Also \eqref{eq:tail-explicit} gives $m_p(J)<\infty$, since $p<\gamma-1$ makes $\int_{R_J}^\infty z^{p-\gamma}\dd z$ finite and $J\in L^\infty(\R)$ handles the remaining bounded interval.

Let $c_0>0$ be the lower constant in \eqref{eq:global-main}, and fix $\delta\in(0,1)$ with $\delta\le c_0$. At time $s_t=(1-a)t$, define
$R_t=\tfrac12\ee^{rs_t/\gamma}$. For large $t$, $1+R_t\le\ee^{rs_t/\gamma}$, and therefore \eqref{eq:global-main} gives
\[
 u(s_t,x)\ge\delta\text{ for } |x|\le R_t.
\]
Let $y$ solve $y'=f(y)$ with $y(0)=\delta$. Let $v_{R_t}$ solve \eqref{eq:cauchy} with initial datum $\delta\ind_{[-R_t,R_t]}$. By Lemma~\ref{lem:wellposed}, $v_{R_t}(\tau,x)\le u(s_t+\tau,x)$ for $\tau\ge0$. Applying \eqref{eq:stability} to $v_{R_t}$ and $y$ at time $\tau=at$, and then Lemma~\ref{lem:jump-moment}, gives, for $|x|\le R_t/2$,
\begin{align}\label{eq:inner-growing-time}
 u(t,x)&\ge y(at)-\delta\ee^{Lat}\mu_{at}(\{|z|>R_t/2\})\notag\\
       &\ge y(at)-C t\exp\left(La t-\frac{pr(1-a)}\gamma t\right)
        =y(at)-Ct\ee^{-bt}.
\end{align}
For $|x|\le R_t/2$, estimate \eqref{eq:stability} gives
\[
 |v_{R_t}(at,x)-y(at)|\le\delta\ee^{Lat}\mu_{at}(\{|z|>R_t/2\}),
\]
since $|x-z|>R_t$ implies $|z|>R_t/2$. By \eqref{eq:jump-moment} and $R_t=\tfrac12\ee^{r(1-a)t/\gamma}$, the right-hand side is at most
\[
 4^p da\delta m_p(J)\,t
 \exp\left(\left[La-\frac{pr(1-a)}\gamma\right]t\right).
\]
Here $d$ is the dispersal rate in \eqref{eq:cauchy}, $a$ is fixed in \eqref{eq:inner-parameters}, and $m_p(J)=\int_{\R}|z|^pJ(z)\dd z$. Thus the prefactor $4^p\delta\,d\,a\,m_p(J)$ is independent of $t$, as required in \eqref{eq:inner-growing-time}. By the first inequality in \eqref{eq:inner-parameters}, $\ee^{ct}\le R_t/2$ for large $t$.

Choose $0<\eta<-f'(1)$. Differentiability of $f$ at one gives $f(1-z)/z\to-f'(1)>\eta$ as $z\downarrow0$. Since $y(\tau)\to1$ by the proof of Lemma~\ref{lem:amplification}, choose $\tau_0$ so that $f(y(\tau))\ge\eta[1-y(\tau)]$ for $\tau\ge\tau_0$. Multiplying $(1-y)'\le-\eta(1-y)$ by $\ee^{\eta\tau}$ and integrating yields $1-y(\tau)\le(1-y(\tau_0))\ee^{-\eta(\tau-\tau_0)}$ there; boundedness of $1-y$ handles $0\le\tau\le\tau_0$. Also $t\ee^{-bt}\le C_b\ee^{-bt/2}$ because $\sup_{t\ge0}t\ee^{-bt/2}=2/(b\ee)<\infty$. Choose $t_0\ge1$ so that \eqref{eq:inner-growing-time} holds and $\ee^{ct}\le R_t/2$ for $t\ge t_0$. Then \eqref{eq:inner-growing-time} implies
\[
 1-u(t,x)\le C\ee^{-\eta a t}+Ct\ee^{-bt}
 \le C\ee^{-\kappa_ct}
\]
with, for instance, $\kappa_c=\min\{\eta a,b/2\}>0$. On the remaining interval $1\le t\le t_0$, the bound $0\le1-u\le1$ shows that choosing $C_c\ge\ee^{\kappa_ct_0}$ extends \eqref{eq:inner-rate} to all $t\ge1$.
\end{proof}

\paragraph{Travelling solutions.}

\begin{corollary}\label{cor:no-wave}
Under \eqref{eq:J}, \eqref{eq:tail}, and \eqref{eq:F}, the evolution equation in \eqref{eq:cauchy} has no solution of the form $V(x-ct)$ with $c\in\R$, continuous profile $0\le V\le1$, and
\begin{equation}\label{eq:wave-limits}
 \lim_{z\to-\infty}V(z)=1,\qquad
 \lim_{z\to+\infty}V(z)=0.
\end{equation}
\end{corollary}

\begin{proof}
For every $a>0$, evenness and \eqref{eq:tail-explicit} give
\[
 \int_{\R}\ee^{-az}J(z)\dd z
 \ge j_1\int_{R_J}^{\infty}\ee^{ay}y^{-\gamma}\dd y=\infty,
\]
since $\ee^{ay}y^{-\gamma}\ge\ee^{ay/2}$ for large $y$.
After $\tau=dt$, the equation is
$\partial_\tau u=J*u-u+(f/d)(u)$.
By Lemma~\ref{lem:wellposed}, the tangent-line extension of $f/d$
is globally Lipschitz and $C^1$, with derivative $r/d>0$ at zero.
Suppose that $V(x-ct)$ is such a solution, and set $\psi(z)=V(-z)$.
Since $J$ is even, spatial reflection gives the solution
$\psi(x+(c/d)\tau)$ of the rescaled equation,
with $\psi(-\infty)=0$ and $\psi(+\infty)=1$.
Theorem~2 of \cite{Yagisita2009} excludes this travelling wave.
\end{proof}

\section{Moment thresholds and convolution errors}\label{sec:ranges}

Assume \eqref{eq:J} and \eqref{eq:tail} throughout this section; additional hypotheses are stated where used. The threshold $\gamma=2$ for $m_1(J)<\infty$ enters \eqref{eq:truncated-error}. The threshold $\gamma=3$ for $m_2(J)<\infty$ enters \eqref{eq:smooth-error} and the Fourier estimates in \eqref{eq:symbol-orders}.

\subsection{A smooth comparison function}

Recall $W_R$ from \eqref{eq:W-def}. For every $R\ge2$, $W_R$ is positive, smooth and even, with $W_R(0)=1$. We first estimate $W_R$ and its first two derivatives.

\begin{lemma}\label{lem:smooth-properties}
There are positive constants depending only on $\gamma$ such that
\begin{equation}\label{eq:smooth-comparison}
 c\left(\frac R{R+|x|}\right)^\gamma
 \le W_R(x)\le
 C\left(\frac R{R+|x|}\right)^\gamma,
 \qquad \int_{\R}W_R(x)\dd x=c_\gamma R,
\end{equation}
where $c_\gamma=\int_{\R}(1+y^2)^{-\gamma/2}\dd y<\infty$. Moreover,
\begin{equation}\label{eq:smooth-derivatives}
 \frac{W_R'(x)}{W_R(x)}=-\frac{\gamma x}{R^2+x^2},
 \qquad
 \frac{W_R''(x)}{W_R(x)}
 =\frac{\gamma[(\gamma+1)x^2-R^2]}{(R^2+x^2)^2}.
\end{equation}
In particular,
\begin{equation}\label{eq:smooth-second-bound}
 |W_R''(x)|\le\frac{C W_R(x)}{(R+|x|)^2},
 \qquad \|W_R''\|_\infty\le C R^{-2}.
\end{equation}
\end{lemma}

\begin{proof}
The inequalities
$(R+|x|)/\sqrt2\le\sqrt{R^2+x^2}\le R+|x|$
give the first comparison in \eqref{eq:smooth-comparison}.

Indeed, $2(R^2+x^2)-(R+|x|)^2=(R-|x|)^2\ge0$ and $(R+|x|)^2-(R^2+x^2)=2R|x|\ge0$. Raising the resulting inequalities to the power $-\gamma$ and multiplying by $R^\gamma$ gives the claimed bounds for $W_R(x)$.

To compute the integral in \eqref{eq:smooth-comparison}, substitute $x=Ry$ on $[-A,A]$:
\[
 \int_{-A}^{A}W_R(x)\dd x
 =R\int_{-A/R}^{A/R}(1+y^2)^{-\gamma/2}\dd y.
\]
Monotone convergence applies as $A\to\infty$ because the integrands are nonnegative; $c_\gamma\le2+2/(\gamma-1)$. Differentiating $\log W_R$ gives
\[
 (\log W_R)'=-\frac{\gamma x}{R^2+x^2},\qquad
 (\log W_R)''=-\frac{\gamma(R^2-x^2)}{(R^2+x^2)^2}.
\]
The identity $W_R''/W_R=(\log W_R)''+[(\log W_R)']^2$ proves \eqref{eq:smooth-derivatives}. Since $R^2+x^2\ge(R+|x|)^2/2$, the bound $|\gamma[(\gamma+1)x^2-R^2]|\le C(R^2+x^2)$ gives the pointwise bound in \eqref{eq:smooth-second-bound}. Since $W_R\le1$ and $R+|x|\ge R$, it also gives $\|W_R''\|_\infty\le CR^{-2}$.

\end{proof}

Evenness of $J$ permits cancellation of the first-order variation even when $\int_{\R}|z|J(z)\dd z=\infty$. Indeed, for every bounded Borel function $W$,
\begin{equation}\label{eq:central-difference}
 (J*W)(x)-W(x)
 =\frac12\int_{\R}J(z)
 [W(x+z)+W(x-z)-2W(x)]\dd z.
\end{equation}
Since $J\in L^1(\R)$ and $W$ is bounded, every term in \eqref{eq:central-difference} is absolutely integrable, including when $1<\gamma\le2$.

More explicitly, $\int_{\R}J(z)|W(x\pm z)|\dd z\le\|W\|_\infty$. The change of variables $z\mapsto-z$ and the evenness of $J$ therefore equate the two integrals containing $W(x+z)$ and $W(x-z)$; averaging proves \eqref{eq:central-difference}. For $W=W_R$, twice applying the fundamental theorem of calculus gives
\[
 W_R(x+z)+W_R(x-z)-2W_R(x)
 =z^2\int_0^1(1-s)\bigl[W_R''(x+sz)+W_R''(x-sz)\bigr]\dd s.
\]
This identity uses derivatives of $W_R$ on a finite interval; no integral of $|z|J(z)$ is introduced.

\begin{lemma}\label{lem:smooth-upper}
Under \eqref{eq:J} and \eqref{eq:tail},
\begin{equation}\label{eq:smooth-upper}
 |(J*W_R)(x)-W_R(x)|\le C\omega_\gamma(R)W_R(x)
 \text{ for } R\ge2,\ x\in\R,
\end{equation}
where $\omega_\gamma$ is defined in \eqref{eq:omega-def}.
\end{lemma}

\begin{proof}
Recall $V_2$ and $T_J$ from \eqref{eq:tail-moments}.
Since $W_R$ is bounded and $J\in L^1(\R)$, the central differences
in \eqref{eq:central-difference} are absolutely integrable.
For $|x|\le2R$, one has $W_R(x)\ge5^{-\gamma/2}$.
Taylor's formula and \eqref{eq:smooth-second-bound} bound
$|W_R(x+z)+W_R(x-z)-2W_R(x)|$ by $CR^{-2}z^2$
on $|z|\le R$ and by four elsewhere. Hence
\begin{equation}\label{eq:smooth-core}
 \frac{|(J*W_R)(x)-W_R(x)|}{W_R(x)}
 \le C\left[R^{-2}V_2(R)+T_J(R)\right]
 \le C\omega_\gamma(R)
\end{equation}
for sufficiently large $R$, by \eqref{eq:truncated-moment}
and \eqref{eq:tail-mass}.

Now let $X=|x|>2R$. For $|z|\le X/2$ and $0\le s\le1$,
\[
 |x\pm sz|\ge X/2,\qquad
 \frac{W_R(x\pm sz)}{W_R(x)}
 =\left(\frac{R^2+X^2}{R^2+|x\pm sz|^2}\right)^{\gamma/2}
 \le2^\gamma.
\]
By \eqref{eq:smooth-second-bound}, $|W_R''(x\pm sz)|\le CX^{-2}W_R(x)$ for $|z|\le X/2$ and $0\le s\le1$. Taylor's formula and \eqref{eq:central-difference} therefore bound the contribution of $|z|\le X/2$, divided by $W_R(x)$, by $CX^{-2}V_2(X/2)\le C\omega_\gamma(R)$.
Here \eqref{eq:truncated-moment} gives the bound by
$C\omega_\gamma(X)$, and each branch of $\omega_\gamma$
decreases for sufficiently large arguments.

On $|z|>X/2$, the term containing $-2W_R(x)$ contributes
at most $CT_J(X/2)\le CX^{1-\gamma}$.
By \eqref{eq:J-h} and \eqref{eq:smooth-comparison},
\begin{align}\label{eq:smooth-remote}
 \frac1{W_R(x)}\int_{|z|>X/2}J(z)W_R(x-z)\dd z
 &\le \frac{CX^{-\gamma}}{R^\gamma X^{-\gamma}}
       \int_{\R}W_R(y)\dd y\notag\\
 &\le CR^{1-\gamma}.
\end{align}
Evenness of $J$ gives the same bound for $W_R(x+z)$.
Since $X>2R$ and $R^{1-\gamma}\le C\omega_\gamma(R)$,
these estimates prove \eqref{eq:smooth-upper} for large $R$.

To cover $2\le R\le R_0$ for a fixed sufficiently large $R_0$,
\eqref{eq:smooth-comparison} gives
$ch_\gamma\le W_R\le C_{R_0}h_\gamma$.
By \eqref{eq:J-h} and \eqref{eq:h-convolution},
$|J*W_R-W_R|\le C_{R_0}h_\gamma\le C'_{R_0}W_R$.
The positive continuous function $\omega_\gamma$ has a positive
minimum on $[2,R_0]$, so increasing $C$ proves
\eqref{eq:smooth-upper} on the full stated range.
\end{proof}

For the lower bound, evaluate at $x=0$. Since $W_R(0)=1$, $W_R\le1$, and $J$ has unit mass, $1-(J*W_R)(0)\ge0$.

\begin{lemma}\label{lem:smooth-lower}
There exist $c,C>0$ such that
\begin{equation}\label{eq:smooth-lower}
 c\omega_\gamma(R)\le1-(J*W_R)(0)\le C\omega_\gamma(R)
 \text{ for } R\ge2.
\end{equation}
\end{lemma}

\begin{proof}
Define $G_\gamma(s)=1-(1+s^2)^{-\gamma/2}$ for $s\ge0$. Since
$G_\gamma(s)/s^2\to\gamma/2$ as $s\to0$ and $G_\gamma$ increases to one, there are positive constants $c,C$ such that
\begin{equation}\label{eq:G-comparison}
 c\min\{s^2,1\}\le G_\gamma(s)\le C\min\{s^2,1\}
 \text{ for } s\ge0.
\end{equation}
On $(0,1]$ this follows by extending $G_\gamma(s)/s^2$ continuously to zero and taking its positive minimum and finite maximum. On $[1,\infty)$ it follows from $G_\gamma(1)>0$ and $G_\gamma\le1$.

The properties of $G_\gamma$ used here follow directly from
\[
 G_\gamma(s)=\gamma\int_0^s v(1+v^2)^{-\gamma/2-1}\dd v,
 \qquad G_\gamma'(s)=\gamma s(1+s^2)^{-\gamma/2-1}.
\]
After $v=sw$, division of the first identity by $s^2$ gives a limit $\gamma\int_0^1 w\dd w=\gamma/2$ as $s\downarrow0$, by dominated convergence with majorant $\gamma w$ on $[0,1]$.

Because $W_R(0)=1$ and $J$ has unit mass,
\begin{equation}\label{eq:center-defect}
 1-(J*W_R)(0)=\int_{\R}J(z)G_\gamma(|z|/R)\dd z.
\end{equation}
By \eqref{eq:G-comparison}, $1-(J*W_R)(0)$ is comparable to
$R^{-2}V_2(R)+T_J(R)$. By Lemma~\ref{lem:kernel-moments}, $1-(J*W_R)(0)\sim\omega_\gamma(R)$. The estimate extends to any bounded range of $R\ge2$: $R\mapsto1-(J*W_R)(0)$ is continuous and strictly positive there, since $J$ has positive mass away from zero.

Indeed, for $R_n\to R\ge2$, $G_\gamma(|z|/R_n)\to G_\gamma(|z|/R)$ pointwise, and $J(z)$ dominates the integrands in \eqref{eq:center-defect}; dominated convergence proves continuity. The bounds in \eqref{eq:tail-explicit} give positive mass to $|z|>\max\{R_J,1\}$, where $G_\gamma(|z|/R)>0$, so \eqref{eq:center-defect} is positive for every $R$. Its ratio to the continuous positive $\omega_\gamma(R)$ has a positive minimum and a finite maximum on each compact interval of $R\ge2$.

\end{proof}

\begin{proof}[Proof of Proposition~\ref{thm:smooth-error}]
By Lemma~\ref{lem:smooth-upper}, \eqref{eq:smooth-upper} bounds the supremum in \eqref{eq:smooth-error} by $C\omega_\gamma(R)$, whereas
\[
 \sup_{x\in\R}\frac{|(J*W_R)(x)-W_R(x)|}{W_R(x)}
 \ge\frac{|(J*W_R)(0)-W_R(0)|}{W_R(0)}
 =1-(J*W_R)(0)\ge c\omega_\gamma(R)
\]
by $W_R(0)=1$ and \eqref{eq:smooth-lower}. Taking $R_*=2$ after enlarging the constants proves both claims of Proposition~\ref{thm:smooth-error}.

\end{proof}

The supremum in \eqref{eq:truncated-error} is comparable to $V_1(R)/R+T_J(R)$, whereas the supremum in \eqref{eq:smooth-error} is comparable to $V_2(R)/R^2+T_J(R)$ as $R\to\infty$. The second comparison uses the cancellation in \eqref{eq:central-difference}.

\subsection{Limits of the convolution error}

We next compute the limits of $[1-(J*W_R)(0)]/\omega_\gamma(R)$, with $\omega_\gamma$ defined in \eqref{eq:omega-def}.

\begin{proposition}\label{prop:center-constants}
Assume \eqref{eq:J} and \eqref{eq:tail}. If $1<\gamma<3$ and \eqref{eq:exact-tail} holds, then
\begin{equation}\label{eq:center-constant-heavy}
 \lim_{R\to\infty}R^{\gamma-1}[1-(J*W_R)(0)]
 =2\lambda\int_0^\infty
       [1-(1+s^2)^{-\gamma/2}]s^{-\gamma}\dd s.
\end{equation}
If $\gamma=3$ and \eqref{eq:exact-tail} holds, then
\begin{equation}\label{eq:center-constant-critical}
 \lim_{R\to\infty}\frac{R^2}{\log R}[1-(J*W_R)(0)]=3\lambda.
\end{equation}
If $\gamma>3$, then
\begin{equation}\label{eq:center-constant-finite}
 \lim_{R\to\infty}R^2[1-(J*W_R)(0)]
 =\frac\gamma2\int_{\R}z^2J(z)\dd z.
\end{equation}
\end{proposition}

\begin{proof}
First let $1<\gamma<3$ and fix $A\ge\max\{1,R_J\}$.
By \eqref{eq:G-comparison}, the contribution of $|z|\le A$
to $R^{\gamma-1}[1-(J*W_R)(0)]$ is at most
$CR^{\gamma-3}V_2(A)\to0$, since $V_2(A)\le A^2$.
Evenness and $z=Rs$ give the remaining expression as
\[
 2\int_0^\infty
 \ind_{\{s>A/R\}}R^\gamma J(Rs)G_\gamma(s)\dd s.
\]
For each $s>0$, its integrand tends to
$2\lambda s^{-\gamma}G_\gamma(s)$ by \eqref{eq:exact-tail}.
Equations \eqref{eq:tail-explicit} and \eqref{eq:G-comparison}
give the majorant $Cs^{-\gamma}\min\{s^2,1\}$ on $(0,\infty)$,
integrable because $1<\gamma<3$.
Dominated convergence proves \eqref{eq:center-constant-heavy};
the limiting integral is positive because $G_\gamma(s)>0$ for $s>0$.

Let $\gamma=3$. We first claim that
\begin{equation}\label{eq:V2-limit}
 \frac{V_2(R)}{\log R}\to2\lambda\qquad\text{as }R\to\infty.
\end{equation}
For $0<\varepsilon<\lambda$, \eqref{eq:exact-tail} gives,
for some $A$ and all $R\ge A$,
\[
 V_2(A)+2(\lambda-\varepsilon)\log(R/A)
 \le V_2(R)\le
 V_2(A)+2(\lambda+\varepsilon)\log(R/A).
\]
Dividing by $\log R$ and then letting $\varepsilon\downarrow0$
proves \eqref{eq:V2-limit}. It also follows that
$V_2(\delta R)/\log R\to2\lambda$ for every fixed $\delta>0$.
Given $\varepsilon>0$, choose $\delta\in(0,1)$ such that
$|G_3(s)/s^2-3/2|\le\varepsilon$ for $0<s\le\delta$.
The contribution of $|z|\le\delta R$ in \eqref{eq:center-defect},
multiplied by $R^2/\log R$, then lies between
$(3/2-\varepsilon)V_2(\delta R)/\log R$ and
$(3/2+\varepsilon)V_2(\delta R)/\log R$.
The remaining nonnegative contribution is at most
$CR^2T_J(\delta R)/\log R\le C\delta^{-2}/\log R\to0$
by \eqref{eq:tail-mass}.
Letting $R\to\infty$ first and then $\varepsilon\downarrow0$
proves \eqref{eq:center-constant-critical}.

Finally, for $\gamma>3$, \eqref{eq:moment-exact-criterion}
gives $m_2(J)<\infty$. At each $z\in\R$,
$R^2G_\gamma(|z|/R)\to(\gamma/2)z^2$.
By \eqref{eq:G-comparison},
$J(z)R^2G_\gamma(|z|/R)\le Cz^2J(z)$ on all of $\R$,
and this majorant is integrable.
Dominated convergence in \eqref{eq:center-defect} proves
\eqref{eq:center-constant-finite}.
\end{proof}

\subsection{The Fourier transform}

Equations \eqref{eq:truncated-moment} for $V_2(R)$ and \eqref{eq:tail-mass} for $T_J(R)$ yield \eqref{eq:symbol-orders} for the Fourier multiplier of $\Lop$ in \eqref{eq:operator}; see \cite[Chapter~1]{Andreu2010} for its role in linear diffusion. For $\xi\in\R$, with $i^2=-1$, define
\begin{equation}\label{eq:Fourier-symbol}
 \widehat J(\xi)=\int_{\R}\ee^{-i\xi z}J(z)\dd z,
 \qquad D_J(\xi)=1-\widehat J(\xi)
                =\int_{\R}(1-\cos(\xi z))J(z)\dd z.
\end{equation}
Since $J$ is even and $|\sin(\xi z)|J(z)\le J(z)\in L^1(\R)$, the sine integral vanishes, proving the second identity in \eqref{eq:Fourier-symbol}. In particular, $D_J$ is real and nonnegative. For $v\in L^1(\R)$, the convolution theorem \cite[Proposition~1.2(5)]{Andreu2010} gives
\[
 \widehat{\Lop v}(\xi)=-dD_J(\xi)\widehat v(\xi),
 \qquad \widehat v(\xi)=\int_{\R}\ee^{-i\xi x}v(x)\dd x,
\]
because both $J$ and $v$ belong to $L^1(\R)$.

\begin{proposition}\label{prop:symbol}
Under \eqref{eq:J} and \eqref{eq:tail}, for $0<|\xi|\ll1$,
\begin{equation}\label{eq:symbol-orders}
 D_J(\xi)\sim
 \begin{cases}
 |\xi|^{\gamma-1},&1<\gamma<3,\\
 \xi^2\log(1/|\xi|),&\gamma=3,\\
 \xi^2,&\gamma>3.
 \end{cases}
\end{equation}
If $1<\gamma<3$ and \eqref{eq:exact-tail} holds, then
\begin{equation}\label{eq:symbol-constant-heavy}
 \lim_{\xi\to0}\frac{D_J(\xi)}{|\xi|^{\gamma-1}}
 =A_\gamma:=2\lambda\int_0^\infty(1-\cos s)s^{-\gamma}\dd s\in(0,\infty).
\end{equation}
If $\gamma=3$ and \eqref{eq:exact-tail} holds, then
\begin{equation}\label{eq:symbol-constant-critical}
 \lim_{\xi\to0}\frac{D_J(\xi)}{\xi^2\log(1/|\xi|)}=\lambda.
\end{equation}
If $\gamma>3$, then
\begin{equation}\label{eq:symbol-constant-finite}
 \lim_{\xi\to0}\frac{D_J(\xi)}{\xi^2}
 =\frac12 m_2(J).
\end{equation}
\end{proposition}

\begin{proof}
Let $R=|\xi|^{-1}$. The upper bound $1-\cos s\le\min\{s^2/2,2\}$ gives
\[
 D_J(\xi)\le \tfrac12R^{-2}V_2(R)+2T_J(R).
\]
On $|z|\le R$, one has $1-\cos(\xi z)\ge c\xi^2z^2$, because
$(1-\cos s)/s^2$ has a positive continuous extension on $[-1,1]$. Therefore
$D_J(\xi)\ge cR^{-2}V_2(R)$. The estimates \eqref{eq:truncated-moment} for $V_2(R)$ and \eqref{eq:tail-mass} for $T_J(R)$ prove \eqref{eq:symbol-orders}.

For $1<\gamma<3$, assume \eqref{eq:exact-tail} and fix $A\ge\max\{1,R_J\}$. For $\xi\ne0$, the change of variables $s=|\xi|z$ and evenness of $J$ give the normalized contribution from $|z|>A$ as
\[
 2\int_0^\infty\ind_{\{s>A|\xi|\}}|\xi|^{-\gamma}J(s/|\xi|)(1-\cos s)\dd s.
\]
The majorant $C\min\{s^2,1\}s^{-\gamma}$ is integrable on $(0,\infty)$ for $1<\gamma<3$, while \eqref{eq:exact-tail} gives the pointwise limit $2\lambda(1-\cos s)s^{-\gamma}$ for every $s>0$. On $|z|\le A$, $1-\cos(\xi z)\le\xi^2z^2/2$ bounds the normalized contribution by $|\xi|^{3-\gamma}V_2(A)/2\to0$. Dominated convergence therefore proves \eqref{eq:symbol-constant-heavy}. The constant $A_\gamma$ is strictly positive, for example by restricting its defining integral to $[1,2]$.

For $\gamma=3$, assume \eqref{eq:exact-tail}. Fix $0<\varepsilon<1/2$ and choose $\delta\in(0,1)$ such that
\[
 \left|\frac{1-\cos s}{s^2}-\frac12\right|\le\varepsilon
 \quad\text{for }0<|s|\le\delta.
\]
Applying \eqref{eq:V2-limit} with $R=\delta/|\xi|$ gives
\[
 \frac{V_2(\delta/|\xi|)}{\log(1/|\xi|)}\longrightarrow2\lambda.
\]
The contribution of $|z|\le\delta/|\xi|$ to $D_J(\xi)$ in \eqref{eq:Fourier-symbol}, divided by $\xi^2\log(1/|\xi|)$, therefore has lower limit at least $(1-2\varepsilon)\lambda$ and upper limit at most $(1+2\varepsilon)\lambda$. By \eqref{eq:tail-mass}, the nonnegative contribution of $|z|>\delta/|\xi|$ to $D_J(\xi)$, with the same normalization, is at most $C\delta^{-2}/\log(1/|\xi|)\to0$. Letting $\varepsilon\downarrow0$ proves \eqref{eq:symbol-constant-critical} with the order of limits specified.

For $\gamma>3$, \eqref{eq:moment-exact-criterion} gives $m_2(J)<\infty$. By \eqref{eq:J}, $J$ satisfies hypothesis~(H) of \cite[Chapter~1]{Andreu2010} in dimension one. Hence \cite[Lemma~1.8(ii)]{Andreu2010} gives $\widehat J'(0)=0$ and $\widehat J''(0)=-m_2(J)$. Taylor's formula yields $D_J(\xi)=\tfrac12m_2(J)\xi^2+o(\xi^2)$, proving \eqref{eq:symbol-constant-finite}.

\end{proof}

We now rescale $X_t$ so that its characteristic function has a nonconstant limit as $t\to\infty$. The Fourier transform of the probability measure \eqref{eq:mu} is
\begin{equation}\label{eq:mu-Fourier}
 \widehat\mu_t(\xi)=\exp[-dtD_J(\xi)].
\end{equation}
Here $\widehat\mu_t(\xi)=\int_{\R}e^{-i\xi z}\dd\mu_t(z)=\E e^{-i\xi X_t}$. Formula \eqref{eq:mu-Fourier} follows by conditioning on $N_t$ and summing the exponential series, as in the compound Poisson characteristic function in \cite[Chapter~1]{Applebaum2009}. Substitution of the limits in Proposition~\ref{prop:symbol} into \eqref{eq:mu-Fourier} gives the following Fourier limits.

For an analytic derivation, integrate each term of \eqref{eq:mu} against $\ee^{-i\xi z}$. This exchange of the sum and integral is justified by
\[
 \ee^{-dt}+\ee^{-dt}\sum_{n=1}^\infty\frac{(dt)^n}{n!}
 \int_{\R}|\ee^{-i\xi z}|J^{*n}(z)\dd z=1,
\]
where $J^{*0}=\delta_0$ and $J^{*n}(z)\dd z$ has unit mass for $n\ge1$. The convolution identity proved by Fubini gives $\widehat{J^{*n}}=\widehat J^{\,n}$, including $n=0$, and hence
\[
 \widehat\mu_t(\xi)=\ee^{-dt}\sum_{n=0}^\infty
 \frac{(dt\widehat J(\xi))^n}{n!}
 =\exp\{dt[\widehat J(\xi)-1]\}.
\]
The series is absolutely convergent since $|\widehat J(\xi)|\le1$. In probability notation, independence of the jumps $Z_k$ gives $\E[\ee^{-i\xi X_t}\mid N_t=n]=\widehat J(\xi)^n$; summation against the Poisson probabilities gives the same formula. Thus the measure and process interpretations both yield \eqref{eq:mu-Fourier}.

\begin{corollary}\label{cor:linear-scales}
Assume \eqref{eq:J} and \eqref{eq:tail}. If $1<\gamma\le3$, assume in addition \eqref{eq:exact-tail}. Set, for $t>1$,
\begin{equation}\label{eq:linear-scales}
 a_\gamma(t)=
 \begin{cases}
 t^{1/(\gamma-1)},&1<\gamma<3,\\
 (t\log t)^{1/2},&\gamma=3,\\
 t^{1/2},&\gamma>3.
 \end{cases}
\end{equation}
For every fixed $\xi\in\R$,
\begin{equation}\label{eq:linear-transform-limits}
 \lim_{t\to\infty}\widehat\mu_t(\xi/a_\gamma(t))=
 \begin{cases}
 \exp[-dA_\gamma|\xi|^{\gamma-1}],&1<\gamma<3,\\
 \exp[-d\lambda\xi^2/2],&\gamma=3,\\
 \exp[-dm_2(J)\xi^2/2],&\gamma>3.
 \end{cases}
\end{equation}
In every case, $a_\gamma(t)/S_\gamma(t)\to0$ as $t\to\infty$.
\end{corollary}

\begin{proof}
For $\xi=0$, one has $\widehat\mu_t(0)=1$.
Fix $\xi\ne0$. If $1<\gamma<3$, equations
\eqref{eq:symbol-constant-heavy} and \eqref{eq:linear-scales}
give
$tD_J(\xi/a_\gamma(t))\to A_\gamma|\xi|^{\gamma-1}$,
since $t/a_\gamma(t)^{\gamma-1}=1$.
If $\gamma>3$, \eqref{eq:symbol-constant-finite} gives
$tD_J(\xi/a_\gamma(t))\to m_2(J)\xi^2/2$,
since $t/a_\gamma(t)^2=1$.
For $\gamma=3$,
\[
 \frac{t}{a_3(t)^2}\log\frac{a_3(t)}{|\xi|}
 =\frac{\tfrac12\log t+\tfrac12\log\log t-\log|\xi|}{\log t}
 \longrightarrow\frac12.
\]
Together with \eqref{eq:symbol-constant-critical}, this yields
$tD_J(\xi/a_3(t))\to\lambda\xi^2/2$.
Substitution into \eqref{eq:mu-Fourier} and continuity of the
exponential prove \eqref{eq:linear-transform-limits}.
Finally,
$\log[a_\gamma(t)/S_\gamma(t)]=O(\log t)-rt/\gamma\to-\infty$,
so $a_\gamma(t)/S_\gamma(t)\to0$.
\end{proof}

For $1<\gamma<3$, the first limit in \eqref{eq:linear-transform-limits} is the characteristic function of a symmetric stable distribution with index $\gamma-1$. For $\gamma=3$ and $\gamma>3$, \eqref{eq:linear-transform-limits} gives characteristic functions of centered Gaussian distributions; see \cite[Chapter~3]{Durrett2019}.

\subsection{The five ranges of the exponent}

We collect the consequences in the five ranges separated by the thresholds $\gamma=2$ for $m_1(J)<\infty$ and $\gamma=3$ for $m_2(J)<\infty$. All assertions about $u$ use the hypotheses of Theorem~\ref{thm:main}. Finiteness of $M_{p,q}(t)$ is considered for each fixed $t>0$, and its growth concerns $t\to\infty$. The suprema in \eqref{eq:truncated-error} and \eqref{eq:smooth-error} are estimated as $R\to\infty$, and $D_J(\xi)$ is estimated as $\xi\to0$.

\paragraph{The range $1<\gamma<2$.}
Here $m_1(J)=m_2(J)=\infty$. By Corollary~\ref{cor:moments}, $\int_{\R}|x|u(t,x)\dd x=\int_{\R}x^2u(t,x)\dd x=\infty$ for every $t>0$. By \eqref{eq:moment-growth} with $(p,q)=(0,1)$, $\mathcal M(t)\sim\ee^{rt/\gamma}$ as $t\to\infty$. For each $\theta\in(0,1)$, the comparison $|x|\sim\ee^{rt/\gamma}$ holds uniformly for $x\in E_\theta(t)$ as $t\to\infty$, by \eqref{eq:all-points-annulus}. Both suprema in \eqref{eq:truncated-error} and \eqref{eq:smooth-error} have order $R^{1-\gamma}$. The orders coincide because $V_1(R)/R$, $V_2(R)/R^2$, and $T_J(R)$ all have order $R^{1-\gamma}$.

\paragraph{The endpoint $\gamma=2$.}
The truncated moments $V_1(R)$ and $V_2(R)$ of $J$ grow logarithmically and linearly, respectively:
$V_1(R)\sim\log R$ and $V_2(R)\sim R$. Thus the supremum in \eqref{eq:truncated-error} has order $R^{-1}\log R$, whereas both $1-(J*W_R)(0)$ and $\sup_{x\in\R}|(J*W_R)(x)-W_R(x)|/W_R(x)$ have order $R^{-1}$. For each $t>0$, $M_{1,1}(t;L)$ grows like $\log L$ as $L\to\infty$. The level-set estimate is
\begin{equation}\label{eq:case-two}
 \log|x|=\frac r2t+O_\theta(1)
 \text{ for } x\in E_\theta(t)\text{ as } t\to\infty.
\end{equation}
Thus \eqref{eq:case-two} gives a bounded remainder, while $V_1(R)\sim\log R$ produces the factor $\log R$ in \eqref{eq:truncated-error}.

\paragraph{The range $2<\gamma<3$.}
Here $m_1(J)<\infty$ and $m_2(J)=\infty$. By \eqref{eq:moment-growth}, $M_{1,1}(t)\sim\ee^{2rt/\gamma}$ as $t\to\infty$, while \eqref{eq:truncated-moment} gives $V_2(R)\sim R^{3-\gamma}$ as $R\to\infty$. The supremum in \eqref{eq:smooth-error} has order $R^{1-\gamma}$, whereas the supremum in \eqref{eq:truncated-error} has order $R^{-1}$ as $R\to\infty$. By Theorem~\ref{thm:main}\textup{(ii)}, the level-set scale is $\ee^{rt/\gamma}$.

\paragraph{The endpoint $\gamma=3$.}
Here $m_1(J)<\infty$, and $V_2(R)\sim\log R$. For each $t>0$, $M_{2,1}(t)=\infty$ and $M_{2,1}(t;L)$ grows like $\log L$ as $L\to\infty$. The quantity $\sup_{x\in\R}|(J*W_R)(x)-W_R(x)|/W_R(x)$ has order $R^{-2}\log R$, and $D_J(\xi)$ has order $\xi^2\log(1/|\xi|)$. Nevertheless,
\begin{equation}\label{eq:case-three}
 \log|x|=\frac r3t+O_\theta(1)
 \text{ for } x\in E_\theta(t)\text{ as } t\to\infty.
\end{equation}
The logarithmic factors in the supremum in \eqref{eq:smooth-error} and in $D_J(\xi)$ in \eqref{eq:symbol-orders} arise from $V_2(R)\sim\log R$, while \eqref{eq:case-three} has a bounded remainder.

\paragraph{The range $\gamma>3$.}
Here $m_1(J)<\infty$ and $m_2(J)<\infty$. By \eqref{eq:moment-growth},
$\int_{\R}|x|u(t,x)\dd x\sim\ee^{2rt/\gamma}$ and
$\int_{\R}x^2u(t,x)\dd x\sim\ee^{3rt/\gamma}$. By Proposition~\ref{thm:smooth-error}, the relative convolution error $\sup_{x\in\R}|(J*W_R)(x)-W_R(x)|/W_R(x)$ for $W_R$ in \eqref{eq:W-def} has order $R^{-2}$ as $R\to\infty$. The finiteness of $\int_{\R}z^2J(z)\dd z$ gives the quadratic behavior of $D_J(\xi)$ in \eqref{eq:symbol-constant-finite}, and Theorem~\ref{thm:main}\textup{(ii)} gives exponential level-set propagation.

The table records the orders of $V_1(R)$, $V_2(R)$ and the suprema in \eqref{eq:truncated-error} and \eqref{eq:smooth-error} as $R\to\infty$.
\begin{center}
\small
\begin{tabular}{@{}lllll@{}}
\toprule
Range & $V_1(R)$ & $V_2(R)$ & Error for $\Psi_R$ & Error for $W_R$\\
\midrule
$1<\gamma<2$ & $R^{2-\gamma}$ & $R^{3-\gamma}$ & $R^{1-\gamma}$ & $R^{1-\gamma}$\\
$\gamma=2$ & $\log R$ & $R$ & $R^{-1}\log R$ & $R^{-1}$\\
$2<\gamma<3$ & $1$ & $R^{3-\gamma}$ & $R^{-1}$ & $R^{1-\gamma}$\\
$\gamma=3$ & $1$ & $\log R$ & $R^{-1}$ & $R^{-2}\log R$\\
$\gamma>3$ & $1$ & $1$ & $R^{-1}$ & $R^{-2}$\\
\bottomrule
\end{tabular}
\end{center}

The last two columns refer to the suprema in \eqref{eq:truncated-error} and \eqref{eq:smooth-error}, respectively. The values of $V_1(R)$ and $V_2(R)$ follow from \eqref{eq:truncated-moment}, and the claims about $M_{1,1}(t;L)$ and $M_{2,1}(t;L)$ at the endpoints follow from \eqref{eq:truncated-solution-moments}.

In the columns for $V_k(R)$, $k=1,2$, an entry $1$ means that $V_k(R)$ converges to the finite positive integral $\int_{\R}|z|^kJ(z)\dd z$ as $R\to\infty$.

Direct integration along the nonlinear exponential scale gives, for sufficiently large $t$,
\begin{equation}\label{eq:integrated-smooth-error}
 \int_t^\infty\omega_\gamma(S_\gamma(s))\dd s=
 \begin{cases}
 \dfrac{\gamma}{r(\gamma-1)}\ee^{-r(\gamma-1)t/\gamma},&1<\gamma<3,\\[2mm]
 \left(\dfrac t2+\dfrac{3}{4r}\right)\ee^{-2rt/3},&\gamma=3,\\[2mm]
 \dfrac{\gamma}{2r}\ee^{-2rt/\gamma},&\gamma>3.
 \end{cases}
\end{equation}
Thus $\int_t^\infty\omega_\gamma(S_\gamma(s))\dd s<\infty$ for every $\gamma>1$ and all sufficiently large $t$.

To verify \eqref{eq:integrated-smooth-error}, take $t$ so large that $S_\gamma(t)\ge2$, and substitute $S_\gamma(s)=\ee^{rs/\gamma}$ in \eqref{eq:omega-def}. For $1<\gamma<3$ and $\gamma>3$ this gives $\int_t^\infty\ee^{-as}\dd s=\ee^{-at}/a$, with $a=r(\gamma-1)/\gamma$ and $a=2r/\gamma$, respectively. For $\gamma=3$ the integrand is $(rs/3)\ee^{-2rs/3}$; integration by parts yields
\[
 \frac r3\int_t^\infty s\ee^{-2rs/3}\dd s
 =\left(\frac t2+\frac{3}{4r}\right)\ee^{-2rt/3},
\]
since $s\ee^{-2rs/3}\to0$ as $s\to\infty$.

The proof of Theorem~\ref{thm:main} uses $\Phi_R$ and $\Psi_R$ with \eqref{eq:relative} and $\int_4^\infty e_\gamma(R)\dd R/R<\infty$, where $e_\gamma$ is defined in \eqref{eq:basic-error}.

\section{Spatial tails and examples}\label{sec:tail}\label{sec:examples}

In this section, \eqref{eq:exact-tail} is assumed in Lemmas~\ref{lem:tail-space}--\ref{lem:uniform-fixed-time} and Proposition~\ref{prop:linear-tail}. Proposition~\ref{prop:nonlinear-mass-identity} uses the hypotheses of Theorem~\ref{thm:main}. We prove Theorem~\ref{thm:tail-coefficient} by taking the spatial limit at fixed time and then studying $\mathcal A(t)$ as $t\to\infty$. Fixed-time tail asymptotics for fractional diffusion appear in \cite[Lemma~5.1]{RoquejoffreTarfulea2017}.

Since $|x|^\gamma h_\gamma(x)\to1$, condition \eqref{eq:exact-tail} is equivalent to $J(x)/h_\gamma(x)\to\lambda$ at both ends. We compute the spatial tails of $u(t,\cdot)$ by taking limits of $u(t,x)/h_\gamma(x)$ at the two ends of $\R$.

\subsection{Tail limits for a convolution}\label{subsec:tail-convolution}

We use \cite[Proposition~2.2]{FinkelshteinTkachovKesten2018} to prove \eqref{eq:tail-convolution} for signed functions. For background on subexponential densities, see \cite[Chapter~4]{FossKorshunovZachary2013}. To keep track of the two spatial limits, define the following subspace of $\Yg$ from \eqref{eq:Y}:
\begin{equation}\label{eq:X}
 \Xg=\left\{v\in\Yg:
 \ell_+(v):=\lim_{x\to+\infty}\frac{v(x)}{h_\gamma(x)}
 \text{ and }
 \ell_-(v):=\lim_{x\to-\infty}\frac{v(x)}{h_\gamma(x)}
 \text{ exist in }\R\right\}.
\end{equation}
We equip $\Xg$ with the norm $\|\cdot\|_\gamma$ and allow its elements to change sign.

\begin{lemma}\label{lem:tail-space}
The space $\Xg$ is a closed linear subspace of $\Yg$, hence a Banach space. The maps $\ell_+$ and $\ell_-$ are continuous linear functionals satisfying
$|\ell_\pm(v)|\le\|v\|_\gamma$. For each fixed $A>0$,
\begin{equation}\label{eq:h-shift-limit}
 \lim_{x\to\pm\infty}\sup_{|z|\le A}
 \left|\frac{h_\gamma(x-z)}{h_\gamma(x)}-1\right|=0.
\end{equation}
\end{lemma}

\begin{proof}
Suppose $v_n\in\Xg$ and $v_n\to v$ in $\Yg$. The quotients $v_n/h_\gamma$ converge uniformly on $\R$ to $v/h_\gamma$. Their limiting values at either end form a Cauchy sequence because
$|\ell_\pm(v_n)-\ell_\pm(v_m)|\le\|v_n-v_m\|_\gamma$. Denote the limits of these sequences by $a_\pm$. Given $\varepsilon>0$, choose $n$ with
$\|v-v_n\|_\gamma<\varepsilon$ and $|\ell_\pm(v_n)-a_\pm|<\varepsilon$. For $x$ sufficiently far toward the chosen end,
$|v_n(x)/h_\gamma(x)-\ell_\pm(v_n)|<\varepsilon$. Adding the three bounds gives
$|v(x)/h_\gamma(x)-a_\pm|<3\varepsilon$. Thus $v\in\Xg$. 

For $v,w\in\Xg$ and $a,b\in\R$, taking the limits in \eqref{eq:X} gives
\[
 \ell_\pm(av+bw)=a\ell_\pm(v)+b\ell_\pm(w),\qquad
 |\ell_\pm(v)|\le\sup_{x\in\R}\frac{|v(x)|}{h_\gamma(x)}=\|v\|_\gamma.
\]
Thus $\Xg$ is linear and $|\ell_\pm(v)-\ell_\pm(w)|\le\|v-w\|_\gamma$, which proves continuity of $\ell_\pm$. Completeness of the closed subspace $\Xg$ follows from completeness of $\Yg$ in Lemma~\ref{lem:weight}.

To prove \eqref{eq:h-shift-limit}, take $|x|>2A$ and $|z|\le A$. The inequality $||x-z|-|x||\le A$ gives
\[
 \left|\frac{1+|x-z|}{1+|x|}-1\right|\le\frac{A}{1+|x|},
 \qquad
 \left|\frac{h_\gamma(x-z)}{h_\gamma(x)}-1\right|
 \le\gamma2^{\gamma+1}\frac{A}{1+|x|}.
\]
The second inequality is the mean value theorem for $s\mapsto s^{-\gamma}$ on $[1/2,3/2]$, where its derivative has absolute value at most $\gamma2^{\gamma+1}$. This proves \eqref{eq:h-shift-limit} at both ends.
\end{proof}

\begin{lemma}\label{lem:tail-convolution}
Convolution by $J$ is a bounded operator on $\Xg$. For each $v\in\Xg$,
\begin{equation}\label{eq:tail-convolution}
 \ell_\pm(J*v)=\ell_\pm(v)+\lambda\int_{\R}v(y)\dd y.
\end{equation}
\end{lemma}

\begin{proof}
By \eqref{eq:J-weight-bound}, convolution by $J$ is bounded on $\Yg$. We derive \eqref{eq:tail-convolution} from \cite[Proposition~2.2]{FinkelshteinTkachovKesten2018}, checking its hypotheses for $h_\gamma$.

On $(0,\infty)$, $h_\gamma(x)=(1+x)^{-\gamma}$ is decreasing, integrable and log-convex. Choose $a\in(1/\gamma,1)$. As $x\to\infty$,
\[
 \sup_{|z|\le x^a}\left|\frac{h_\gamma(x+z)}{h_\gamma(x)}-1\right|\longrightarrow0,
 \qquad xh_\gamma(x^a)\longrightarrow0.
\]
The first limit follows from $x^a/(1+x)\to0$, and the second from $1-a\gamma<0$. Thus \cite[Lemma~2.4]{FinkelshteinTkachovKesten2018} shows that $h_\gamma$ is weakly subexponential at $+\infty$, as required by \cite[Proposition~2.2]{FinkelshteinTkachovKesten2018}. Its additional bound $h_\gamma(x+z)\le h_\gamma(x)$ for $x,z>0$ follows from monotonicity.

For $v\in\Xg$, set $v^+(x)=\max\{v(x),0\}$ and $v^-(x)=\max\{-v(x),0\}$. Then $v^\pm\in L^1(\R)\cap\Xg$, since $0\le v^\pm\le\|v\|_\gamma h_\gamma$ and $\gamma>1$. Moreover, $\ell_+(v^\pm)=\max\{\pm\ell_+(v),0\}$ and $J(x)/h_\gamma(x)\to\lambda$. Applying \cite[Proposition~2.2]{FinkelshteinTkachovKesten2018} with reference function $h_\gamma$ to $J$ and each of $v^+$ and $v^-$ gives, using $\int_\R J=1$,
\[
 \lim_{x\to+\infty}\frac{(J*v^\pm)(x)}{h_\gamma(x)}
 =\ell_+(v^\pm)+\lambda\int_\R v^\pm(y)\dd y.
\]
Subtracting proves \eqref{eq:tail-convolution} at $+\infty$. For $\check v(x)=v(-x)$, evenness of $J$ gives $(J*\check v)(x)=(J*v)(-x)$, $\ell_+(\check v)=\ell_-(v)$ and $\int_\R\check v=\int_\R v$. Applying the positive-tail formula in \eqref{eq:tail-convolution} to $\check v$ proves the limit at $-\infty$. Hence $J*v\in\Xg$, and the norm bound in \eqref{eq:J-weight-bound} completes the proof.
\end{proof}

Since $\Lop v=d(J*v-v)$, \eqref{eq:tail-convolution} yields $\ell_\pm(\Lop v)=d\lambda\int_{\R}v(y)\dd y$ for $v\in\Xg$.

\subsection{Tail limits for the solution}

\begin{lemma}\label{lem:tail-evolution}
The solution of \eqref{eq:cauchy} belongs to $C^1([0,\infty);\Xg)$. The coefficients $\ell_\pm(u(t,\cdot))$ from \eqref{eq:X} satisfy
\begin{equation}\label{eq:tail-ODE}
 \frac{\dd}{\dd t}\ell_\pm(u(t,\cdot))
 =r\ell_\pm(u(t,\cdot))+d\lambda\mathcal M(t),
 \qquad \ell_\pm(u(0,\cdot))=0.
\end{equation}
Consequently,
\begin{equation}\label{eq:tail-ODE-solution}
 \ell_\pm(u(t,\cdot))
 =d\lambda\ee^{rt}\int_0^t\ee^{-rs}\mathcal M(s)\dd s.
\end{equation}
\end{lemma}

\begin{proof}
We use the Lipschitz extension of $f$ from Lemma~\ref{lem:wellposed}. If $L_f$ is its global Lipschitz constant, then $\|f(v)-f(w)\|_\gamma\le L_f\|v-w\|_\gamma$ for $v,w\in\Yg$. To see that composition preserves $\Xg$, let $v\in\Xg$. Since $v(x)\to0$ and $f$ is differentiable at zero,
\[
 f(v(x))=rv(x)+o(|v(x)|)\qquad\text{as }x\to\pm\infty.
\]

For a direct $\varepsilon$ estimate, differentiability of the extended $f$ at zero gives $\eta>0$ such that $|f(s)-rs|\le\varepsilon|s|$ for $|s|\le\eta$. Since $v\in C_0(\R)$, $|v(x)|\le\eta$ for all sufficiently large $|x|$. Consequently
\[
 \frac{|f(v(x))-rv(x)|}{h_\gamma(x)}\le\varepsilon\|v\|_\gamma
\]
there. Hence $[f(v(x))-rv(x)]/h_\gamma(x)\to0$ as $x\to\pm\infty$. Thus
\begin{equation}\label{eq:reaction-tail-limit}
 \ell_\pm(f(v))=r\ell_\pm(v).
\end{equation}
Formula \eqref{eq:reaction-tail-limit} also holds for signed $v$, since the extension satisfies $f(s)=rs$ for $s<0$.

By Lemma~\ref{lem:tail-convolution} and \eqref{eq:reaction-tail-limit}, the vector field $v\mapsto\Lop v+f(v)$ preserves $\Xg$. Moreover, \eqref{eq:J-weight-bound} and the bound with $L_f$ give
\[
 \|\Lop v+f(v)-\Lop w-f(w)\|_\gamma
 \le[d(C_JC_h+1)+L_f]\|v-w\|_\gamma.
\]
Since $u_0\in\Xg$ and $\ell_+(u_0)=\ell_-(u_0)=0$ by \eqref{eq:data}, the Cauchy--Lipschitz theorem gives a global $C^1$ solution in $\Xg$. The continuous embedding into $C_0(\R)$ and uniqueness in Lemma~\ref{lem:wellposed} identify it with $u$.

Since $u\in C^1([0,\infty);\Xg)$ and $\ell_\pm$ is bounded and linear by Lemma~\ref{lem:tail-space}, one has $[\ell_\pm(u(t,\cdot))]^{\prime}=\ell_\pm(u_t(t,\cdot))$. Apply $\ell_\pm$ to $u_t=\Lop u+f(u)$ in $\Xg$. Formula \eqref{eq:tail-convolution} gives
\[
 \ell_\pm(\Lop u)
 =d\left[\ell_\pm(u)+\lambda\mathcal M-\ell_\pm(u)\right]
 =d\lambda\mathcal M.
\]
The identity $\ell_\pm(\Lop u)=d\lambda\mathcal M$ and \eqref{eq:reaction-tail-limit} prove \eqref{eq:tail-ODE}. 

For example, the error in applying $\ell_\pm$ to the difference quotient is at most
\[
 \left\|\frac{u(t+h,\cdot)-u(t,\cdot)}h-u_t(t,\cdot)\right\|_\gamma\longrightarrow0.
\]
Continuity of $\mathcal M$ makes the right-hand side of \eqref{eq:tail-ODE} continuous. Thus
\[
 [\ee^{-rt}\ell_\pm(u(t,\cdot))]'=d\lambda\ee^{-rt}\mathcal M(t).
\]
Integration from zero to $t$, using $\ell_\pm(u_0)=0$, gives \eqref{eq:tail-ODE-solution}.
\end{proof}

\begin{lemma}\label{lem:uniform-fixed-time}
For every $T>0$, the spatial limits
\begin{equation}\label{eq:uniform-fixed-time}
 \frac{u(t,x)}{\ee^{rt}J(x)}
 \to d\int_0^t\ee^{-rs}\mathcal M(s)\dd s
 \qquad\text{as }x\to\pm\infty
\end{equation}
are uniform for $0\le t\le T$.
\end{lemma}

\begin{proof}
By Lemma~\ref{lem:tail-evolution}, the path $t\mapsto u(t,\cdot)$ is continuous from $[0,T]$ into $\Xg$. Hence $\{u(t,\cdot):0\le t\le T\}$ is compact in $\|\cdot\|_\gamma$. We claim that
$u(t,x)/h_\gamma(x)\to\ell_\pm(u(t,\cdot))$ as $x\to\pm\infty$, uniformly for $0\le t\le T$.

Fix $\varepsilon>0$. Choose $v_1,\ldots,v_N$ in the image such that each $u(t,\cdot)$, $0\le t\le T$, is within $\varepsilon$ of at least one $v_j$ in $\|\cdot\|_\gamma$. At the selected end, choose a common $A$ so large that
$|v_j(x)/h_\gamma(x)-\ell_\pm(v_j)|<\varepsilon$ beyond $A$ for every $j$. For each $t$, select $j$ with $\|u(t,\cdot)-v_j\|_\gamma<\varepsilon$. The functional bound in Lemma~\ref{lem:tail-space} gives
\[
 \left|\frac{u(t,x)}{h_\gamma(x)}-\ell_\pm(u(t,\cdot))\right|
 <3\varepsilon
\]
beyond the same $A$. Hence $u(t,x)/h_\gamma(x)\to\ell_\pm(u(t,\cdot))$ as $x\to\pm\infty$, uniformly for $0\le t\le T$.

To check division uniformly, set
\[
 K_T=\sup_{0\le t\le T}\|u(t,\cdot)\|_\gamma<\infty.
\]
For large $|x|$, $J(x)/h_\gamma(x)\ge\lambda/2$. Since $|\ell_\pm(u(t,\cdot))|\le K_T$ and $\ee^{-rt}\le1$,
\[
\begin{split}
 &\left|\frac{u(t,x)}{\ee^{rt}J(x)}-
 \frac{\ee^{-rt}}\lambda\ell_\pm(u(t,\cdot))\right|\\
 &\quad\le\frac2\lambda
 \left|\frac{u(t,x)}{h_\gamma(x)}-\ell_\pm(u(t,\cdot))\right|
 +\frac{2K_T}{\lambda^2}\left|\frac{J(x)}{h_\gamma(x)}-\lambda\right|.
\end{split}
\]
Both terms tend to zero uniformly for $0\le t\le T$. By \eqref{eq:tail-ODE-solution}, $\ee^{-rt}\ell_\pm(u(t,\cdot))/\lambda=\mathcal A(t)$, which proves \eqref{eq:uniform-fixed-time}.
\end{proof}

\begin{proof}[Proof of Theorem~\ref{thm:tail-coefficient}]
Lemmas~\ref{lem:tail-evolution} and~\ref{lem:uniform-fixed-time} give \eqref{eq:tail-formula} and \eqref{eq:tail-uniform-main}. By \eqref{eq:mass-order},
\begin{equation}\label{eq:mass-weighted-sharp}
 c\ee^{-\beta_\gamma t}
 \le\ee^{-rt}\mathcal M(t)
 \le C\ee^{-\beta_\gamma t}
 \text{ for } t\ge1,
 \qquad\beta_\gamma=r\left(1-\frac1\gamma\right)>0.
\end{equation}
Continuity of $\mathcal M$ controls its integral on $[0,1]$, while the upper bound in \eqref{eq:mass-weighted-sharp} is integrable on $[1,\infty)$. Hence $\mathcal A_\infty$ in \eqref{eq:A-infty} is finite, and positivity of $\mathcal M$ gives $\mathcal A_\infty>0$. Since
\[
 \mathcal A_\infty-\mathcal A(t)
 =d\int_t^\infty\ee^{-rs}\mathcal M(s)\dd s,
 \qquad
 \int_t^\infty\ee^{-\beta_\gamma s}\dd s
 =\frac{\ee^{-\beta_\gamma t}}{\beta_\gamma},
\]
integration of \eqref{eq:mass-weighted-sharp} gives
\begin{equation}\label{eq:A-two-sided}
 \frac{dc}{\beta_\gamma}\ee^{-\beta_\gamma t}
 \le\mathcal A_\infty-\mathcal A(t)
 \le\frac{dC}{\beta_\gamma}\ee^{-\beta_\gamma t}
 \text{ for } t\ge1.
\end{equation}
This proves \eqref{eq:A-infty}. Moreover, $\mathcal A'(t)=d\ee^{-rt}\mathcal M(t)>0$ by \eqref{eq:A-def}.
\end{proof}

\subsection{The linear equation and a mass identity}\label{subsec:linear-tail}

Proposition~\ref{prop:linear-tail} computes the limit in \eqref{eq:linear-tail} for the linear equation \eqref{eq:linearized}. We then compare $dtM_0$ with $\mathcal A(t)$ in \eqref{eq:A-def}.

\begin{proposition}\label{prop:linear-tail}
Let $v$ solve
\begin{equation}\label{eq:linearized}
 v_t=d(J*v-v)+rv,\qquad v(0,x)=u_0(x),
 \qquad t>0,\ x\in\R.
\end{equation}
Assume \eqref{eq:J}, \eqref{eq:tail}, \eqref{eq:data}, and \eqref{eq:exact-tail}. Then, for every fixed $t\ge0$,
\begin{equation}\label{eq:linear-tail}
 \lim_{x\to\pm\infty}\frac{v(t,x)}{\ee^{rt}J(x)}
 =dt\int_{\R}u_0(y)\dd y.
\end{equation}
For each integer $n\ge1$, the convolution density $J^{*n}$ has mass one, belongs to $\Xg$, and satisfies
\begin{equation}\label{eq:convolution-power-tail}
 \ell_+(J^{*n})=\ell_-(J^{*n})=n\lambda.
\end{equation}
\end{proposition}

\begin{proof}
By \eqref{eq:J-weight-bound}, the series \eqref{eq:P} converges in operator norm on $\Xg$, uniformly on bounded time intervals: its $n$th term has norm at most $\ee^{-dt}(dtC_JC_h)^n/n!$. Thus $v(t,\cdot)=\ee^{rt}P_tu_0$ belongs to $C^1([0,\infty);\Xg)$. Mass preservation in Lemma~\ref{lem:semigroup} gives $\int_\R v(t,x)\dd x=\ee^{rt}M_0$. Applying the bounded functional $\ell_\pm$ to the equation and using \eqref{eq:tail-convolution}, we obtain
\[
 \frac{\dd}{\dd t}\bigl[\ee^{-rt}\ell_\pm(v(t,\cdot))\bigr]
 =d\lambda M_0,
 \qquad \ell_\pm(v(0,\cdot))=0.
\]
Hence $\ell_\pm(v(t,\cdot))=d\lambda t\ee^{rt}M_0$. Since $J/h_\gamma\to\lambda$ at both ends, this proves \eqref{eq:linear-tail}.

For \eqref{eq:convolution-power-tail}, \eqref{eq:J-h} and \eqref{eq:exact-tail} give $J\in\Xg$ and $\ell_\pm(J)=\lambda$. Lemma~\ref{lem:semigroup} gives $\int_\R J^{*n}=1$ for every $n\ge1$. Lemma~\ref{lem:tail-convolution} therefore yields $J^{*(n+1)}\in\Xg$ and
\[
 \ell_\pm(J^{*(n+1)})=\ell_\pm(J^{*n})+\lambda.
\]
Induction proves \eqref{eq:convolution-power-tail}.
\end{proof}

For each fixed $n$, \eqref{eq:convolution-power-tail} agrees with \cite[Theorem~2.1]{FinkelshteinTkachovKesten2018}. By \eqref{eq:A-infty}, the nonlinear coefficient $\mathcal A(t)$ converges to $\mathcal A_\infty$, while the linear coefficient in \eqref{eq:linear-tail} is $dtM_0$.

\paragraph{The nonlinear remainder.}

To compute $\ee^{-rt}\mathcal M(t)$ and $\int_0^\infty\ee^{-rt}\mathcal M(t)\dd t$, define
\begin{equation}\label{eq:nonlinear-remainder}
 G(s)=rs-f(s),\qquad
 \mathcal G(t)=\int_{\R}G(u(t,x))\dd x.
\end{equation}
By \eqref{eq:F}, $G$ is nonnegative on $[0,1]$, and hence $\mathcal G(t)\ge0$. The integral is finite and continuous in $t$, since $0\le G(s)\le rs$ for $0\le s\le1$, $G$ is Lipschitz on $[0,1]$, and Lemma~\ref{lem:weight} gives continuity of $u$ in $\Yg\hookrightarrow L^1(\R)$.

\begin{proposition}\label{prop:nonlinear-mass-identity}
Under the hypotheses of Theorem~\ref{thm:main},
\begin{equation}\label{eq:remainder-total}
 \int_0^\infty\ee^{-rt}\mathcal G(t)\dd t=M_0,
 \qquad M_0=\int_{\R}u_0(x)\dd x,
\end{equation}
and, for every $t\ge0$,
\begin{equation}\label{eq:weighted-mass-remainder}
 \ee^{-rt}\mathcal M(t)
 =\int_t^\infty\ee^{-rs}\mathcal G(s)\dd s.
\end{equation}
Moreover,
\begin{equation}\label{eq:A-remainder}
 d\int_0^\infty\ee^{-rt}\mathcal M(t)\dd t
 =d\int_0^\infty t\ee^{-rt}\mathcal G(t)\dd t\in(0,\infty).
\end{equation}
If \eqref{eq:exact-tail} also holds, both integrals in \eqref{eq:A-remainder} equal $\mathcal A_\infty$ defined in \eqref{eq:A-infty}.
\end{proposition}

\begin{proof}
By \eqref{eq:mass-equation} and \eqref{eq:nonlinear-remainder},
\begin{equation}\label{eq:weighted-mass-derivative}
 \frac{\dd}{\dd t}[\ee^{-rt}\mathcal M(t)]
 =-\ee^{-rt}\mathcal G(t).
\end{equation}
Integrating this continuous identity over $[t,T]$ gives
\[
 \int_t^T\ee^{-rs}\mathcal G(s)\dd s
 =\ee^{-rt}\mathcal M(t)-\ee^{-rT}\mathcal M(T).
\]
By \eqref{eq:mass-order}, the final term tends to zero.
Since $\mathcal G\ge0$, monotone convergence as $T\to\infty$
proves \eqref{eq:weighted-mass-remainder};
taking $t=0$ also gives \eqref{eq:remainder-total}.

The integrand $\ee^{-rs}\mathcal G(s)$ is nonnegative on
$\{(t,s):0\le t\le s<\infty\}$. Thus Tonelli's theorem gives
\[
 \int_0^\infty\ee^{-rt}\mathcal M(t)\dd t
 =\int_0^\infty\int_t^\infty\ee^{-rs}\mathcal G(s)\dd s\dd t
 =\int_0^\infty s\ee^{-rs}\mathcal G(s)\dd s,
\]
where the last equality uses $\int_0^s\dd t=s$.
The first integral is finite by continuity on $[0,1]$ and
\[
 \ee^{-rt}\mathcal M(t)
 \le C\ee^{-r(1-1/\gamma)t}\text{ for }t\ge1,
\]
which follows from \eqref{eq:mass-order}.
Moreover, $\int_0^\infty\ee^{-rt}\mathcal M(t)\dd t>0$
because $\mathcal M(t)>0$.
Multiplication by $d$ proves \eqref{eq:A-remainder}.
Under \eqref{eq:exact-tail}, the common value is
$\mathcal A_\infty$ by \eqref{eq:A-infty}.
\end{proof}

\subsection{Examples}

For each $\gamma>1$, the function
\begin{equation}\label{eq:example-kernel}
 J_\gamma(x)=\frac{\gamma-1}{2}(1+|x|)^{-\gamma}
\end{equation}
is a continuous, bounded, even probability density with $J_\gamma(0)>0$. The normalization $\int_{\R}J_\gamma(x)\dd x=1$ follows from
$2\int_0^\infty(1+x)^{-\gamma}\dd x=2/(\gamma-1)$. It satisfies \eqref{eq:tail} and also
\[
 \lim_{|x|\to\infty}|x|^\gamma J_\gamma(x)=\frac{\gamma-1}{2}.
\]
For this example, direct integration gives
\begin{equation}\label{eq:example-moments}
 \int_{\R}|x|J_\gamma(x)\dd x=\frac1{\gamma-2}
 \text{ for } \gamma>2,\qquad
 \int_{\R}x^2J_\gamma(x)\dd x
 =\frac2{(\gamma-2)(\gamma-3)}\text{ for } \gamma>3.
\end{equation}
By evenness of $J_\gamma$, each integral in \eqref{eq:example-moments} is twice its integral over $[0,\infty)$. Substituting $y=1+x$ and expanding $x=y-1$ and $x^2=(y-1)^2$ then reduces \eqref{eq:example-moments} to power integrals over $[1,\infty)$. Outside the displayed ranges, Lemma~\ref{lem:kernel-moments} gives divergence.

The computations of $m_1(J_\gamma)$ and $m_2(J_\gamma)$ are
\[
\begin{aligned}
 m_1(J_\gamma)
 &=(\gamma-1)\left(\frac1{\gamma-2}-\frac1{\gamma-1}\right)
 =\frac1{\gamma-2} &&\text{ for }\gamma>2,\\
 m_2(J_\gamma)
 &=(\gamma-1)\left(\frac1{\gamma-3}-\frac2{\gamma-2}+\frac1{\gamma-1}\right)
 =\frac2{(\gamma-2)(\gamma-3)} &&\text{ for }\gamma>3.
\end{aligned}
\]
Each term follows from $\int_1^\infty y^{-a}\dd y=1/(a-1)$ for $a>1$; this is the range in which the corresponding antiderivative has a vanishing limit at infinity.

For example, the choices $\gamma=3/2,2,5/2,3,4$ give respective Cauchy level-set scales
\begin{equation}\label{eq:five-examples}
 \ee^{2rt/3},\qquad \ee^{rt/2},\qquad
 \ee^{2rt/5},\qquad \ee^{rt/3},\qquad \ee^{rt/4}.
\end{equation}
By Theorem~\ref{thm:main}\textup{(ii)}, for each listed value of $\gamma$ and each fixed $\theta\in(0,1)$, $|x|$ is comparable to the associated expression in \eqref{eq:five-examples}, uniformly for $x\in E_\theta(t)$ as $t\to\infty$.

Smooth kernels are equally admissible. For example,
\begin{equation}\label{eq:example-smooth}
 J(x)=c_\gamma^{-1}(1+x^2)^{-\gamma/2},
 \qquad c_\gamma=\int_{\R}(1+x^2)^{-\gamma/2}\dd x,
\end{equation}
satisfies \eqref{eq:J} and \eqref{eq:exact-tail} with $\lambda=c_\gamma^{-1}$. By Lemma~\ref{lem:kernel-moments}, $\int_{\R}|x|J(x)\dd x<\infty$ exactly when $\gamma>2$, and $\int_{\R}x^2J(x)\dd x<\infty$ exactly when $\gamma>3$.

An oscillatory example satisfying \eqref{eq:tail-explicit} is
\begin{equation}\label{eq:oscillatory-kernel}
\begin{gathered}
 J_{\rm osc}(x)=Z^{-1}(1+|x|)^{-\gamma}[2+\sin(\log(1+|x|))],\\
 Z=\int_{\R}(1+|x|)^{-\gamma}[2+\sin(\log(1+|x|))]\dd x.
\end{gathered}
\end{equation}
The number $Z$ is finite and positive, since the bracket is between one and three. Thus $J_{\rm osc}$ satisfies \eqref{eq:J} and \eqref{eq:tail-explicit}. For integers $n\ge1$, set
$x_n=\exp(2\pi n+\pi/2)-1$ and
$y_n=\exp(2\pi n+3\pi/2)-1$. As $n\to\infty$,
\[
 x_n^\gamma J_{\rm osc}(x_n)\to3/Z,
 \qquad y_n^\gamma J_{\rm osc}(y_n)\to1/Z.
\]
Thus $|x|^\gamma J_{\rm osc}(x)$ has two distinct subsequential limits. Theorem~\ref{thm:main}, Proposition~\ref{thm:smooth-error}, Corollary~\ref{cor:moments}, and \eqref{eq:arrival} apply to $J_{\rm osc}$ in \eqref{eq:oscillatory-kernel}.

For reactions, a simple family is
\begin{equation}\label{eq:reaction-example}
 f(s)=rs(1-s^a),\qquad 0\le s\le1,
 \qquad r>0,\quad a>0.
\end{equation}
This function is $C^1$ on $[0,1]$, satisfies $f'(0)=r$, $f'(1)=-ra<0$, and has nonlinear remainder $rs-f(s)=rs^{1+a}$. It obeys \eqref{eq:F} with $\sigma=\min\{a,1\}$ and $K=r$. Hence Theorem~\ref{thm:main} and Proposition~\ref{prop:interior-rate} both apply. This family includes reactions beyond the quadratic logistic case.

For \eqref{eq:reaction-example}, differentiation gives $f'(s)=r[1-(a+1)s^a]$ for $s>0$, which tends to $r$ at zero. Since $a\ge\sigma=\min\{a,1\}$ and $0\le s\le1$, one has $0\le rs-f(s)=rs^{1+a}\le rs^{1+\sigma}$, verifying \eqref{eq:F}.

A $C^1$ KPP reaction can have a remainder larger than every power $s^{1+\sigma}$. For example, the function defined for $0<s<\ee^{-2}$ by
\begin{equation}\label{eq:slow-remainder-example}
 f(s)=rs\left(1-\frac1{\log(1/s)}\right),\qquad f(0)=0,
\end{equation}
satisfies $0<f(s)\le rs$ for $0<s<\ee^{-2}$, and has derivative $r$ at zero with $f'$ continuous there.

Writing the logarithm explicitly, differentiation for $s>0$ gives
\[
 f'(s)=r\left[1-\frac1{\log(1/s)}-\frac1{\log(1/s)^2}\right]\longrightarrow r
 \text{ as }s\downarrow0.
\]
The equality $f'(0)=r$ follows from $f(s)/s\to r$. But
\[
 \frac{rs-f(s)}{s^{1+\sigma}}
 =\frac{r}{s^\sigma\log(1/s)}\to\infty
 \qquad\text{as }s\downarrow0
\]
for every $\sigma>0$. To define a reaction on $[0,1]$ that agrees with \eqref{eq:slow-remainder-example} near zero, choose $\chi\in C^\infty([0,1])$ with $0\le\chi\le1$, $\chi=1$ on $[0,\ee^{-3}]$, and $\chi=0$ on $[\ee^{-2},1]$. For $0<s<\ee^{-2}$ set $q(s)=\chi(s)/\log(1/s)$, and set $q(s)=0$ for $s\ge\ee^{-2}$. Define $\widetilde f(0)=0$ and $\widetilde f(s)=rs[1-q(s)-(1-\chi(s))s]$ for $0<s\le1$. For $0<s<\ee^{-2}$, the identity $1-q(s)-(1-\chi(s))s=\chi(s)[1-1/\log(1/s)]+(1-\chi(s))(1-s)$ shows that this factor lies strictly between zero and one. Thus $\widetilde f\in C^1([0,1])$, $\widetilde f(0)=\widetilde f(1)=0$, and $0<\widetilde f(s)\le rs$ for $0<s<1$. Moreover, $\widetilde f$ agrees with \eqref{eq:slow-remainder-example} near zero and equals $rs(1-s)$ for $s\ge\ee^{-2}$. The displayed divergence shows that $\widetilde f$ violates the power remainder bound in \eqref{eq:F} for every $\sigma>0$.

\paragraph{Spatial rescaling.}

The exponent $r/\gamma$ is independent of $d$, but the constants may depend on $d$, $J$, and $u_0$. The following rescaling records their dependence on the spatial unit.

\begin{proposition}\label{prop:spatial-scaling}
For the solution $u$ of \eqref{eq:cauchy}, fix $\ell>0$ and define
\begin{equation}\label{eq:spatial-scaling}
 J_\ell(x)=\ell^{-1}J(x/\ell),\qquad
 u_{0,\ell}(x)=u_0(x/\ell),\qquad
 u_\ell(t,x)=u(t,x/\ell).
\end{equation}
For $0<\theta<1$, write $E_{\theta,\ell}(t)=\{x\in\R:u_\ell(t,x)=\theta\}$ and $\mathcal M_\ell(t)=\int_{\R}u_\ell(t,x)\dd x$. Then $u_\ell$ solves \eqref{eq:cauchy} with kernel $J_\ell$ and initial function $u_{0,\ell}$, with the same $d$ and $f$. The tail exponent is still $\gamma$, and
\begin{equation}\label{eq:scaled-levels}
 E_{\theta,\ell}(t)=\ell E_\theta(t),\qquad
 \mathcal M_\ell(t)=\ell\mathcal M(t).
\end{equation}
If \eqref{eq:exact-tail} holds, define $\mathcal A_{\infty,\ell}=d\int_0^\infty e^{-rt}\mathcal M_\ell(t)\dd t$. Then $J_\ell$ has tail constant
$\lambda_\ell=\ell^{\gamma-1}\lambda$ and
$\mathcal A_{\infty,\ell}=\ell\mathcal A_\infty$.
\end{proposition}

\begin{proof}
The change of variables $y=\ell z$ gives
\[
 \int_{\R}J_\ell(x-y)u_\ell(t,y)\dd y
 =\int_{\R}J(x/\ell-z)u(t,z)\dd z
 =(J*u(t,\cdot))(x/\ell).
\]
Differentiating $u_\ell$ in time proves the equation; the initial condition follows from the definition of $u_\ell$.

The convolution substitution is valid since $0\le u_\ell\le1$ and $\int_\R J_\ell=1$, so its absolute integral is finite. The equality
$\partial_tu_\ell(t,x)=u_t(t,x/\ell)$ then gives the differential equation pointwise. Moreover, for every $x$,
\[
 x\in E_{\theta,\ell}(t)\ \Longleftrightarrow\ x/\ell\in E_\theta(t),\qquad
 \int_\R u_\ell(t,x)\dd x=\ell\int_\R u(t,y)\dd y.
\]
The second equality uses $x=\ell y$ and the integrability of $u(t,\cdot)$ from Lemma~\ref{lem:weight}. This proves \eqref{eq:scaled-levels}. The tail limit of $J_\ell$ follows from
$|x|^\gamma J_\ell(x)=\ell^{\gamma-1}|x/\ell|^\gamma J(x/\ell)$. The formula for $\mathcal A_{\infty,\ell}$ follows from its defining time integral and the mass identity in \eqref{eq:scaled-levels}.

Normalization follows from $\int_\R\ell^{-1}J(x/\ell)\dd x=1$; evenness, nonnegativity, boundedness, continuity and $J_\ell(0)>0$ follow directly from \eqref{eq:spatial-scaling}. The two bounds in \eqref{eq:tail-explicit} become
\[
 j_1\ell^{\gamma-1}|x|^{-\gamma}\le J_\ell(x)
 \le j_2\ell^{\gamma-1}|x|^{-\gamma}
 \text{ for }|x|\ge\ell R_J.
\]
Under \eqref{eq:exact-tail}, pulling the positive constant $\ell$ through the finite defining time integral gives
$\mathcal A_{\infty,\ell}=d\int_0^\infty\ee^{-rt}\ell\mathcal M(t)\dd t=\ell\mathcal A_\infty$.
\end{proof}

By \eqref{eq:scaled-levels}, rescaling multiplies every level-set position by $\ell$ and preserves the exponential rate $r/\gamma$.

\paragraph{Half-line data.} The distinction between compact support and an occupied half-line already appears in one convolution. Under \eqref{eq:exact-tail}, for the compact datum in \eqref{eq:data},
\begin{equation}\label{eq:compact-one-convolution}
 \lim_{x\to+\infty}\frac{(J*u_0)(x)}{J(x)}=M_0.
\end{equation}
Indeed, \eqref{eq:exact-tail} gives $J(x-y)/J(x)\to1$ uniformly for $y\in\supp u_0$ as $x\to+\infty$. Integrating this ratio against $u_0(y)\dd y$ on its compact support proves \eqref{eq:compact-one-convolution}.

For justification, fix $R_u$ with $\supp u_0\subset[-R_u,R_u]$ and factor
\[
 \frac{J(x-y)}{J(x)}=
 \frac{J(x-y)}{h_\gamma(x-y)}
 \frac{h_\gamma(x-y)}{h_\gamma(x)}
 \frac{h_\gamma(x)}{J(x)}.
\]
The first and third factors tend uniformly to $\lambda$ and $1/\lambda$, respectively, and the middle factor tends uniformly to one by \eqref{eq:h-shift-limit}. The product is therefore uniformly bounded on this support for large $x$. A constant multiple of $u_0(y)\ind_{[-R_u,R_u]}(y)$ is an integrable dominating function, so dominated convergence proves \eqref{eq:compact-one-convolution}. For the half-line datum $\phi_0=\ind_{(-\infty,0]}$, however,
\begin{equation}\label{eq:halfline-one-convolution}
 (J*\phi_0)(x)=\int_x^\infty J(z)\dd z,
 \qquad
 \lim_{x\to+\infty}x^{\gamma-1}(J*\phi_0)(x)
 =\frac{\lambda}{\gamma-1}.
\end{equation}
The limit follows by integrating the inequalities
$(\lambda-\varepsilon)z^{-\gamma}\le J(z)\le(\lambda+\varepsilon)z^{-\gamma}$ for large $z$ and then letting $\varepsilon\downarrow0$.

Equation \eqref{eq:halfline-one-convolution} gives $(J*\phi_0)(x)\sim x^{1-\gamma}$ as $x\to+\infty$, whereas \eqref{eq:compact-one-convolution} gives $(J*u_0)(x)/J(x)\to M_0$. The distinction between the two initial geometries is treated in \cite{FinkelshteinTkachov2019}.

\subsection{Comparison with free boundaries}

Let $g_{\rm fb}(t)$ and $h_{\rm fb}(t)$ denote the left and right endpoints in the spreading solutions studied in \cite{DuNiJEMS}, and let $u_{\rm fb}$ denote the corresponding density. In the finite-speed cases, $c_0$ is the speed determined by the semi-wave and its flux condition.

The free-boundary model in \cite{DuNiJEMS} prescribes, in addition to the equation on $(g_{\rm fb}(t),h_{\rm fb}(t))$, the right-boundary flux
\begin{equation}\label{eq:fb-flux-comparison}
 h_{\rm fb}'(t)=\mu\int_{g_{\rm fb}(t)}^{h_{\rm fb}(t)}
 u_{\rm fb}(t,x)\int_{h_{\rm fb}(t)}^\infty J(x-y)\dd y\dd x,
 \qquad \mu>0.
\end{equation}
The density $u_{\rm fb}$ vanishes at both boundaries and outside this interval.

The middle column below records \cite[Theorems~1.1--1.3]{DuNiJEMS}, always in the spreading case, with the kernel hypothesis \eqref{eq:J} and the initial data specified in that paper. Their basic reaction hypothesis is $f\in C^1([0,\infty))$, $f(0)=f(1)=0$, $f(s)>0$ for $0<s<1$, $f'(0)>0>f'(1)$, and $f(s)/s$ nonincreasing for $s>0$. The basic reaction conditions listed above are those of \cite[Theorem~1.2]{DuNiJEMS}. Theorems~1.1 and~1.3 of \cite{DuNiJEMS} additionally assume that $f'$ is locally Lipschitz on $[0,\infty)$; Theorem~1.3 of \cite{DuNiJEMS} further requires $[f(s)/s]'<0$ for $s>0$. The first column identifies the applicable theorem in \cite{DuNiJEMS}; the last column follows from Theorem~\ref{thm:main}\textup{(ii)}.
\begin{center}
\small
\begin{tabularx}{\textwidth}{@{}lXX@{}}
\toprule
Range / JEMS theorem & Free-boundary quantity & Cauchy level sets\\
\midrule
$1<\gamma<2$ / 1.2 & $h_{\rm fb}(t)\sim t^{1/(\gamma-1)}$ & $|x_\theta^\pm(t)|\sim\ee^{rt/\gamma}$\\[1mm]
$\gamma=2$ / 1.2 & $h_{\rm fb}(t)\sim t\log t$ & $|x_\theta^\pm(t)|\sim\ee^{rt/2}$\\[1mm]
$2<\gamma<3$ / 1.3 & $c_0t-h_{\rm fb}(t)\sim t^{3-\gamma}$ & $|x_\theta^\pm(t)|\sim\ee^{rt/\gamma}$\\[1mm]
$\gamma=3$ / 1.3 & $c_0t-h_{\rm fb}(t)\sim\log t$ & $|x_\theta^\pm(t)|\sim\ee^{rt/3}$\\[1mm]
$\gamma>3$ / 1.1 & $h_{\rm fb}(t)-c_0t=O(1)$ & $|x_\theta^\pm(t)|\sim\ee^{rt/\gamma}$\\
\bottomrule
\end{tabularx}
\end{center}
The notation $\sim$ in both columns means comparison by positive constants, as defined after \eqref{eq:tail}. In the last row, $O(1)$ means that $h_{\rm fb}(t)-c_0t$ remains bounded as $t\to\infty$.

Under an exact kernel tail limit, Du and Ni \cite{DuNiMathAnn} determine the free-boundary acceleration constants. Theorem~\ref{thm:tail-coefficient} gives the fixed-time spatial coefficient for \eqref{eq:cauchy}.


\begin{thebibliography}{99}

\bibitem{AlfaroCoville2017}
M. Alfaro and J. Coville,
Propagation phenomena in monostable integro-differential equations: acceleration or not?,
\emph{J. Differential Equations} \textbf{263} (2017), 5727--5758.

\bibitem{Andreu2010}
F. Andreu-Vaillo, J. M. Maz\'on, J. D. Rossi and J. J. Toledo-Melero,
\emph{Nonlocal Diffusion Problems}, Mathematical Surveys and Monographs, vol.~165,
American Mathematical Society, Providence, RI, and Real Sociedad Matem\'atica Espa\~nola, Madrid, 2010.

\bibitem{Applebaum2009}
D. Applebaum,
\emph{L\'evy Processes and Stochastic Calculus}, 2nd ed.,
Cambridge Studies in Advanced Mathematics, vol.~116,
Cambridge University Press, Cambridge, 2009.

\bibitem{AronsonWeinberger1975}
D. G. Aronson and H. F. Weinberger,
Nonlinear diffusion in population genetics, combustion, and nerve pulse propagation,
in \emph{Partial Differential Equations and Related Topics}, Lecture Notes in Mathematics, vol.~446,
Springer, Berlin, 1975, pp.~5--49.

\bibitem{AronsonWeinberger1978}
D. G. Aronson and H. F. Weinberger,
Multidimensional nonlinear diffusion arising in population genetics,
\emph{Adv. Math.} \textbf{30} (1978), 33--76.

\bibitem{BouinEtAl2018}
E. Bouin, J. Garnier, C. Henderson and F. Patout,
Thin front limit of an integro-differential Fisher--KPP equation with fat-tailed kernels,
\emph{SIAM J. Math. Anal.} \textbf{50} (2018), 3365--3394.

\bibitem{Boutillon2024}
N. Boutillon,
Large deviations and the emergence of a logarithmic delay in a nonlocal linearised Fisher--KPP equation,
\emph{Nonlinear Anal.} \textbf{240} (2024), 113465.

\bibitem{Bramson1978}
M. D. Bramson,
Maximal displacement of branching Brownian motion,
\emph{Comm. Pure Appl. Math.} \textbf{31} (1978), 531--581.

\bibitem{Bramson1983}
M. Bramson,
Convergence of solutions of the Kolmogorov equation to travelling waves,
\emph{Mem. Amer. Math. Soc.} \textbf{44} (1983), no.~285, iv+190 pp.

\bibitem{CabreCoulonRoquejoffre2012}
X. Cabr\'e, A.-C. Coulon and J.-M. Roquejoffre,
Propagation in Fisher--KPP type equations with fractional diffusion in periodic media,
\emph{C. R. Math. Acad. Sci. Paris} \textbf{350} (2012), 885--890.

\bibitem{CabreRoquejoffre2009}
X. Cabr\'e and J.-M. Roquejoffre,
Front propagation in Fisher--KPP equations with fractional diffusion,
\emph{C. R. Math. Acad. Sci. Paris} \textbf{347} (2009), 1361--1366.

\bibitem{CabreRoquejoffre2013}
X. Cabr\'e and J.-M. Roquejoffre,
The influence of fractional diffusion in Fisher--KPP equations,
\emph{Comm. Math. Phys.} \textbf{320} (2013), 679--722.

\bibitem{CaoDuLiLi2019}
J.-F. Cao, Y. Du, F. Li and W.-T. Li,
The dynamics of a Fisher--KPP nonlocal diffusion model with free boundaries,
\emph{J. Funct. Anal.} \textbf{277} (2019), 2772--2814.

\bibitem{CovilleDavilaMartinez2008}
J. Coville, J. D\'avila and S. Mart\'inez,
Nonlocal anisotropic dispersal with monostable nonlinearity,
\emph{J. Differential Equations} \textbf{244} (2008), 3080--3118.

\bibitem{CuiLiXu2025}
T.-L. Cui, W.-T. Li and W.-B. Xu,
Accelerating or not in the spatial propagation of nonlocal dispersal cooperative reducible systems,
\emph{J. Differential Equations} \textbf{443} (2025), 113519.

\bibitem{DuLiZhou2021}
Y. Du, F. Li and M. Zhou,
Semi-wave and spreading speed of the nonlocal Fisher--KPP equation with free boundaries,
\emph{J. Math. Pures Appl.} \textbf{154} (2021), 30--66.

\bibitem{DuNiJEMS}
Y. Du and W. Ni,
Rate of propagation for the Fisher--KPP equation with nonlocal diffusion and free boundaries,
\emph{J. Eur. Math. Soc.} \textbf{27} (2025), 1267--1319.

\bibitem{DuNiMathAnn}
Y. Du and W. Ni,
Exact rate of accelerated propagation in the Fisher--KPP equation with nonlocal diffusion and free boundaries,
\emph{Math. Ann.} \textbf{389} (2024), 2931--2958.

\bibitem{Durrett2019}
R. Durrett,
\emph{Probability: Theory and Examples}, 5th ed.,
Cambridge University Press, Cambridge, 2019.

\bibitem{FinkelshteinEtAl2019}
D. Finkelshtein, Y. Kondratiev and P. Tkachov,
Accelerated front propagation for monostable equations with nonlocal diffusion: multidimensional case,
\emph{J. Elliptic Parabol. Equ.} \textbf{5} (2019), 423--471.

\bibitem{FinkelshteinTkachov2019}
D. Finkelshtein and P. Tkachov,
Accelerated nonlocal nonsymmetric dispersion for monostable equations on the real line,
\emph{Appl. Anal.} \textbf{98} (2019), 756--780.

\bibitem{FinkelshteinTkachovKesten2018}
D. Finkelshtein and P. Tkachov,
Kesten's bound for subexponential densities on the real line and its multi-dimensional analogues,
\emph{Adv. Appl. Probab.} \textbf{50} (2018), 373--395.

\bibitem{Fisher1937}
R. A. Fisher,
The wave of advance of advantageous genes,
\emph{Ann. Eugenics} \textbf{7} (1937), 355--369.

\bibitem{FossKorshunovZachary2013}
S. Foss, D. Korshunov and S. Zachary,
\emph{An Introduction to Heavy-Tailed and Subexponential Distributions}, 2nd ed.,
Springer Series in Operations Research and Financial Engineering,
Springer, New York, 2013.

\bibitem{Garnier2011}
J. Garnier,
Accelerating solutions in integro-differential equations,
\emph{SIAM J. Math. Anal.} \textbf{43} (2011), 1955--1974.

\bibitem{Graham2022}
C. Graham,
The Bramson correction for integro-differential Fisher--KPP equations,
\emph{Commun. Math. Sci.} \textbf{20} (2022), 563--596.

\bibitem{HamelRoques2010}
F. Hamel and L. Roques,
Fast propagation for KPP equations with slowly decaying initial conditions,
\emph{J. Differential Equations} \textbf{249} (2010), 1726--1745.

\bibitem{JinZhou2026}
Z. Jin and T. Zhou,
Spreading speeds for nonlocal Fisher--KPP equations with time-dependent asymmetric kernels,
preprint, arXiv:2606.15226, 2026.

\bibitem{KotLewisDriessche1996}
M. Kot, M. A. Lewis and P. van den Driessche,
Dispersal data and the spread of invading organisms,
\emph{Ecology} \textbf{77} (1996), 2027--2042.

\bibitem{KPP1937}
A. N. Kolmogorov, I. G. Petrovsky and N. S. Piskunov,
\'Etude de l'\'equation de la diffusion avec croissance de la quantit\'e de mati\`ere
et son application \`a un probl\`eme biologique,
\emph{Bull. Univ. \'Etat Moscou S\'er. Internat. A} \textbf{1} (1937), 1--26.

\bibitem{LiAge2026}
N. Li,
Sharp rate of accelerating propagation in a nonlocal model arising in population dynamics,
\emph{Nonlinear Anal. Real World Appl.} \textbf{89} (2026), 104526.

\bibitem{LiCooperative2026}
N. Li,
Sharp rates of accelerating propagation for nonlocal dispersal cooperative systems,
\emph{J. Evol. Equ.} \textbf{26} (2026), article~5.

\bibitem{LiPanSu2025}
N. Li, Y. Pan and Y. Su,
Sharp rate of the accelerating propagation for a recursive system,
\emph{Stud. Appl. Math.} \textbf{154} (2025), e70029.

\bibitem{MedlockKot2003}
J. Medlock and M. Kot,
Spreading disease: integro-differential equations old and new,
\emph{Math. Biosci.} \textbf{184} (2003), 201--222.

\bibitem{MeleardMirrahimi2015}
S. M\'el\'eard and S. Mirrahimi,
Singular limits for reaction-diffusion equations with fractional Laplacian and local or nonlocal nonlinearity,
\emph{Comm. Partial Differential Equations} \textbf{40} (2015), 957--993.

\bibitem{Roquejoffre2024}
J.-M. Roquejoffre,
\emph{The Dynamics of Front Propagation in Nonlocal Reaction--Diffusion Equations},
Lecture Notes on Mathematical Modelling in the Life Sciences,
Springer, Cham, 2024.

\bibitem{RoquejoffreTarfulea2017}
J.-M. Roquejoffre and A. Tarfulea,
Gradient estimates and symmetrization for Fisher--KPP front propagation with fractional diffusion,
\emph{J. Math. Pures Appl.} \textbf{108} (2017), 399--424.

\bibitem{RoquejoffreZhang2026}
J.-M. Roquejoffre and M. Zhang,
Large time behaviour in the multidimensional nonlocal Fisher--KPP equation,
preprint, arXiv:2609.32051v1, 2026.

\bibitem{Schumacher1980}
K. Schumacher,
Travelling-front solutions for integro-differential equations. I,
\emph{J. Reine Angew. Math.} \textbf{316} (1980), 54--70.

\bibitem{SouganidisTarfulea2019}
P. E. Souganidis and A. Tarfulea,
Front propagation for integro-differential KPP reaction--diffusion equations in periodic media,
\emph{NoDEA Nonlinear Differential Equations Appl.} \textbf{26} (2019), article~29.

\bibitem{TaoEtAl2025}
W. Tao, W.-T. Li, J.-W. Sun and W.-B. Xu,
Interface dynamics in nonlocal dispersal Fisher--KPP equations,
\emph{Acta Math. Sci. Ser. B} \textbf{45} (2025), 1774--1813.

\bibitem{Weinberger1982}
H. F. Weinberger,
Long-time behavior of a class of biological models,
\emph{SIAM J. Math. Anal.} \textbf{13} (1982), 353--396.

\bibitem{WeinbergerZhao2010}
H. F. Weinberger and X.-Q. Zhao,
An extension of the formula for spreading speeds,
\emph{Math. Biosci. Eng.} \textbf{7} (2010), 187--194.

\bibitem{XuLiRuan2021}
W.-B. Xu, W.-T. Li and S. Ruan,
Spatial propagation in nonlocal dispersal Fisher--KPP equations,
\emph{J. Funct. Anal.} \textbf{280} (2021), 108957.

\bibitem{Yagisita2009}
H. Yagisita,
Existence and nonexistence of traveling waves for a nonlocal monostable equation,
\emph{Publ. Res. Inst. Math. Sci.} \textbf{45} (2009), 925--953.

\end{thebibliography}
\end{document}